\documentclass[11pt,reqno]{amsart}

\usepackage[margin=0.8in]{geometry}
\usepackage{amsmath,amssymb,amsthm,mathtools}
\usepackage{mathrsfs}
\usepackage{enumitem}
\usepackage{booktabs}
\usepackage{xcolor}
\usepackage[colorlinks=true,linkcolor=blue!50!black,citecolor=green!40!black,urlcolor=blue!50!black]{hyperref}

\numberwithin{equation}{section}

\newtheorem{theorem}{Theorem}[section]
\newtheorem{proposition}[theorem]{Proposition}
\newtheorem{lemma}[theorem]{Lemma}
\newtheorem{corollary}[theorem]{Corollary}
\theoremstyle{definition}
\newtheorem{definition}[theorem]{Definition}
\newtheorem{assumption}[theorem]{Assumption}
\theoremstyle{remark}
\newtheorem{remark}[theorem]{Remark}

\newcommand{\eps}{\varepsilon}
\newcommand{\R}{\mathbb{R}}
\newcommand{\N}{\mathbb{N}}
\newcommand{\E}{\mathbb{E}}
\renewcommand{\P}{\mathbb{P}}
\newcommand{\Ocal}{\mathcal{O}}
\newcommand{\LL}{\mathbb{L}^2}
\newcommand{\HH}{\mathbb{H}^1}
\newcommand{\HHH}{\mathbb{H}^2}
\newcommand{\Hz}{H_0}
\newcommand{\Ast}{\mathcal{A}}
\newcommand{\An}{\mathrm{A}}
\newcommand{\DA}{D(\mathrm{A})}
\newcommand{\Bt}{\widetilde{B}}
\newcommand{\KT}{\mathcal{K}_T}
\newcommand{\LHS}{\mathscr{L}_2}
\newcommand{\LHSd}{\mathscr{L}_2(\ell^2,\mathbb{H}^1)}
\newcommand{\Xs}{\mathfrak{X}}
\newcommand{\Fc}{\mathcal{F}}
\newcommand{\Ec}{\mathfrak{E}}
\newcommand{\Dc}{\mathfrak{D}}
\newcommand{\Lin}{\mathfrak{L}}
\newcommand{\Gam}{\Gamma}
\newcommand{\Sb}{\mathbb{S}}
\newcommand{\Ac}{\mathcal{A}^{\mathrm{c}}}
\newcommand{\dd}{\,\mathrm{d}}

\newcommand{\Law}{\operatorname{Law}}
\newcommand{\abs}[1]{\lvert #1\rvert}
\newcommand{\norm}[1]{\lVert #1\rVert}
\newcommand{\hone}[1]{\lVert #1\rVert_{\HH}}
\newcommand{\hda}[1]{\lvert #1\rvert_{\DA}}
\newcommand{\dual}[2]{\langle #1,#2\rangle}
\newcommand{\scal}[2]{( #1,#2)}
\newcommand{\sone}[2]{( #1,#2)_{\HH}}

\newcommand{\ft}{\tilde f}
\newcommand{\Ft}{\widetilde F}

\begin{document}

\title[CLT and moderate deviations for stochastic nematic flows]
{Gaussian fluctuations and moderate deviations for\\
two-dimensional stochastic nematic liquid crystal flows}

\author{Akash Ashirbad Panda$^{\dagger}$}
\author{Utpal Manna}

\thanks{%
$^{\dagger}$Corresponding author: Akash Ashirbad Panda,
Department of Mathematics, Indian Institute of Technology Indore,
Simrol, Indore 453552, Madhya Pradesh, India.
Email: \texttt{akashpanda@iiti.ac.in}.}

\thanks{%
Utpal Manna,
School of Mathematics, Indian Institute of Science Education and Research
Thiruvananthapuram,
Vithura, Thiruvananthapuram 695551, Kerala, India.
Email: \texttt{manna.utpal@iisertvm.ac.in}.}

\date{\today}

\begin{abstract}
We study the small-noise behaviour of a two-dimensional stochastic nematic
liquid crystal model, where the Navier--Stokes equations with multiplicative
noise are coupled with a Ginzburg--Landau equation for the director field
driven by geometric Stratonovich noise. While large deviations for stochastic
nematic liquid crystal models have been studied earlier, the central limit
and moderate deviation regimes have not been considered. We first prove a
central limit theorem in the natural energy space and show that the rescaled
fluctuations converge, in every $L^p$-norm and with rate $\sqrt{\eps}$, to the
solution of a linear stochastic system. The limiting process is centred
Gaussian. We then establish a moderate deviation principle with a quadratic
rate function determined by the corresponding linearised skeleton equation,
which coincides with the Cameron--Martin form of the Gaussian limit. We also
show that the moderate deviation principle holds uniformly over bounded sets
of initial data. The main arguments use a cancellation between the coupling
terms, a one-sided estimate for the polynomial potential, and energy estimates
with deterministic Gronwall coefficients, thereby avoiding exponential moment
estimates.
\end{abstract}

\maketitle


\smallskip


 \tableofcontents

\section{Introduction}\label{sec:intro}

\subsection{The model}

Nematic liquid crystals are anisotropic fluids whose molecules have orientational order but no positional order. Their dynamics are described by the fluid velocity $u$ and the \emph{director} $d$, which represents the preferred molecular orientation. The continuum theory was developed by Ericksen \cite{Ericksen61} and Leslie \cite{Leslie68}. Later, Lin and Liu \cite{LinLiu95} introduced a simplified model which retains the main coupling between the velocity and the director field and replaces the constraint $\abs d=1$ by a Ginzburg--Landau penalisation.

Following \cite{BHR19a,BMP20}, we consider a general radial polynomial
$f(d)=\ft(\abs d^2)d$ and study the following two-dimensional stochastic system:
\begin{align}
  \dd u^\eps + \big[(u^\eps\cdot\nabla)u^\eps-\Delta u^\eps+\nabla p^\eps\big]\dd t
  &= -\nabla\cdot\big(\nabla d^\eps\odot\nabla d^\eps\big)\dd t
  +\sqrt\eps\,\sigma(u^\eps)\dd W_1, \label{eq:intro-u}\\
  \nabla\cdot u^\eps &=0, \label{eq:intro-div}\\
  \dd d^\eps+(u^\eps\cdot\nabla)d^\eps\dd t
  &=\big[\Delta d^\eps-f(d^\eps)\big]\dd t
  +\sqrt\eps\sum_{k\ge1}(d^\eps\times h_k)\circ\dd W_2^k,
  \label{eq:intro-d}
\end{align}
in $(0,T]\times\Ocal$, where $\Ocal\subset\R^2$ is a bounded domain with smooth boundary. We impose
\begin{equation}\label{eq:intro-bc}
  u^\eps=0,\qquad
  \frac{\partial d^\eps}{\partial n}=0
  \quad\text{on }\partial\Ocal,\qquad
  u^\eps(0)=u_0,\quad d^\eps(0)=d_0 .
\end{equation}
Here $u^\eps:\Ocal\to\R^2$, $d^\eps:\Ocal\to\R^3$, and $p^\eps$ denotes the pressure. Moreover,
\[
[\nabla d\odot\nabla d]_{ij}
=\sum_{k=1}^3\partial_{x_i}d^{(k)}\partial_{x_j}d^{(k)}
\]
is the Ericksen stress. The process $W_1$ is a Wiener process acting on the velocity, while $\{W_2^k\}_{k\ge1}$ are independent real Brownian motions, independent of $W_1$. We assume
\[
\sum_{k\ge1}
\left(
\norm{h_k}_{L^\infty}
+\norm{\nabla h_k}_{L^3}
\right)^2<\infty.
\]
All physical constants are normalised to one.

The velocity equation is perturbed by multiplicative It\^o noise, whereas the director equation contains geometric Stratonovich noise. Since
\[
d\cdot(d\times h_k)=0,
\]
the noise $d\times h_k$ acts tangentially to the sphere and, by itself, preserves the length of the director. This form of noise is also natural when considering the limit in which the Ginzburg--Landau penalisation is removed and the constraint $\abs d=1$ is recovered; see \cite{BHR19a,BMP20,BDR20}.

For a single director noise, well-posedness was studied in \cite{BHR19a}, while Brze\'zniak, Manna and Panda \cite{BMP20} established a Freidlin--Wentzell large deviation principle (LDP) \cite{FW12}. In the present paper we allow countably many geometric noises and study the fluctuations of the solution in the small-noise limit.

\subsection{Small-noise regimes and main results}

Let
\[
X^\eps=(u^\eps,d^\eps),\qquad X^0=(u^0,d^0),
\]
where $X^0$ denotes the solution of the deterministic system obtained by setting $\eps=0$. Different asymptotic regimes arise depending on the scale of $X^\eps-X^0$.

\begin{table}[ht]
\centering
\renewcommand{\arraystretch}{1.2}
\begin{tabular}{@{}llll@{}}
\toprule
Regime & Rescaled process & Main question & Result \\ \midrule
Law of large numbers
& $X^\eps-X^0$
& convergence and rate
& Theorem~\ref{thm:LLN}
\\
Central limit theorem
& $(X^\eps-X^0)/\sqrt\eps$
& Gaussian limit
& Theorems~\ref{thm:linear}, \ref{thm:CLT}
\\
Moderate deviations
& $(X^\eps-X^0)/(\sqrt\eps\,a(\eps))$
& LDP with speed $a(\eps)^2$
& Theorems~\ref{thm:MDP}, \ref{thm:uniform}
\\
Large deviations
& $X^\eps$
& LDP with speed $\eps^{-1}$
& \cite{BMP20}
\\ \bottomrule
\end{tabular}
\caption{Small-noise regimes for \eqref{eq:intro-u}--\eqref{eq:intro-bc}. In the moderate deviation regime,
$a(\eps)\to\infty$ and $\sqrt\eps\,a(\eps)\to0$.}
\label{tab:regimes}
\end{table}

The LDP of \cite{BMP20} describes deviations of order one and is governed by a nonlinear control problem. The central limit theorem (CLT) describes the usual fluctuations of order $\sqrt\eps$. Between these two scales lies the moderate deviation regime
\[
\sqrt\eps\,a(\eps),\qquad
a(\eps)\to\infty,\qquad
\sqrt\eps\,a(\eps)\to0,
\]
whose rate function is determined by the linearisation of the system around $X^0$.

Large deviations for stochastic liquid crystal models have been considered in
\cite{ZhangZhou18,BMP20,Medjo23,Medjo24,PSZZ22}. To the best of our knowledge, the CLT and moderate deviation principle (MDP) have not previously been established for the stochastic nematic system considered here.

We work in the energy space
\[
\KT=
\big[C([0,T];H)\cap L^2(0,T;V)\big]
\times
\big[C([0,T];\HH)\cap L^2(0,T;\DA)\big].
\]
Our first result gives the quantitative law of large numbers
\[
\E\norm{X^\eps-X^0}_{\KT}^{2p}\le C_p\eps^p,
\qquad p\ge1.
\]
We then prove that
\[
\frac{X^\eps-X^0}{\sqrt\eps}
\longrightarrow X^1
\quad\text{in }L^{2p}(\Omega;\KT),
\]
where $X^1$ is the solution of the linearised stochastic system along $X^0$. More precisely,
\[
\E\Big\|
\frac{X^\eps-X^0}{\sqrt\eps}-X^1
\Big\|_{\KT}^{2p}
\le C_p\eps^p.
\]
The limit $X^1$ is a centred Gaussian random variable in $\KT$.

Let $\Gam(q)$ denote the solution of the corresponding deterministic linearised system with control
\[
q=(\theta,\rho)\in L^2(0,T;\Hz\times\ell^2).
\]
We identify the covariance and Cameron--Martin space of $X^1$ through $\Gam$. In particular,
\[
\operatorname{Var}\ell(X^1)=\norm{\ell\circ\Gam}^2,
\qquad \ell\in\KT^*.
\]

For every $a(\eps)\to\infty$ such that $\sqrt\eps\,a(\eps)\to0$, we further prove that
\[
\frac{X^\eps-X^0}{\sqrt\eps\,a(\eps)}
\]
satisfies an LDP in $\KT$ with speed $a(\eps)^2$ and rate function
\[
I(z)=
\inf\left\{
\frac12\int_0^T
\big(
\abs{\theta(t)}_{\Hz}^2+\abs{\rho(t)}_{\ell^2}^2
\big)\dd t:
\Gam(\theta,\rho)=z
\right\}.
\]
Thus the MDP rate function is precisely the Cameron--Martin quadratic form of the Gaussian limit. We also prove that these results are uniform over bounded sets of initial data. In particular, the MDP holds uniformly in the Freidlin--Wentzell sense and as an equicontinuous uniform Laplace principle; over compact sets it also holds in the Dembo--Zeitouni sense.

\subsection{Relation to the literature and main difficulties}

The stochastic Ginzburg--Landau approximation of nematic liquid crystals has been studied in several works. Brze\'zniak, Hausenblas and Razafimandimby established existence of martingale solutions, pathwise uniqueness in two dimensions and a maximum principle in
\cite{BHR14,BHR19a}; strong solutions were further studied in \cite{BHR21}. For the constrained stochastic Ericksen--Leslie system, local smooth solutions were obtained in \cite{BHR19b}, while the Ginzburg--Landau approximation was studied in \cite{BDR20}; see also \cite{DHP19}. The corresponding system driven by pure jump noise in Marcus canonical form was studied in \cite{BMP19}.

Large deviation principles for liquid crystal models were proved in
\cite{ZhangZhou18,BMP20,Medjo23,Medjo24}. They can also be obtained within general frameworks for fully locally monotone SPDEs \cite{RSZ24,PSZZ22}. A recent framework for large and moderate deviations of hydrodynamic SPDEs with Poisson noise is given in \cite{GangulySundar26}.

For related hydrodynamic models, CLTs and MDPs have been studied for the stochastic Navier--Stokes equations \cite{SritharanSundar06,ChueshovMillet10,WZZ15}, primitive equations \cite{ZhangZhou17,Slavik22}, second grade fluids \cite{ZZZ16}, and hydrodynamic systems with jump noise \cite{LiShangWang25}; see also \cite{WangZhang15} for reaction--diffusion equations. Uniform large deviation principles over bounded sets were developed in \cite{SBD19,Salins19}, and a uniform MDP for a class of semilinear SPDEs was recently obtained in \cite{JS25}. For a recent uniform LDP for anisotropic stochastic Navier--Stokes equations, see \cite{SunYuLiu26}.

Our proof uses the weak convergence approach of Budhiraja and Dupuis
\cite{BD00,BDM08,BD19}, in the form of the criterion of Matoussi, Sabbagh and Zhang \cite{MSZ21}; for jump noise, see \cite{BDG16}.

The analysis of the nematic system has some additional difficulties. The Ericksen stress is quadratic in $\nabla d$, the polynomial $f$ is not globally Lipschitz, and the Stratonovich noise produces an It\^o correction that has to be controlled at each asymptotic scale. The main estimates rely on an exact cancellation between the Ericksen stress and the transport term in the director equation, together with a one-sided estimate for the polynomial potential. By expanding around the deterministic solution $X^0$, the resulting Gronwall coefficients are deterministic. This allows us to obtain the quantitative CLT and the uniform MDP without requiring exponential moment estimates.

We also allow countably many director noises. The corresponding well-posedness is proved in Appendix~\ref{app:wellposed} using the theory of fully locally monotone SPDEs \cite{RSZ24,RSZ26}.

\subsection{Organisation}

The paper is organised as follows. Section~\ref{sec:prelim} contains the functional setting and the main estimates used throughout the paper. In Section~\ref{sec:results}, we state the assumptions, introduce the linearised system and give the main results. Section~\ref{sec:LLN} proves the quantitative law of large numbers. The linearised system and the Gaussian limit are studied in Section~\ref{sec:linear}, and the CLT is proved in Section~\ref{sec:CLT}. Sections~\ref{sec:MDP} and \ref{sec:uniform} are devoted to the MDP and its uniform version. Some further remarks are given in Section~\ref{sec:concl}. The moment estimates and the well-posedness for countably many director noises are proved in the appendices.

\section{Preliminaries}\label{sec:prelim}

This section fixes the functional setting and collects the analytic tools used throughout the paper. Throughout, $\Ocal\subset\R^2$ is a bounded domain with $C^\infty$ boundary and $T>0$ is fixed.

\subsection{Function spaces and linear operators}\label{sec:spaces}
\emph{Lebesgue and Sobolev spaces.} For $p\in[1,\infty]$ and $k\in\N$ we write $L^p=L^p(\Ocal)$ and $H^k=W^{k,2}(\Ocal)$ for the Lebesgue and Sobolev spaces of scalar or vector-valued functions; the target space will be clear from the context. For $\R^3$-valued functions we write $\LL=L^2(\Ocal;\R^3)$, $\HH=H^1(\Ocal;\R^3)$ and $\HHH=H^2(\Ocal;\R^3)$.

\emph{Velocity.} Let $\mathscr V=\{v\in C^\infty_c(\Ocal;\R^2):\nabla\cdot v=0\}$, let $H$ be the closure of $\mathscr V$ in $L^2(\Ocal;\R^2)$ and $V$ its closure in $H^1(\Ocal;\R^2)$. On $H$ we use the $L^2$ inner product $\scal\cdot\cdot$ and norm $\abs\cdot$; on $V$ we use $\norm v=\abs{\nabla v}$, which is equivalent to the $H^1$ norm by the Poincar\'e inequality. Identifying $H$ with its dual we obtain the Gelfand triple
\[
  V\subset H\subset V',
\]
and we denote by $\dual\cdot\cdot$ the duality between $V'$ and $V$. The Stokes operator $\Ast\colon V\to V'$ is defined by $\dual{\Ast u}{v}=\scal{\nabla u}{\nabla v}$; see e.g.\ \cite{Temam,ConstantinFoias88} for these spaces and operators.

\emph{Director.} On $\LL$ we also write $\scal\cdot\cdot$ and $\abs\cdot$. On $\HH$ we use the norm $\hone d^2=\abs d^2+\abs{\nabla d}^2$ and the inner product $\sone{d}{e}=\scal de+\scal{\nabla d}{\nabla e}$. The Neumann Laplacian
\[
  \An d=-\Delta d,\qquad \DA=\Big\{d\in\HHH:\ \frac{\partial d}{\partial n}=0\ \text{on }\partial\Ocal\Big\},
\]
is a non-negative self-adjoint operator in $\LL$ with compact resolvent, and we equip $\DA$ with the graph norm $\hda d^2=\abs d^2+\abs{\An d}^2$. By elliptic regularity \cite{ADN59} there is $C_\Ocal$ such that
\begin{equation}\label{eq:elliptic}
  \norm{d}_{H^2}\le C_\Ocal\,\hda d,\qquad d\in\DA .
\end{equation}
For $g\in\HH$ and $e\in\DA$ Green's formula gives $\scal{g}{\An e}=\scal{\nabla g}{\nabla e}$, hence $\sone ge=\scal{g}{(I+\An)e}$. Consequently, if $\DA$ is normed by $\abs{(I+\An)\cdot}$, which is equivalent to $\hda\cdot$, we obtain a second Gelfand triple, with pivot space $\HH$:
\begin{equation}\label{eq:triple-d}
  \DA\subset\HH\subset\LL\cong(\DA)',\qquad {}_{(\DA)'}\langle g,e\rangle_{\DA}=\scal{g}{(I+\An)e}.
\end{equation}
All embeddings in both triples are dense and compact.

\emph{Summary of norms.} Thus $\abs\cdot$ and $\scal\cdot\cdot$ denote the norm and inner product of both $H$ and $\LL$, $\norm\cdot$ the norm of $V$, $\hone\cdot$ and $\sone\cdot\cdot$ the norm and inner product of $\HH$, and $\hda\cdot$ the graph norm of $\DA$. In all energy estimates, the energy of a pair $(v,e)$ is $\abs v^2+\hone e^2$ and its dissipation is $\norm v^2+\abs{\nabla e}^2+\abs{\An e}^2$.

\subsection{Two-dimensional inequalities}\label{sec:ineq}
\begin{lemma}\label{lem:ineq}
There is a constant $C$, and for every $q\in[1,\infty)$ a constant $C_q$, depending only on $\Ocal$, such that for all $v\in V$, $\varphi\in H^1(\Ocal)$ and $d\in\DA$,
\begin{align}
  \norm v_{L^4}&\le C\abs v^{1/2}\norm v^{1/2}, \label{eq:lady}\\
  \norm\varphi_{L^4}&\le C\abs\varphi^{1/2}\norm\varphi_{H^1}^{1/2},\qquad \norm\varphi_{L^q}\le C_q\norm\varphi_{H^1}, \label{eq:GN}\\
  \norm{\nabla d}_{L^4}&\le C\abs{\nabla d}^{1/2}\hda d^{1/2}. \label{eq:GNgrad}
\end{align}
\end{lemma}

\begin{proof}
\eqref{eq:lady} is Ladyzhenskaya's inequality \cite{Ladyzhenskaya69}, and \eqref{eq:GN} consists of the Gagliardo--Nirenberg inequality \cite{Nirenberg59} and the Sobolev embedding $H^1(\Ocal)\subset L^q(\Ocal)$, $q<\infty$, in dimension two; see e.g.\ \cite{Temam}. Applying the first inequality in \eqref{eq:GN} to $\varphi=\partial_{x_i}d$ and using $\norm{\partial_{x_i}d}_{H^1}\le\norm d_{H^2}$ and \eqref{eq:elliptic} gives \eqref{eq:GNgrad}.
\end{proof}

\subsection{The nonlinear terms and the energy cancellation}\label{sec:nonlinear}
Let
\[
  b(u,v,w)=\sum_{i,j=1}^2\int_\Ocal u^{(i)}\partial_{x_i}v^{(j)}w^{(j)}\dd x
\]
and define $B\colon V\times V\to V'$ by $\dual{B(u,v)}{w}=b(u,v,w)$; we write $B(u)=B(u,u)$. For $u\in V$ and $d\in\HH$ let $\Bt(u,d)=(u\cdot\nabla)d$. Finally, for $d_1,d_2\in\DA$ let $M(d_1,d_2)\in V'$ be defined by
\begin{equation}\label{eq:M}
  \dual{M(d_1,d_2)}{v}=-\sum_{i,j=1}^2\int_\Ocal\partial_{x_i}d_1\cdot\partial_{x_j}d_2\,\partial_{x_j}v^{(i)}\dd x,\qquad v\in V,
\end{equation}
and $M(d)=M(d,d)$. For smooth fields, $\dual{M(d)}{v}=\scal{\nabla\cdot(\nabla d\odot\nabla d)}{v}$, so that $M(d)$ is the Leray projection of the Ericksen stress term in \eqref{eq:intro-u}.

\begin{lemma}\label{lem:nonlinear}
There is a constant $C$, depending only on $\Ocal$, such that the following hold.
\begin{enumerate}[label=\textup{(\roman*)}]
\item $b(u,v,w)=-b(u,w,v)$ and $b(u,v,v)=0$ for $u,v,w\in V$.
\item For $u,v,w\in V$: $\abs{b(w,u,w)}\le C\abs w\norm w\norm u$ and $\abs{B(u,v)}_{V'}\le C\abs u^{1/2}\norm u^{1/2}\abs v^{1/2}\norm v^{1/2}$.
\item For $u\in V$ and $d\in\DA$: $\abs{\Bt(u,d)}\le C\abs u^{1/2}\norm u^{1/2}\abs{\nabla d}^{1/2}\hda d^{1/2}$. For $u\in V$ and $d,\varphi\in\HH$: $\abs{\scal{\Bt(u,d)}{\varphi}}\le C\norm u\abs{\nabla d}\hone\varphi$ and $\scal{\Bt(u,d)}{d}=0$.
\item For $d_1,d_2\in\DA$: $\abs{M(d_1,d_2)}_{V'}\le C\abs{\nabla d_1}^{1/2}\hda{d_1}^{1/2}\abs{\nabla d_2}^{1/2}\hda{d_2}^{1/2}$.
\end{enumerate}
\end{lemma}

\begin{proof}
(i) is classical \cite[Ch.~II]{Temam}. (ii) follows from (i), H\"older's inequality and \eqref{eq:lady}: $\abs{b(w,u,w)}\le\norm w_{L^4}^2\norm u$ and $\abs{b(u,v,w)}=\abs{b(u,w,v)}\le\norm u_{L^4}\norm w\norm v_{L^4}$. (iii) The first bound follows from $\abs{\Bt(u,d)}\le\norm u_{L^4}\norm{\nabla d}_{L^4}$, \eqref{eq:lady} and \eqref{eq:GNgrad}; the second from $\abs{\scal{\Bt(u,d)}{\varphi}}\le\norm u_{L^4}\abs{\nabla d}\norm\varphi_{L^4}$ and \eqref{eq:GN}. Finally $\scal{\Bt(u,d)}{d}=\frac12\int_\Ocal u\cdot\nabla\abs d^2\dd x=0$, since $\nabla\cdot u=0$ and $u|_{\partial\Ocal}=0$. (iv) By H\"older's inequality $\abs{\dual{M(d_1,d_2)}v}\le\norm{\nabla d_1}_{L^4}\norm{\nabla d_2}_{L^4}\norm v$, and we conclude with \eqref{eq:GNgrad}.
\end{proof}

The following identity expresses the exchange of energy between the flow and the director. It removes the cubic terms from all energy estimates in this paper.

\begin{lemma}[Cancellation]\label{lem:cancel}
For every $v\in V$ and $d\in\DA$,
\begin{equation}\label{eq:cancel}
  \dual{M(d,d)}{v}+\scal{\Bt(v,d)}{\An d}=0 .
\end{equation}
Consequently $\dual{M(d,d)}{v}+\scal{\Bt(v,d)}{(I+\An)d}=0$.
\end{lemma}

\begin{proof}
Let first $v\in\mathscr V$. Since $v$ vanishes near $\partial\Ocal$ and $d\in\HHH$, integration by parts in $x_j$ gives
\[
  \dual{M(d,d)}{v}=\sum_{i,j}\int_\Ocal\partial_{x_j}\big(\partial_{x_i}d\cdot\partial_{x_j}d\big)v^{(i)}\dd x
  =\sum_{i}\int_\Ocal v^{(i)}\partial_{x_i}\Big(\tfrac12\abs{\nabla d}^2\Big)\dd x+\sum_i\int_\Ocal v^{(i)}\,\partial_{x_i}d\cdot\Delta d\dd x .
\]
The first integral vanishes because $\nabla\cdot v=0$, and the second equals $\scal{(v\cdot\nabla)d}{\Delta d}=-\scal{\Bt(v,d)}{\An d}$. For general $v\in V$ we argue by density, since by Lemma~\ref{lem:nonlinear}(iii),(iv) both sides of \eqref{eq:cancel} are continuous in $v\in V$ for fixed $d\in\DA$. The last statement follows from $\scal{\Bt(v,d)}{d}=0$.
\end{proof}

\begin{remark}
Identity \eqref{eq:cancel} is what makes the energy $\abs u^2+\abs{\nabla d}^2+\int_\Ocal\Ft(\abs d^2)\dd x$ of the deterministic system non-increasing \cite{LinLiu95}. It concerns the diagonal $M(d,d)$ only: the bilinear analogue $\dual{M(d_1,d_2)}{v}=-\scal{\Bt(v,d_1)}{\An d_2}$ is false in general. This is why, in the fluctuation equations, the quadratic term $M(e,e)$ cancels while the cross terms $M(e,d^0)+M(d^0,e)$ must be estimated (Lemma~\ref{lem:coercive}).
\end{remark}

\subsection{The polynomial potential}\label{sec:potential}
\begin{assumption}\label{ass:F}
$N\in\N=\{1,2,\dots\}$ and $\Ft\colon\R\to\R$ is a polynomial of degree $N+1$ with positive leading coefficient and $\Ft(0)=0$. We set $\ft=\Ft'$ and
\[
  f(d)=\ft(\abs d^2)\,d,\qquad d\in\R^3 .
\]
\end{assumption}
Thus $f(d)=\nabla_d\big[\frac12\Ft(\abs d^2)\big]$; the Ginzburg--Landau case corresponds to $\Ft(r)=\frac12(r-1)^2-\frac12$ and $N=1$. Neither $f$ nor $Df$ is globally Lipschitz. The following lemma shows that $f$ is nevertheless one-sided Lipschitz, which is the property used in all energy estimates.

\begin{lemma}\label{lem:f}
Under Assumption~\ref{ass:F} the following hold.
\begin{enumerate}[label=\textup{(\roman*)}]
\item $Df(d)=\ft(\abs d^2)I+2\ft'(\abs d^2)\,d\otimes d$, a symmetric matrix with eigenvalues $\ft(r)$ \textup{(}on $d^\perp$\textup{)} and $\ft(r)+2r\ft'(r)$ \textup{(}in the direction of $d$\textup{)}, where $r=\abs d^2$.
\item There exists $\kappa\ge0$ such that $Df(d)\xi\cdot\xi\ge-\kappa\abs\xi^2$ for all $d,\xi\in\R^3$. In particular
\[
  \big(f(a)-f(b)\big)\cdot(a-b)\ge-\kappa\abs{a-b}^2,\qquad a,b\in\R^3 .
\]
\item There is $C$ such that for all $a,b,d\in\R^3$: 
$$
\abs{f(d)}\le C(1+\abs d^{2N+1}), \quad \abs{Df(d)}\le C(1+\abs d^{2N}), \quad \abs{D^2f(d)}\le C(1+\abs d^{2N-1})
$$
and
$$
\abs{Df(a)-Df(b)}\le C(1+\abs a^{2N-1}+\abs b^{2N-1})\abs{a-b}\,.
$$
\end{enumerate}
\end{lemma}

\begin{proof}
(i) is a direct computation. (ii) The polynomials $\ft(r)$ and $\ft(r)+2r\ft'(r)$ have degree $N$ and leading coefficients $(N+1)c_{N+1}$ and $(2N+1)(N+1)c_{N+1}$, where $c_{N+1}>0$ is the leading coefficient of $\Ft$. Hence both are bounded from below on $[0,\infty)$, and we let $\kappa$ be the negative part of the smaller of the two lower bounds. The second statement follows from $(f(a)-f(b))\cdot(a-b)=\int_0^1Df(b+s(a-b))(a-b)\cdot(a-b)\dd s$. (iii) $f$ is a polynomial map of degree $2N+1$; for the last estimate use the mean value theorem and $\abs{b+s(a-b)}\le\max(\abs a,\abs b)$.
\end{proof}

In function spaces, the one-sided bound holds exactly at the $\LL$ level, and at the $\HH$ level up to a superlinear remainder that is controlled by the $\DA$-norm of the reference function.

\begin{lemma}\label{lem:fspace}
There is a constant $C$, depending only on $\Ocal$ and $\Ft$, such that for all $d,D\in\HH$ and $d^0,e\in\DA$:
\begin{enumerate}[label=\textup{(\alph*)}]
\item $\abs{f(d)}\le C(1+\hone d^{2N+1})$ and $\abs{Df(d)D}\le C(1+\hone d^{2N})\hone D$;
\item $\scal{f(d^0+e)-f(d^0)}{e}\ge-\kappa\abs e^2$;
\item with $\phi(d^0)=(1+\hone{d^0}^{2N})(1+\hda{d^0})$,
\[
  \scal{f(d^0+e)-f(d^0)}{\An e}\ge-\kappa\abs{\nabla e}^2-C\phi(d^0)\big(\hone e^2+\hone e^{2N+1}\big);
\]
\item the Taylor remainder $\mathcal R(d,D)=f(d+D)-f(d)-Df(d)D=\int_0^1(1-r)D^2f(d+rD)[D,D]\dd r$ satisfies
\[
  \abs{\mathcal R(d,D)}\le C\big(1+\hone d^{2N-1}+\hone D^{2N-1}\big)\hone D^2 .
\]
\end{enumerate}
\end{lemma}

\begin{proof}
(a) By Lemma~\ref{lem:f}(iii), H\"older's inequality and \eqref{eq:GN}, $\abs{Df(d)D}\le C\norm{1+\abs d^{2N}}_{L^4}\norm D_{L^4}\le C(1+\hone d^{2N})\hone D$; the bound on $f(d)$ is similar.

(b) is Lemma~\ref{lem:f}(ii) applied pointwise.

(c) Since $\DA\subset\HHH\subset C(\overline\Ocal)$, the functions $f(d^0+e)$ and $f(d^0)$ belong to $\HH$ with $\partial_{x_i}f(d)=Df(d)\partial_{x_i}d$, and Green's formula yields
\[
  \scal{f(d^0+e)-f(d^0)}{\An e}=\sum_{i=1}^2\scal{Df(d^0+e)\partial_{x_i}e}{\partial_{x_i}e}+\sum_{i=1}^2\scal{\big[Df(d^0+e)-Df(d^0)\big]\partial_{x_i}d^0}{\partial_{x_i}e}.
\]
The first sum is $\ge-\kappa\abs{\nabla e}^2$ by Lemma~\ref{lem:f}(ii). By Lemma~\ref{lem:f}(iii) and $\abs{d^0+e}\le\abs{d^0}+\abs e$, the second sum is bounded in absolute value by $C(I_1+I_2)$, where
\[
  I_1=\int_\Ocal\big(1+\abs{d^0}^{2N-1}\big)\abs e\abs{\nabla d^0}\abs{\nabla e}\dd x,\qquad I_2=\int_\Ocal\abs e^{2N}\abs{\nabla d^0}\abs{\nabla e}\dd x .
\]
By H\"older's inequality with exponents $(8,8,4,2)$, \eqref{eq:GN} and \eqref{eq:GNgrad},
\begin{align*}
  I_1&\le\norm{1+\abs{d^0}^{2N-1}}_{L^8}\norm e_{L^8}\norm{\nabla d^0}_{L^4}\abs{\nabla e}\\
  &\le C\big(1+\hone{d^0}^{2N-1}\big)\hone{d^0}^{1/2}\hda{d^0}^{1/2}\hone e^2\le C\phi(d^0)\hone e^2,
\end{align*}
and, with exponents $(4,4,2)$ applied to $\abs e^{2N}$, $\abs{\nabla d^0}$ and $\abs{\nabla e}$,
\[
  I_2\le\norm e_{L^{8N}}^{2N}\norm{\nabla d^0}_{L^4}\abs{\nabla e}\le C\phi(d^0)\hone e^{2N+1}.
\]

(d) The integral formula is Taylor's formula. Since $\abs{d+rD}\le\abs d+\abs D$, Lemma~\ref{lem:f}(iii) gives
\[
  \abs{\mathcal R(d,D)}\le C\norm{(1+\abs d^{2N-1}+\abs D^{2N-1})\abs D^2}_{L^2}\le C\norm{1+\abs d^{2N-1}+\abs D^{2N-1}}_{L^4}\norm D_{L^8}^2,
\]
and we conclude with \eqref{eq:GN}.
\end{proof}

\subsection{The stochastic setting and the geometric noise}\label{sec:noise}
\emph{Probability space and Wiener processes.} Let $(\Omega,\Fc,(\Fc_t)_{t\in[0,T]},\P)$ be a filtered probability space satisfying the usual conditions. We use the following convention for both noises; see \cite[Sections~4.1--4.2]{DPZ14} and \cite[Section~2.5]{LiuRockner15}.
\begin{itemize}[leftmargin=*]
\item \emph{Velocity noise.} Let $\Hz\subset H$ be a separable Hilbert space such that the inclusion $\iota\colon\Hz\hookrightarrow H$ is Hilbert--Schmidt, and let $Q=\iota\iota^*$. Then $Q$ is a non-negative trace-class operator on $H$, and $Q^{1/2}H=\Hz$ with $\abs{Q^{1/2}h}_{\Hz}=\abs h$ for $h\in(\ker Q)^\perp$. Let $W_1$ be an $H$-valued $Q$-Wiener process. Equivalently, $W_1(t)=\sum_j\beta_j(t)e_j$, where $(e_j)$ is an orthonormal basis of $\Hz$, the $\beta_j$ are independent real Brownian motions, and the series converges in $L^2(\Omega;C([0,T];H))$. Stochastic integrals $\int_0^t\Psi\dd W_1=\sum_j\int_0^t\Psi e_j\dd\beta_j$ are defined for progressively measurable integrands $\Psi$ with values in the space $\LHS=\LHS(\Hz,H)$ of Hilbert--Schmidt operators from $\Hz$ to $H$.
\item \emph{Director noise.} Let $W_2=(W_2^k)_{k\in\N}$ be a sequence of independent real Brownian motions, independent of the $\beta_j$. For progressively measurable $\Psi$ with values in $\LHSd=\LHS(\ell^2,\HH)$ we set $\int_0^t\Psi\dd W_2=\sum_k\int_0^t\Psi\mathrm e_k\dd W_2^k$, where $(\mathrm e_k)$ is the canonical basis of $\ell^2$. To regard $W_2$ as a process with continuous paths, we fix a Hilbert space $U_2\supset\ell^2$ such that the inclusion $\ell^2\hookrightarrow U_2$ is Hilbert--Schmidt (for instance a weighted $\ell^2$ space with summable weights); then $W_2(t)=\sum_kW_2^k(t)\mathrm e_k$ converges in $L^2(\Omega;C([0,T];U_2))$.
\end{itemize}
We write $W=(W_1,W_2)$, which is a Wiener process with values in $U_1=H\times U_2$ and reproducing kernel Hilbert space $\mathcal U=\Hz\times\ell^2$.

\emph{The geometric noise.} Let $(h_k)_{k\in\N}$ be vector fields $h_k\colon\Ocal\to\R^3$ such that
\begin{equation}\label{eq:hnorm}
  \mathsf h:=\sum_{k\ge1}\big(\norm{h_k}_{L^\infty}+\norm{\nabla h_k}_{L^3}\big)^2<\infty .
\end{equation}
For $d\in\LL$ and $\rho=(\rho_k)\in\ell^2$ we set
\[
  G_k(d)=d\times h_k,\qquad G(d)\rho=\sum_{k\ge1}\rho_k\,G_k(d),\qquad G^2(d)=\sum_{k\ge1}G_k\big(G_k(d)\big)=\sum_{k\ge1}(d\times h_k)\times h_k ,
\]
and $G(d)\dd W_2=\sum_kG_k(d)\dd W_2^k$. As usual, the Stratonovich equation \eqref{eq:intro-d} is understood through its It\^o form. Since each $G_k$ is linear, the Stratonovich correction is $\frac12\sum_kG_k(G_k(d))$, that is,
\[
  \sum_{k\ge1}G_k(d)\circ\dd W_2^k:=G(d)\dd W_2+\tfrac12G^2(d)\dd t .
\]
A single noise $(d\times h)\circ\dd W_2$, as in \cite{BHR19a,BMP20}, corresponds to $h_1=h$ and $h_k=0$ for $k\ge2$.

\begin{lemma}\label{lem:G}
Assume \eqref{eq:hnorm} and let $C_h=C(1+\mathsf h)^2$, with $C$ depending only on $\Ocal$. Then for every $d\in\HH$:
\begin{enumerate}[label=\textup{(\roman*)}]
\item the series defining $G(d)\rho$ and $G^2(d)$ converge in $\HH$, $G(d)\in\LHSd$, and
\[
  \abs{G(d)}_{\LHS(\ell^2,\LL)}^2+\abs{G^2(d)}^2\le C_h\abs d^2,\qquad \abs{G(d)}_{\LHSd}^2+\hone{G^2(d)}^2\le C_h\hone d^2 ;
\]
\item $G_k(d)\cdot d=0$ for every $k$, and $G^2(d)\cdot d=-\sum_k\abs{G_k(d)}^2$ pointwise. In particular $\scal{G(d)\rho}{d}=0$ for every $\rho\in\ell^2$, and
\[
  \scal{G^2(d)}{d}+\abs{G(d)}_{\LHS(\ell^2,\LL)}^2=0 .
\]
\end{enumerate}
\end{lemma}

\begin{proof}
(i) Write $\mathsf h_k=\norm{h_k}_{L^\infty}+\norm{\nabla h_k}_{L^3}$. We have $\abs{G_k(d)}\le\norm{h_k}_{L^\infty}\abs d$ and $\nabla(d\times h_k)=\nabla d\times h_k+d\times\nabla h_k$, with 
$$
\abs{d\times\nabla h_k}\le\norm d_{L^6}\norm{\nabla h_k}_{L^3}\le C\hone d\norm{\nabla h_k}_{L^3}
$$
 by \eqref{eq:GN}; hence $\hone{G_k(d)}\le C\mathsf h_k\hone d$. In the same way, $\nabla((d\times h_k)\times h_k)$ is a sum of terms bounded pointwise by $\abs{\nabla d}\abs{h_k}^2$ and $\abs d\abs{\nabla h_k}\abs{h_k}$, so $\hone{G_k(G_k(d))}\le C\mathsf h_k^2\hone d$ and $\abs{G_k(G_k(d))}\le\mathsf h_k^2\abs d$. Summing over $k$ and using \eqref{eq:hnorm} gives (i), since $\abs{G(d)}_{\LHSd}^2=\sum_k\hone{G_k(d)}^2$. (ii) is elementary vector algebra: $(d\times h_k)\cdot d=0$ and $((d\times h_k)\times h_k)\cdot d=-(d\times h_k)\cdot(d\times h_k)$.
\end{proof}

Part (ii) reflects the fact that the noise alone preserves $\abs d$ pointwise (Section~\ref{sec:intro}): in the It\^o formula for $\abs d^2$, the It\^o correction $\scal{G^2(d)}{d}$ and the quadratic variation $\abs{G(d)}_{\LHS(\ell^2,\LL)}^2$ cancel exactly.

\subsection{It\^o's formula and compactness}\label{sec:ito}
We use It\^o's formula for the square of the norm in the two Gelfand triples of Section~\ref{sec:spaces}.

\begin{lemma}\label{lem:Ito}
\begin{enumerate}[label=\textup{(\roman*)},leftmargin=*]
\item Let $u_0\in H$, and let $u\in L^2(0,T;V)$, $g\in L^2(0,T;V')$ and $\Psi\in L^2(0,T;\LHS)$ be progressively measurable and such that $u(t)=u_0+\int_0^tg\dd s+\int_0^t\Psi\dd W_1$ in $V'$ for all $t$, a.s. Then $u$ has a continuous $H$-valued version and, a.s.\ for all $t\in[0,T]$,
\[
  \abs{u(t)}^2=\abs{u_0}^2+2\int_0^t\dual{g}{u}\dd s+2\int_0^t\scal{u}{\Psi\dd W_1}+\int_0^t\abs\Psi_{\LHS}^2\dd s .
\]
\item Let $d_0\in\HH$, and let $d\in L^2(0,T;\DA)$, $g\in L^2(0,T;\LL)$ and $\psi\in L^2(0,T;\LHSd)$ be progressively measurable and such that $d(t)=d_0+\int_0^tg\dd s+\int_0^t\psi\dd W_2$ in $\LL$ for all $t$, a.s. Then $d$ has a continuous $\HH$-valued version and, a.s.\ for all $t\in[0,T]$,
\[
  \hone{d(t)}^2=\hone{d_0}^2+2\int_0^t\scal{g}{(I+\An)d}\dd s+2\int_0^t\sone{\psi\dd W_2}{d}+\int_0^t\abs\psi_{\LHSd}^2\dd s .
\]
\end{enumerate}
\end{lemma}

\begin{proof}
Both statements are instances of the It\^o formula for the square of the norm in a Gelfand triple, see \cite{Pardoux75,KrylovRozovskii81,Krylov13} and \cite[Section~4.2]{LiuRockner15}. For (ii) we use the triple \eqref{eq:triple-d}, in which the duality between $\LL\cong(\DA)'$ and $\DA$ is $\scal{g}{(I+\An)d}$.
\end{proof}

\begin{lemma}[Aubin--Lions--Simon]\label{lem:AL}
Let $X_0\subset X\subset X_1$ be Banach spaces, with $X_0$ compactly embedded in $X$ and $X$ continuously embedded in $X_1$. Then $\{v\in L^2(0,T;X_0):\partial_tv\in L^2(0,T;X_1)\}$ is compactly embedded in $L^2(0,T;X)$.
\end{lemma}

This is \cite[Corollary~4]{Simon87}. We use it with $(X_0,X,X_1)=(V,H,V')$ and $(\DA,\HH,\LL)$.
\section{Setting, well-posedness and main results}\label{sec:results}

We first fix the assumptions and the functional formulation (Section~\ref{sec:setting}) and recall well-posedness and moment bounds (Section~\ref{sec:wp}). We then introduce the linearised dynamics around the deterministic flow (Section~\ref{sec:lindyn}) and state the main results (Section~\ref{sec:main}). Section~\ref{sec:coercive} proves the coercivity estimate on which all proofs rest.

\subsection{Standing assumptions and abstract formulation}\label{sec:setting}
Throughout the paper we assume Assumption~\ref{ass:F} and the following.
\begin{assumption}\label{ass:main}
\begin{enumerate}[label=\textup{(A\arabic*)},leftmargin=*]
\item\label{A:sigma} $\sigma\colon H\to\LHS(\Hz,H)$ is measurable and there are $K,L>0$ such that for all $u,v\in H$
\[
  \abs{\sigma(u)}_{\LHS}^2\le K(1+\abs u^2),\qquad \abs{\sigma(u)-\sigma(v)}_{\LHS}^2\le L\abs{u-v}^2 .
\]
\item\label{A:h} $h_k\in W^{1,3}(\Ocal;\R^3)\cap L^\infty(\Ocal;\R^3)$ for every $k\in\N$, and $\mathsf h<\infty$, see \eqref{eq:hnorm}.
\item\label{A:init} $(u_0,d_0)\in H\times\HH$.
\end{enumerate}
\end{assumption}
For a single director noise, these are the hypotheses of \cite{BMP20}. The saturation condition $\abs{d_0}=1$ is allowed but not needed.

Using the operators of Section~\ref{sec:prelim} and the It\^o form of the Stratonovich integral (Section~\ref{sec:noise}), the system \eqref{eq:intro-u}--\eqref{eq:intro-bc} reads, for $\eps\in[0,1]$,
\begin{align}
  \dd u^\eps+\big[\Ast u^\eps+B(u^\eps)+M(d^\eps)\big]\dd t&=\sqrt\eps\,\sigma(u^\eps)\dd W_1, \label{eq:u}\\
  \dd d^\eps+\big[\An d^\eps+\Bt(u^\eps,d^\eps)+f(d^\eps)\big]\dd t&=\frac\eps2G^2(d^\eps)\dd t+\sqrt\eps\,G(d^\eps)\dd W_2, \label{eq:d}
\end{align}
with $u^\eps(0)=u_0$ and $d^\eps(0)=d_0$. We write $X^\eps=(u^\eps,d^\eps)$. For $\eps=0$ the system is deterministic, and we denote its solution by $X^0=(u^0,d^0)$.

\emph{Convention on constants.} The constants $C$, $C_p$, \dots\ may change from line to line. Unless stated otherwise, they may depend on $p$, $T$, $\Ocal$, $\Ft$, $K$, $L$, $\mathsf h$ and the initial datum $(u_0,d_0)$, but never on $\eps$. Section~\ref{sec:uniform} makes their dependence on the initial datum explicit.

\subsection{Path space, well-posedness and moment bounds}\label{sec:wp}
The natural path space is the separable Banach space
\[
  \KT=\big[C([0,T];H)\cap L^2(0,T;V)\big]\times\big[C([0,T];\HH)\cap L^2(0,T;\DA)\big],
\]
\[
  \norm{(v,e)}_{\KT}^2=\sup_{t\in[0,T]}\big(\abs{v(t)}^2+\hone{e(t)}^2\big)+\int_0^T\big(\norm{v(t)}^2+\hda{e(t)}^2\big)\dd t .
\]
In energy estimates we use, for $(v,e)\in\KT$ and $t\in[0,T]$, the energy functional
\begin{equation}\label{eq:energy-functional}
  \Ec_t(v,e)=\sup_{s\le t}\big(\abs{v(s)}^2+\hone{e(s)}^2\big)+\int_0^t\big(\norm{v}^2+\abs{\nabla e}^2+\abs{\An e}^2\big)\dd s .
\end{equation}
It is equivalent to the norm of $\KT$: $(1+T)^{-1}\norm{(v,e)}_{\KT}^2\le\Ec_T(v,e)\le(1+T)\norm{(v,e)}_{\KT}^2$.

\begin{definition}\label{def:solution}
Let $\eps\in[0,1]$ and let $(\Omega,\Fc,(\Fc_t),\P,W)$ be as in Section~\ref{sec:noise}. A \emph{solution} of \eqref{eq:u}--\eqref{eq:d} is a progressively measurable process $X^\eps=(u^\eps,d^\eps)$ with $\P$-a.s.\ paths in $\KT$ such that, a.s.\ for all $t\in[0,T]$, the integral form of \eqref{eq:u} holds in $V'$ and that of \eqref{eq:d} holds in $\LL$.
\end{definition}

\begin{theorem}\label{thm:wellposed}
Let Assumptions~\ref{ass:F} and \ref{ass:main} hold.
\begin{enumerate}[label=\textup{(\roman*)}]
\item For $\eps=0$ there is a unique $X^0=(u^0,d^0)\in\KT$ solving \eqref{eq:u}--\eqref{eq:d}.
\item For every $\eps\in(0,1]$, pathwise uniqueness holds for \eqref{eq:u}--\eqref{eq:d}, and there is a Borel map
\[
  J^\eps\colon C([0,T];H\times U_2)\to\KT
\]
such that, on any stochastic basis carrying a Wiener process $W$ as above, $X^\eps=J^\eps(W)$ is the unique solution.
\end{enumerate}
\end{theorem}

For a single director noise, part (i) is \cite[Theorem~5.5]{BMP20} with zero control; for the deterministic system see also \cite{LinLiu95,LinLinWang10,Hong11}. Part (ii) follows from the existence of martingale solutions and pathwise uniqueness \cite{BHR19a} together with the Yamada--Watanabe theorem \cite{YamadaWatanabe71,Ondrejat04}; see \cite[Theorems~4.2--4.4]{BMP20}. For countably many director noises, existence follows from the theory of R\"ockner, Shang and Zhang \cite{RSZ24,RSZ26} for SPDEs with fully local monotone coefficients, which extends the locally monotone framework of \cite{LiuRockner10}, and uniqueness from a direct energy estimate; see Appendix~\ref{app:wellposed}. We write $X^{\eps}$, $X^0$ and $J^\eps$ when the initial datum is fixed, and $X^{\eps,x}$, $X^{0,x}$ and $J^\eps(x,\cdot)$ when the dependence on $x=(u_0,d_0)$ matters.

The following moment bounds, uniform in $\eps$, are proved in Appendix~\ref{app:moments}.

\begin{proposition}\label{prop:moments}
For every $p\ge1$ there is $C_p$, independent of $\eps\in[0,1]$ and of the initial datum, such that
\[
  \E\norm{X^\eps}_{\KT}^{2p}\le C_p\big(1+\abs{u_0}^2+\hone{d_0}^2\big)^p .
\]
In particular $\norm{X^0}_{\KT}^2\le C_1(1+\abs{u_0}^2+\hone{d_0}^2)$.
\end{proposition}

\subsection{The linearised dynamics}\label{sec:lindyn}
All fluctuation equations in this paper are linearisations of \eqref{eq:u}--\eqref{eq:d} around the deterministic flow $X^0$. The flow $X^0$ enters the estimates only through the function
\begin{equation}\label{eq:g}
  g(t)=1+\norm{u^0(t)}^2+\abs{u^0(t)}^2\norm{u^0(t)}^2+\abs{\nabla d^0(t)}^2\hda{d^0(t)}^2+\big(1+\hone{d^0(t)}^{4N}\big)\big(1+\hda{d^0(t)}\big),
\end{equation}
which belongs to $L^1(0,T)$ because $X^0\in\KT$.

\emph{The linearised operator.} For $t\in[0,T]$ and $(U,D)\in V\times\DA$ let
\begin{align*}
  \Lin_1(t)(U,D)&=\Ast U+B(U,u^0(t))+B(u^0(t),U)+M(D,d^0(t))+M(d^0(t),D)\in V',\\
  \Lin_2^\flat(t)(U,D)&=\An D+\Bt(U,d^0(t))+\Bt(u^0(t),D)\in\LL,\\
  \Lin_2(t)(U,D)&=\Lin_2^\flat(t)(U,D)+Df(d^0(t))D\in\LL .
\end{align*}
Thus $\Lin=(\Lin_1,\Lin_2)$ is the derivative at $X^0(t)$ of the drift of the deterministic system. The operator $\Lin^\flat=(\Lin_1,\Lin_2^\flat)$ omits the potential; it is used when the difference $f(d^0+e)-f(d^0)$ is kept in nonlinear form and estimated with Lemma~\ref{lem:fspace}.

\emph{The linearised SPDE.} The Gaussian limit $X^1=(U,D)$ solves the linearisation of \eqref{eq:u}--\eqref{eq:d} along $X^0$, driven by the noise frozen at $X^0$:
\begin{align}
  \dd U+\big[\Ast U+B(U,u^0)+B(u^0,U)+M(D,d^0)+M(d^0,D)\big]\dd t&=\sigma(u^0)\dd W_1, \label{eq:lin-u}\\
  \dd D+\big[\An D+\Bt(U,d^0)+\Bt(u^0,D)+Df(d^0)D\big]\dd t&=G(d^0)\dd W_2, \label{eq:lin-d}
\end{align}
with $U(0)=0$ and $D(0)=0$, that is, $\dd X^1+\Lin(t)X^1\dd t=(\sigma(u^0)\dd W_1,G(d^0)\dd W_2)$. Since the noise coefficients are deterministic, the It\^o and Stratonovich forms of \eqref{eq:lin-d} coincide, and the It\^o correction $\frac\eps2G^2(d^\eps)$ leaves no trace at this scale.

\emph{The skeleton map.} Let $\mathcal U=\Hz\times\ell^2$ be the reproducing kernel Hilbert space of $W=(W_1,W_2)$. For $q=(\theta,\rho)\in L^2(0,T;\mathcal U)$, let $\Gam(q)=(U^q,D^q)$ denote the solution of the deterministic \emph{skeleton} system
\begin{align}
  \partial_tU^q+\Ast U^q+B(U^q,u^0)+B(u^0,U^q)+M(D^q,d^0)+M(d^0,D^q)&=\sigma(u^0)\theta, \label{eq:skel-u}\\
  \partial_tD^q+\An D^q+\Bt(U^q,d^0)+\Bt(u^0,D^q)+Df(d^0)D^q&=G(d^0)\rho, \label{eq:skel-d}
\end{align}
with $U^q(0)=0$ and $D^q(0)=0$; this is \eqref{eq:lin-u}--\eqref{eq:lin-d} with the noise replaced by the control $q$. By Proposition~\ref{prop:skeleton}, $\Gam\colon L^2(0,T;\mathcal U)\to\KT$ is a bounded linear operator. It determines both the law of $X^1$ (Theorem~\ref{thm:linear}) and the moderate deviations of $X^\eps$ (Theorem~\ref{thm:MDP}).

\subsection{Main results}\label{sec:main}
\emph{Law of large numbers and central limit theorem.}

\begin{theorem}[Quantitative law of large numbers]\label{thm:LLN}
For every $p\ge1$ there is $C_p$ such that for all $\eps\in(0,1]$
\[
  \E\norm{X^\eps-X^0}_{\KT}^{2p}\le C_p\,\eps^p .
\]
\end{theorem}

\begin{theorem}[The linearised system]\label{thm:linear}
The following assertions hold.
\begin{enumerate}[label=\textup{(\roman*)}]
\item The system \eqref{eq:lin-u}--\eqref{eq:lin-d} has a unique solution $X^1=(U,D)$ with paths in $\KT$. For every $p\ge1$, $\E\norm{X^1}_{\KT}^{2p}\le C_p$.
\item For every $\ell\in\KT^*$, let $k_\ell\in L^2(0,T;\mathcal U)$ be the Riesz representative of the bounded linear functional $\ell\circ\Gam$. Then
\[
  \ell(X^1)=\int_0^T\big\langle k_\ell(t),\dd W(t)\big\rangle_{\mathcal U}\quad\text{a.s.}
\]
In particular, $X^1$ is a centred Gaussian random variable in $\KT$, and $\operatorname{Var}\ell(X^1)=\norm{\ell\circ\Gam}^2=\sup\{\ell(\Gam q)^2:\norm q_{L^2(0,T;\mathcal U)}\le1\}$.
\end{enumerate}
\end{theorem}

\begin{theorem}[Central limit theorem with rate]\label{thm:CLT}
For every $p\ge1$ there is $C_p$ such that for all $\eps\in(0,1]$
\[
  \E\Big\|\frac{X^\eps-X^0}{\sqrt\eps}-X^1\Big\|_{\KT}^{2p}\le C_p\,\eps^{p}.
\]
\end{theorem}

\begin{corollary}\label{cor:CLT}
For every $p\ge1$ the following hold.
\begin{enumerate}[label=\textup{(\roman*)}]
\item $(X^\eps-X^0)/\sqrt\eps\to X^1$ in $L^{2p}(\Omega;\KT)$, and hence in probability and in law on $\KT$; the limit law is the centred Gaussian measure $\Law(X^1)$.
\item The $p$-Wasserstein distance \cite{Villani09} on $\KT$ satisfies
\[
  \mathcal W_p\Big(\Law\big((X^\eps-X^0)/\sqrt\eps\big),\,\Law(X^1)\Big)\le C_p\sqrt\eps .
\]
\item $X^\eps=X^0+\sqrt\eps X^1+\mathcal Q^\eps$ with $\big(\E\norm{\mathcal Q^\eps}_{\KT}^{2p}\big)^{1/2p}\le C_p\,\eps$.
\end{enumerate}
\end{corollary}

\begin{proof}
(i) and (iii) are restatements of Theorem~\ref{thm:CLT}, with $\mathcal Q^\eps=\sqrt\eps\big((X^\eps-X^0)/\sqrt\eps-X^1\big)$. For (ii), the two random variables are defined on the same probability space, so 
$$
\mathcal W_p\le(\E\norm{(X^\eps-X^0)/\sqrt\eps-X^1}_{\KT}^{p})^{1/p};
$$
 we conclude by H\"older's inequality and Theorem~\ref{thm:CLT}.
\end{proof}

\emph{Moderate deviations.} Let $a\colon(0,1]\to(0,\infty)$ satisfy
\begin{equation}\label{eq:a}
  a(\eps)\to\infty\qquad\text{and}\qquad\delta_\eps:=\sqrt\eps\,a(\eps)\to0\qquad\text{as }\eps\downarrow0,
\end{equation}
and define the moderately rescaled fluctuations
\[
  Y^\eps=\frac{X^\eps-X^0}{\sqrt\eps\,a(\eps)} .
\]
\begin{theorem}[Moderate deviation principle]\label{thm:MDP}
Assume \eqref{eq:a}. The family $\{Y^\eps\}_{\eps\in(0,1]}$ satisfies a large deviation principle on $\KT$ with speed $a(\eps)^2$ and good rate function
\[
  I(z)=\inf\Big\{\frac12\int_0^T\big(\abs{\theta(t)}_{\Hz}^2+\abs{\rho(t)}_{\ell^2}^2\big)\dd t:\ q=(\theta,\rho)\in L^2(0,T;\mathcal U),\ \Gam(q)=z\Big\},
\]
with $\inf\emptyset=+\infty$. That is, for every closed set $F\subset\KT$ and every open set $O\subset\KT$,
\[
  \limsup_{\eps\to0}\frac1{a(\eps)^2}\log\P(Y^\eps\in F)\le-\inf_{z\in F}I(z),\qquad
  \liminf_{\eps\to0}\frac1{a(\eps)^2}\log\P(Y^\eps\in O)\ge-\inf_{z\in O}I(z).
\]
\end{theorem}

\begin{remark}\label{rem:CM}
\begin{enumerate}[label=(\alph*),leftmargin=*]
\item Since $\Gam$ is linear and bounded, $I$ is a quadratic form: $I(z)=\frac12\norm z_{\mathscr H}^2$ for $z\in\mathscr H=\Gam(L^2(0,T;\mathcal U))$, with $\norm z_{\mathscr H}=\inf\{\norm q:\Gam q=z\}$, and $I(z)=+\infty$ otherwise. By Theorem~\ref{thm:linear}(ii), the one-dimensional projections of $I$ are those of the Gaussian law of $X^1$, as for sums of independent random vectors \cite{deAcosta92}: for $\ell\in\KT^*$ with $\operatorname{Var}\ell(X^1)>0$,
\[
  \inf\{I(z):\ell(z)=r\}=\frac{r^2}{2\operatorname{Var}\ell(X^1)}.
\]
More precisely, $\mathscr H$ is the Cameron--Martin space of $\Law(X^1)$ (Corollary~\ref{cor:gaussian}).
\item The skeleton \eqref{eq:skel-u}--\eqref{eq:skel-d} is the linearisation at $q=0$ of the nonlinear skeleton equation of the LDP, whose solution for $q=0$ is $X^0$ (see \cite[Section~5.4]{BMP20} for a single director noise). In the normalisation of the Budhiraja--Dupuis representation \cite{BD00} used here, both the LDP and the MDP rate functions carry the factor $\frac12$. The rate function displayed in \cite[(5.12)]{BMP20} omits this factor, which appears to be a misprint.
\end{enumerate}
\end{remark}

As an application, we obtain sharp moderate-deviation asymptotics for scalar observables of the flow. A function $\Phi\colon\KT\to\R$ is \emph{$C^{1,1}$ at $X^0$} if there are $\ell_\Phi\in\KT^*$ and $C_\Phi<\infty$ with 
$$
\abs{\Phi(X^0+z)-\Phi(X^0)-\ell_\Phi(z)}\le C_\Phi\norm z_{\KT}^2
$$
 for $\norm z_{\KT}\le1$. Examples include continuous linear functionals, the elastic energy at the final time, $\Phi(u,d)=\frac12\abs{\nabla d(T)}^2$, and the total viscous dissipation, $\Phi(u,d)=\int_0^T\norm{u(t)}^2\dd t$.

\begin{corollary}\label{cor:observables}
Assume \eqref{eq:a}, let $\Phi$ be $C^{1,1}$ at $X^0$ and set $s_\Phi^2=\operatorname{Var}\ell_\Phi(X^1)$. If $s_\Phi>0$, then for every $r>0$
\[
  \lim_{\eps\to0}\frac1{a(\eps)^2}\log\P\Big(\Phi(X^\eps)-\Phi(X^0)\ge r\sqrt\eps\,a(\eps)\Big)=-\frac{r^2}{2s_\Phi^2}.
\]
If $s_\Phi=0$, the left-hand side equals $-\infty$ for every $r>0$.
\end{corollary}

\emph{Uniformity in the initial data.} In Section~\ref{sec:uniform} we show that Theorems~\ref{thm:LLN}, \ref{thm:linear}(i), \ref{thm:CLT} and \ref{thm:MDP} hold uniformly over bounded sets of initial data in $H\times\HH$ (Corollary~\ref{cor:unif-LLN-CLT} and Theorem~\ref{thm:uniform}), and we derive uniform Gaussian tube estimates (Corollary~\ref{cor:tube}).

\subsection{Coercivity of the linearised operator}\label{sec:coercive}
The following lemma is the basic energy estimate for all fluctuation equations. Once Lemma~\ref{lem:cancel} has removed the cubic terms, it controls every remaining coefficient by the deterministic function $g$ of \eqref{eq:g}.

\begin{lemma}\label{lem:coercive}
There is a constant $C_*$, depending only on $\Ocal$ and $\Ft$, such that for a.e.\ $t\in[0,T]$ and all $(U,D)\in V\times\DA$:
\begin{enumerate}[label=\textup{(\roman*)}]
\item the following inequality holds, and remains true with $\Lin_2$ replaced by $\Lin_2^\flat$:
\[
  \dual{\Lin_1(t)(U,D)}{U}+\scal{\Lin_2(t)(U,D)}{(I+\An)D}\ge\frac12\big(\norm U^2+\abs{\nabla D}^2+\abs{\An D}^2\big)-C_*g(t)\big(\abs U^2+\hone D^2\big);
\]
\item for every $(U,D)\in L^\infty(0,T;H\times\HH)\cap L^2(0,T;V\times\DA)$, with $\norm{(U,D)}_{\KT}$ computed with the essential supremum in time,
\[
  \int_0^T\Big(\abs{\Lin_1(t)(U,D)}_{V'}^2+\abs{\Lin_2(t)(U,D)}^2\Big)\dd t\le C_*\Big(1+\int_0^Tg\dd t\Big)\norm{(U,D)}_{\KT}^2 .
\]
\end{enumerate}
\end{lemma}

\begin{proof}
(i) We have $\dual{\Ast U}{U}=\norm U^2$, $\scal{\An D}{(I+\An)D}=\abs{\nabla D}^2+\abs{\An D}^2$ and, by Lemma~\ref{lem:nonlinear}(i),(iii), $\dual{B(u^0,U)}{U}=0$ and $\scal{\Bt(u^0,D)}{D}=0$. Let $\eta>0$. We drop the argument $t$ and recall that $\hda D^2=\abs D^2+\abs{\An D}^2$. By Lemma~\ref{lem:nonlinear} and Young's inequality,
\begin{align*}
  \abs{\dual{B(U,u^0)}{U}}&\le C\abs U\norm U\norm{u^0}\le\eta\norm U^2+C_\eta\norm{u^0}^2\abs U^2,\\
  \abs{\dual{M(D,d^0)+M(d^0,D)}{U}}&\le C\norm U\abs{\nabla D}^{1/2}\hda D^{1/2}\abs{\nabla d^0}^{1/2}\hda{d^0}^{1/2}\\
  &\le\eta\norm U^2+\eta\hda D^2+C_\eta\abs{\nabla d^0}^2\hda{d^0}^2\abs{\nabla D}^2,\\
  \abs{\scal{\Bt(U,d^0)}{(I+\An)D}}&\le C\norm U\abs{\nabla d^0}\hone D+C\abs U^{1/2}\norm U^{1/2}\abs{\nabla d^0}^{1/2}\hda{d^0}^{1/2}\abs{\An D}\\
  &\le2\eta\norm U^2+\eta\abs{\An D}^2+C_\eta\big(\abs{\nabla d^0}^2\hone D^2+\abs{\nabla d^0}^2\hda{d^0}^2\abs U^2\big),\\
  \abs{\scal{\Bt(u^0,D)}{\An D}}&\le C\abs{u^0}^{1/2}\norm{u^0}^{1/2}\abs{\nabla D}^{1/2}\hda D^{3/2}\le\eta\hda D^2+C_\eta\abs{u^0}^2\norm{u^0}^2\abs{\nabla D}^2 .
\end{align*}
For the potential, Lemma~\ref{lem:f}(ii) and Lemma~\ref{lem:fspace}(a) give
\[
  \scal{Df(d^0)D}{(I+\An)D}\ge-\kappa\abs D^2-\eta\abs{\An D}^2-C_\eta\big(1+\hone{d^0}^{4N}\big)\hone D^2 .
\]
Since $\abs{\nabla d^0}^2\le\hone{d^0}^2\le 1+\hone{d^0}^{4N}$, every coefficient above is bounded by a multiple of $g$. With $\eta=1/16$ the absorbed terms add up to at most $\frac14\norm U^2+\frac14\abs{\An D}^2+\frac18\abs D^2$, and $\frac18\abs D^2\le\frac18g\hone D^2$ because $g\ge1$. This yields (i) for $\Lin_2$; for $\Lin_2^\flat$ simply omit the estimate of the potential.

(ii) By Lemma~\ref{lem:nonlinear}(ii)--(iv) and Lemma~\ref{lem:fspace}(a),
\begin{align*}
  \abs{\Lin_1(U,D)}_{V'}^2&\le C\big(\norm U^2+\abs U\norm U\abs{u^0}\norm{u^0}+\abs{\nabla D}\hda D\abs{\nabla d^0}\hda{d^0}\big)\\
  &\le C\big(\norm U^2+\hda D^2+g\,(\abs U^2+\abs{\nabla D}^2)\big),\\
  \abs{\Lin_2(U,D)}^2&\le C\big(\hda D^2+\abs U\norm U\abs{\nabla d^0}\hda{d^0}+\abs{u^0}\norm{u^0}\abs{\nabla D}\hda D+(1+\hone{d^0}^{4N})\hone D^2\big)\\
  &\le C\big(\norm U^2+\hda D^2+g\,(\abs U^2+\hone D^2)\big).
\end{align*}
Integrating over $[0,T]$ and bounding $\abs U^2+\hone D^2$ by its supremum gives (ii).
\end{proof}
\section{Quantitative law of large numbers}\label{sec:LLN}

In this section we prove Theorem~\ref{thm:LLN}. We derive an energy inequality for $X^\eps-X^0$ with deterministic coefficients (Section~\ref{sec:LLN-energy}), localise the superlinear term produced by the potential (Section~\ref{sec:LLN-local}), and remove the localisation with the moment bounds of Proposition~\ref{prop:moments} (Section~\ref{sec:LLN-proof}).

\subsection{The equation for the deviation}\label{sec:LLN-eq}
Fix $\eps\in(0,1]$ and set
\[
  w=u^\eps-u^0,\qquad e=d^\eps-d^0 .
\]
Since $B$, $\Bt$ and $M$ are bilinear,
\begin{align*}
  B(u^\eps)-B(u^0)&=B(w,u^0)+B(u^0,w)+B(w,w),\\
  M(d^\eps)-M(d^0)&=M(e,d^0)+M(d^0,e)+M(e,e),\\
  \Bt(u^\eps,d^\eps)-\Bt(u^0,d^0)&=\Bt(w,d^0)+\Bt(u^0,e)+\Bt(w,e).
\end{align*}
Hence, with the operators $\Lin_1$ and $\Lin_2^\flat$ of Section~\ref{sec:lindyn}, the pair $(w,e)$ satisfies $w(0)=0$, $e(0)=0$ and
\begin{align}
  \dd w+\big[\Lin_1(t)(w,e)+B(w,w)+M(e,e)\big]\dd t&=\sqrt\eps\,\sigma(u^\eps)\dd W_1, \label{eq:w}\\
  \dd e+\big[\Lin_2^\flat(t)(w,e)+\Bt(w,e)+f(d^0+e)-f(d^0)\big]\dd t&=\frac\eps2G^2(d^\eps)\dd t+\sqrt\eps\,G(d^\eps)\dd W_2. \label{eq:e}
\end{align}
\emph{It\^o's formula applies to $(w,e)$.} Let $(v,d)$ be either $X^\eps$ or $X^0$; its paths lie in $\KT$. By Lemma~\ref{lem:nonlinear}(ii)--(iv),
\[
  \abs{B(v)}_{V'}\le C\abs v\norm v,\qquad \abs{M(d)}_{V'}\le C\abs{\nabla d}\hda d,\qquad \abs{\Bt(v,d)}^2\le C\abs v\norm v\abs{\nabla d}\hda d,
\]
so that $\Ast v+B(v)+M(d)\in L^2(0,T;V')$ and $\An d+\Bt(v,d)\in L^2(0,T;\LL)$; moreover $f(d)\in L^\infty(0,T;\LL)$ by Lemma~\ref{lem:fspace}(a), and $G^2(d)\in L^\infty(0,T;\HH)$ by Lemma~\ref{lem:G}(i). Similarly, $\sigma(u^\eps)\in L^\infty(0,T;\LHS)$ by \ref{A:sigma} and $G(d^\eps)\in L^\infty(0,T;\LHSd)$ by Lemma~\ref{lem:G}(i). Hence, almost surely, the drifts of \eqref{eq:w} and \eqref{eq:e} belong to $L^2(0,T;V')$ and $L^2(0,T;\LL)$, and their diffusion coefficients to $L^2(0,T;\LHS)$ and $L^2(0,T;\LHSd)$. Lemma~\ref{lem:Ito} therefore applies to
\[
  \mathcal Y(t)=\abs{w(t)}^2+\hone{e(t)}^2,\qquad \Dc(t)=\norm{w(t)}^2+\abs{\nabla e(t)}^2+\abs{\An e(t)}^2 .
\]
Note that $\Ec_t(w,e)=\sup_{s\le t}\mathcal Y(s)+\int_0^t\Dc\dd s$, see \eqref{eq:energy-functional}.

\subsection{The energy inequality}\label{sec:LLN-energy}
Let
\[
  \Pi^\eps=1+\abs{u^\eps}^2+\hone{d^\eps}^2 .
\]

\begin{lemma}\label{lem:LLN-energy}
There is a constant $C$, depending only on $\Ocal$, $\Ft$, $K$ and $\mathsf h$, such that a.s.\ for all $t\in[0,T]$
\begin{equation}\label{eq:LLN-energy}
  \mathcal Y(t)+\int_0^t\Dc\dd s\le C\int_0^tg(s)\Big(\mathcal Y(s)+\hone{e(s)}^{2N+1}\Big)\dd s+C\eps\int_0^t\Pi^\eps(s)\dd s+\mathcal M(t),
\end{equation}
where $\mathcal M$ is the continuous local martingale
\[
  \mathcal M(t)=2\sqrt\eps\int_0^t\scal{w}{\sigma(u^\eps)\dd W_1}+2\sqrt\eps\int_0^t\sone{G(d^\eps)\dd W_2}{e} .
\]
\end{lemma}

\begin{proof}
By Lemma~\ref{lem:Ito}, applied to \eqref{eq:w} in the triple $V\subset H\subset V'$ and to \eqref{eq:e} in the triple $\DA\subset\HH\subset\LL$, and since $\mathcal Y(0)=0$,
\begin{align*}
  \mathcal Y(t)&=-2\int_0^t\Big[\dual{\Lin_1(w,e)}{w}+\scal{\Lin_2^\flat(w,e)}{(I+\An)e}\Big]\dd s\\
  &\quad-2\int_0^t\Big[\dual{B(w,w)}{w}+\dual{M(e,e)}{w}+\scal{\Bt(w,e)}{(I+\An)e}\Big]\dd s\\
  &\quad-2\int_0^t\scal{f(d^0+e)-f(d^0)}{(I+\An)e}\dd s+\eps\int_0^t\scal{G^2(d^\eps)}{(I+\An)e}\dd s\\
  &\quad+\eps\int_0^t\Big[\abs{\sigma(u^\eps)}_{\LHS}^2+\abs{G(d^\eps)}_{\LHSd}^2\Big]\dd s+\mathcal M(t).
\end{align*}
We estimate the integrands for a.e.\ $s\in[0,t]$, where $d^0(s),e(s)\in\DA$.

\emph{Quadratic terms.} The second bracket vanishes: $\dual{B(w,w)}{w}=0$ by Lemma~\ref{lem:nonlinear}(i), and $\dual{M(e,e)}{w}+\scal{\Bt(w,e)}{(I+\An)e}=0$ by Lemma~\ref{lem:cancel} with $v=w$ and $d=e$.

\emph{Linear terms.} By Lemma~\ref{lem:coercive}(i) for $\Lin^\flat$, the first bracket is bounded below by $\frac12\Dc-C_*g\mathcal Y$, so the first line is at most $-\int_0^t\Dc\dd s+2C_*\int_0^tg\,\mathcal Y\dd s$.

\emph{Potential.} By Lemma~\ref{lem:fspace}(b),(c),
\[
  -2\scal{f(d^0+e)-f(d^0)}{(I+\An)e}\le2\kappa\hone e^2+2C\phi(d^0)\big(\hone e^2+\hone e^{2N+1}\big).
\]
Since $1+x^{2N}\le2(1+x^{4N})$ for $x\ge0$, the definition \eqref{eq:g} gives $\phi(d^0)\le2g$; together with $g\ge1$, the right-hand side is at most $Cg\big(\mathcal Y+\hone e^{2N+1}\big)$.

\emph{It\^o correction.} Since $G^2(d^\eps)\in\HH$ and $e\in\DA$, Green's formula gives $\scal{G^2(d^\eps)}{(I+\An)e}=\sone{G^2(d^\eps)}{e}$. By Lemma~\ref{lem:G}(i) and Young's inequality,
\[
  \eps\scal{G^2(d^\eps)}{(I+\An)e}\le\eps C_h^{1/2}\hone{d^\eps}\hone e\le C\eps\,\Pi^\eps+C g\,\mathcal Y .
\]

\emph{Quadratic variation.} By \ref{A:sigma} and Lemma~\ref{lem:G}(i), $\abs{\sigma(u^\eps)}_{\LHS}^2+\abs{G(d^\eps)}_{\LHSd}^2\le K(1+\abs{u^\eps}^2)+C_h\hone{d^\eps}^2\le C\,\Pi^\eps$.

Inserting these bounds, moving $\int_0^t\Dc\dd s$ to the left-hand side, and noting that each constant depends only on $\Ocal$, $\Ft$, $K$ and $\mathsf h$, we obtain \eqref{eq:LLN-energy}.
\end{proof}

\subsection{Localisation and moment estimates}\label{sec:LLN-local}
The term $\hone e^{2N+1}$ in \eqref{eq:LLN-energy} is superlinear in $\mathcal Y$, so Gronwall's lemma cannot be applied directly. We localise it. Let
\begin{equation}\label{eq:tau-LLN}
  \tau=\inf\{t\in[0,T]:\hone{e(t)}\ge1\}\wedge T .
\end{equation}
Since $e$ is adapted with continuous paths in $\HH$ and $e(0)=0$, $\tau$ is a stopping time and $\tau>0$. On $[0,\tau]$ we have $\hone e\le1$, hence
\[
  \hone e^{2N+1}\le\hone e^2\le\mathcal Y\qquad\text{on }[0,\tau].
\]
Let $Z$ be the localised energy
\[
  Z(t)=\sup_{s\le t}\mathcal Y(s\wedge\tau)+\int_0^{t\wedge\tau}\Dc\dd s,\qquad t\in[0,T].
\]
It is adapted, continuous and nondecreasing, and $Z(T)=\Ec_\tau(w,e)\le\Ec_T(w,e)$.

\begin{lemma}\label{lem:LLN-local}
For every $p\ge1$ there is $C_p$, independent of $\eps$, such that
\[
  \E Z(T)^p\le C_p\,\eps^p .
\]
In particular $\P(\tau<T)\le C_q\,\eps^q$ for every $q\ge1$.
\end{lemma}

\begin{proof}
\emph{Step 1: a pathwise bound.} Let $t\in[0,T]$ and $s\le t$. For $r\le s\wedge\tau$ we have $\mathcal Y(r)=\mathcal Y(r\wedge\tau)\le Z(r)$, so
\[
  \int_0^{s\wedge\tau}g\big(\mathcal Y+\hone e^{2N+1}\big)\dd r\le2\int_0^{s\wedge\tau}g\,\mathcal Y\dd r\le2\int_0^tg(r)Z(r)\dd r .
\]
Evaluating \eqref{eq:LLN-energy} at $s\wedge\tau$, both $\mathcal Y(s\wedge\tau)$ and $\int_0^{s\wedge\tau}\Dc\dd r$ are bounded by the right-hand side. Taking the supremum over $s\le t$ and adding the two bounds, we obtain
\begin{equation}\label{eq:LLN-Z}
  Z(t)\le C\int_0^tg(r)Z(r)\dd r+C\eps\int_0^T\Pi^\eps\dd s+2\sup_{s\le t}\abs{\mathcal M(s\wedge\tau)} .
\end{equation}

\emph{Step 2: the $p$-th power.} By the elementary inequality $(a+b+c)^p\le3^{p-1}(a^p+b^p+c^p)$ and H\"older's inequality with respect to the finite measure $g(r)\dd r$ on $[0,t]$,
\[
  Z(t)^p\le C_p\Big(\int_0^Tg\dd r\Big)^{p-1}\int_0^tg(r)Z(r)^p\dd r+C_p\eps^p\Big(\int_0^T\Pi^\eps\dd s\Big)^p+C_p\sup_{s\le t}\abs{\mathcal M(s\wedge\tau)}^p .
\]

\emph{Step 3: the martingale term.} The quadratic variation of $\mathcal M$ satisfies, by \ref{A:sigma} and Lemma~\ref{lem:G}(i),
\begin{align*}
  \dd\langle\mathcal M\rangle_s&=4\eps\Big(\abs{\sigma(u^\eps)^*w}_{\Hz}^2+\sum_{k\ge1}\sone{G_k(d^\eps)}{e}^2\Big)\dd s\\
  &\le4\eps\Big(\abs{\sigma(u^\eps)}_{\LHS}^2\abs w^2+\abs{G(d^\eps)}_{\LHSd}^2\hone e^2\Big)\dd s\le C\eps\,\Pi^\eps\mathcal Y\dd s .
\end{align*}
Since $\mathcal Y(s)\le Z(t)$ for $s\le t\wedge\tau$, the Burkholder--Davis--Gundy inequality \cite[Chapter~IV]{RevuzYor99} and Young's inequality give, for every $\gamma>0$,
\begin{align*}
  \E\sup_{s\le t}\abs{\mathcal M(s\wedge\tau)}^p&\le C_p\,\E\langle\mathcal M\rangle_{t\wedge\tau}^{p/2}\le C_p\eps^{p/2}\,\E\Big[Z(t)^{p/2}\Big(\int_0^T\Pi^\eps\dd s\Big)^{p/2}\Big]\\
  &\le\gamma\,\E Z(t)^p+C_{p,\gamma}\,\eps^p\,\E\Big(\int_0^T\Pi^\eps\dd s\Big)^p .
\end{align*}

\emph{Step 4: Gronwall's lemma.} Since $Z(T)\le\Ec_T(w,e)\le2(1+T)\big(\norm{X^\eps}_{\KT}^2+\norm{X^0}_{\KT}^2\big)$, Proposition~\ref{prop:moments} shows that $m_p(t)=\E Z(t)^p$ is finite, and that
\[
  \E\Big(\int_0^T\Pi^\eps\dd s\Big)^p\le T^p\,\E\big(1+\norm{X^\eps}_{\KT}^2\big)^p\le C_p .
\]
Taking expectations in Step~2, inserting Step~3 with $\gamma$ so small that $C_p\gamma\le\frac12$, and absorbing $\frac12m_p(t)$ into the left-hand side, we obtain
\[
  m_p(t)\le C_p\Big(\int_0^Tg\dd r\Big)^{p-1}\int_0^tg(r)\,m_p(r)\dd r+C_p\,\eps^p,\qquad t\in[0,T].
\]
The function $m_p$ is bounded and nondecreasing, and $g\in L^1(0,T)$ is \emph{deterministic}. Gronwall's lemma therefore gives
\[
  \E Z(T)^p=m_p(T)\le C_p\,\eps^p\exp\Big(C_p\Big(\int_0^Tg\dd r\Big)^p\Big).
\]

\emph{Step 5: the exit probability.} On $\{\tau<T\}$, continuity of $e$ in $\HH$ gives $\hone{e(\tau)}=1$, hence $Z(T)\ge\mathcal Y(\tau)\ge1$. By Chebyshev's inequality, $\P(\tau<T)\le\P(Z(T)\ge1)\le\E Z(T)^q\le C_q\,\eps^q$.
\end{proof}

\subsection{Proof of Theorem~\ref{thm:LLN}}\label{sec:LLN-proof}
\begin{proof}[Proof of Theorem~\ref{thm:LLN}]
Let $Z_*=\Ec_T(w,e)=\sup_{t\le T}\mathcal Y(t)+\int_0^T\Dc\dd s$ be the unlocalised energy. On $\{\tau=T\}$ we have $Z_*=Z(T)$. Moreover $Z_*\le2(1+T)(\norm{X^\eps}_{\KT}^2+\norm{X^0}_{\KT}^2)$, so $\E Z_*^{2p}\le C_p$ uniformly in $\eps$ by Proposition~\ref{prop:moments}. By the Cauchy--Schwarz inequality and Lemma~\ref{lem:LLN-local}, the latter with exponent $p$ for $Z(T)$ and with $q=2p$ for the exit probability,
\[
  \E Z_*^p\le\E Z(T)^p+\E\big[Z_*^p\mathbf 1_{\{\tau<T\}}\big]\le C_p\eps^p+\big(\E Z_*^{2p}\big)^{1/2}\P(\tau<T)^{1/2}\le C_p\eps^p .
\]
Since $\norm{(w,e)}_{\KT}^2\le(1+T)Z_*$, we conclude that $\E\norm{X^\eps-X^0}_{\KT}^{2p}\le(1+T)^p\,\E Z_*^p\le C_p\eps^p$.
\end{proof}

\begin{remark}\label{rem:LLN}
\begin{enumerate}[label=(\alph*),leftmargin=*]
\item \emph{Dependence on the initial datum.} The constant in Lemma~\ref{lem:LLN-energy} does not depend on the initial datum; this is why the same computation can be reused in the proof of Lemma~\ref{lem:lipschitz}, with $X^0$ replaced by another deterministic solution. The constant in Theorem~\ref{thm:LLN} depends on the initial datum only through $\int_0^Tg\dd t$ and $\abs{u_0}^2+\hone{d_0}^2$, so by Lemma~\ref{lem:unif-bounds} it is uniform over bounded sets of initial data.
\item \emph{Optimality of the rate.} The rate $\eps^p$ is optimal. Indeed, by Theorem~\ref{thm:CLT} and Minkowski's inequality,
\[
  \Big|\eps^{-1/2}\big(\E\norm{X^\eps-X^0}_{\KT}^{2p}\big)^{1/2p}-\big(\E\norm{X^1}_{\KT}^{2p}\big)^{1/2p}\Big|\le C_p\sqrt\eps,
\]
so that $\E\norm{X^\eps-X^0}_{\KT}^{2p}=\eps^p\big(\E\norm{X^1}_{\KT}^{2p}+o(1)\big)$. The limit is positive unless $X^1=0$ a.s. By \eqref{eq:lin-u}--\eqref{eq:lin-d}, $X^1=0$ a.s.\ forces the stochastic integrals $\int_0^\cdot\sigma(u^0)\dd W_1$ and $\int_0^\cdot G(d^0)\dd W_2$ to vanish, hence also their quadratic variations. Thus $X^1=0$ a.s.\ if and only if $\sigma(u^0(t))=0$ and $G(d^0(t))=0$ for a.e.\ $t\in[0,T]$, that is, only when the noise vanishes along the deterministic flow.
\end{enumerate}
\end{remark}
\section{The linearised system and the skeleton equation}\label{sec:linear}

This section studies the skeleton map $\Gam$ and the linearised SPDE \eqref{eq:lin-u}--\eqref{eq:lin-d}. We show that $\Gam$ is a bounded, causal and compact linear operator (Propositions~\ref{prop:skeleton} and \ref{prop:skel-cont}) and prove Theorem~\ref{thm:linear}.

The main point in Theorem~\ref{thm:linear}(ii) is the duality relation
\[
  \E\big[\ell(X^1)\,\mathbb W(q)\big]=\ell\big(\Gam(q)\big),\qquad\ell\in\KT^*,\ q\in L^2(0,T;\mathcal U),
\]
where $\mathbb W(q)$ is the Wiener integral of $q$. A direct proof would rely on a variation-of-constants formula 
$$
X^1(t)=\int_0^t\mathcal T(t,s)\Sigma(s)\dd W(s),
$$
 and hence on the construction of an evolution family $\mathcal T(t,s)$ for the non-autonomous operator $\Lin(t)$, whose coefficients are only integrable in time. We avoid this by proving the relation for Galerkin approximations, for which such a formula is elementary, and passing to the limit.

Throughout this section we write, for $t\in[0,T]$ and $q=(\theta,\rho)\in\mathcal U$,
\[
  \Sigma(t)q=\big(\sigma(u^0(t))\theta,\ G(d^0(t))\rho\big)\in H\times\HH,
\]
so that the skeleton system reads $\partial_t\Gam(q)+\Lin(t)\Gam(q)=\Sigma q$ and the linearised SPDE reads $\dd X^1+\Lin(t)X^1\dd t=\Sigma\dd W$. By \ref{A:sigma} and Lemma~\ref{lem:G}(i), the map $t\mapsto\Sigma(t)$ is continuous from $[0,T]$ into the Hilbert--Schmidt operators $\LHS(\mathcal U,H\times\HH)$, and
\begin{equation}\label{eq:kappa}
  \Sigma_*:=\sup_{t\in[0,T]}\Big(\abs{\sigma(u^0(t))}_{\LHS}^2+\abs{G(d^0(t))}_{\LHSd}^2\Big)\le K\Big(1+\sup_t\abs{u^0(t)}^2\Big)+C_h\sup_t\hone{d^0(t)}^2<\infty .
\end{equation}

\subsection{Galerkin projections and the basic energy estimate}\label{sec:lin-prelim}
\emph{Galerkin projections.} Let $\{\varphi_j\}\subset D(\Ast)$ and $\{\varsigma_j\}\subset\DA$ be orthonormal bases of $H$ and $\LL$ consisting of eigenfunctions of the Stokes operator and of the Neumann Laplacian. Since $\Ocal$ is smooth, they belong to $C^\infty(\overline\Ocal)$. Let $H_n=\operatorname{span}\{\varphi_1,\dots,\varphi_n\}$ and $L_n=\operatorname{span}\{\varsigma_1,\dots,\varsigma_n\}$, let $P_n$ and $\widetilde P_n$ be the orthogonal projections of $H$ onto $H_n$ and of $\LL$ onto $L_n$, and let $\Pi_n=P_n\oplus\widetilde P_n$. We use the following standard facts; see e.g.\ \cite{Lions69,Temam}.
\begin{enumerate}[label=\textup{(P\arabic*)},leftmargin=*]
\item\label{P1} $P_n$ is also the orthogonal projection onto $H_n$ in $V$, with the inner product $\scal{\nabla\cdot}{\nabla\cdot}$. It is a contraction on $H$ and on $V$, it extends by duality to a contraction on $V'$, and $P_nv\to v$ in $H$ and in $V$.
\item\label{P2} $\widetilde P_n$ commutes with $\An$; in particular $(I+\An)L_n\subset L_n$. It is the orthogonal projection onto $L_n$ in $\LL$, $\HH$ and $\DA$, a contraction on each of these spaces, and $\widetilde P_ne\to e$ in each of them.
\end{enumerate}

\emph{The energy estimate.} The following lemma collects the energy computation that is used repeatedly below. For $X=(U,D)$ we write
\[
  \mathcal Y=\abs U^2+\hone D^2,\qquad \Dc=\norm U^2+\abs{\nabla D}^2+\abs{\An D}^2 .
\]

\begin{lemma}\label{lem:lin-energy}
Let $F=(F_1,F_2)\in L^2(0,T;H\times\HH)$, and let $X=(U,D)\in\KT$ with $X(0)=0$ satisfy
\[
  \partial_tU+\Lin_1(t)X=F_1\ \text{in }V',\qquad\partial_tD+\Lin_2(t)X=F_2\ \text{in }\LL,\qquad\text{for a.e.\ }t\in[0,T].
\]
\begin{enumerate}[label=\textup{(\roman*)}]
\item For all $t\in[0,T]$,
\begin{equation}\label{eq:lin-energy}
  \mathcal Y(t)+\int_0^t\Dc\dd s\le2C_*\int_0^tg\,\mathcal Y\dd s+2\int_0^t\Big(\scal{F_1}{U}+\sone{F_2}{D}\Big)\dd s .
\end{equation}
\item If there are $\beta,c\ge0$ such that the last integral in \eqref{eq:lin-energy} is bounded by $\beta+c\int_0^tg\,\mathcal Y\dd s$ for all $t\in[0,T]$, then
\begin{equation}\label{eq:lin-gronwall}
  \Ec_T(X)\le2\beta\exp\Big((2C_*+c)\int_0^Tg\dd t\Big).
\end{equation}
\end{enumerate}
The same conclusions hold for every absolutely continuous $X_n\colon[0,T]\to H_n\times L_n$ with $X_n(0)=0$ that solves the Galerkin system $\partial_tX_n+\Pi_n\Lin(t)X_n=\Pi_nF$.
\end{lemma}

\begin{proof}
(i) By Lemma~\ref{lem:coercive}(ii), $F_1-\Lin_1X\in L^2(0,T;V')$ and $F_2-\Lin_2X\in L^2(0,T;\LL)$, so Lemma~\ref{lem:Ito} applies without noise:
\[
  \mathcal Y(t)=2\int_0^t\Big(\dual{F_1-\Lin_1X}{U}+\scal{F_2-\Lin_2X}{(I+\An)D}\Big)\dd s .
\]
Since $F_2\in\HH$ and $D\in\DA$ for a.e.\ $s$, Green's formula gives $\scal{F_2}{(I+\An)D}=\sone{F_2}{D}$, and Lemma~\ref{lem:coercive}(i) gives $-2\big(\dual{\Lin_1X}{U}+\scal{\Lin_2X}{(I+\An)D}\big)\le-\Dc+2C_*g\mathcal Y$. This proves \eqref{eq:lin-energy}.

(ii) Let $E(t)=\mathcal Y(t)+\int_0^t\Dc\dd s$. By (i) and the assumption, $E(t)\le\beta+(2C_*+c)\int_0^tg\,E\dd s$, so Gronwall's lemma gives $E(t)\le\beta\exp((2C_*+c)\int_0^Tg)$ for all $t$. Since $\Ec_T(X)\le\sup_tE(t)+E(T)$, \eqref{eq:lin-gronwall} follows.

For the Galerkin system, $\frac{\dd}{\dd t}\abs{U_n}^2=2\dual{P_n(F_1-\Lin_1X_n)}{U_n}=2\dual{F_1-\Lin_1X_n}{U_n}$ because $P_nU_n=U_n$, and $\frac{\dd}{\dd t}\hone{D_n}^2=2\scal{\widetilde P_n(F_2-\Lin_2X_n)}{(I+\An)D_n}=2\scal{F_2-\Lin_2X_n}{(I+\An)D_n}$ because $(I+\An)D_n\in L_n$ by \ref{P2}. The rest of the proof is identical.
\end{proof}

\subsection{The skeleton map}\label{sec:skeleton}

\begin{proposition}\label{prop:skeleton}
For every $q=(\theta,\rho)\in L^2(0,T;\mathcal U)$ the system \eqref{eq:skel-u}--\eqref{eq:skel-d} has a unique solution $\Gam(q)=(U^q,D^q)\in\KT$. The map $\Gam\colon L^2(0,T;\mathcal U)\to\KT$ is linear and
\begin{equation}\label{eq:skel-bound}
  \Ec_T(\Gam(q))\le C_\Gam\norm q_{L^2(0,T;\mathcal U)}^2,\qquad C_\Gam=2\exp\Big(C\int_0^Tg\dd t\Big),
\end{equation}
where $C$ depends only on $\Ocal$, $\Ft$, $K$, $\mathsf h$ and $\sup_t(\abs{u^0(t)}+\hone{d^0(t)})$. Moreover $\Gam$ is causal: the restriction of $\Gam(q)$ to $[0,t]$ depends only on $q|_{[0,t]}$.
\end{proposition}

\begin{proof}
\emph{Step 1: a priori estimate.} Let $X=(U,D)\in\KT$ solve \eqref{eq:skel-u}--\eqref{eq:skel-d}, that is, the system of Lemma~\ref{lem:lin-energy} with $F=\Sigma q$. By \eqref{eq:kappa}, Young's inequality and $g\ge1$,
\[
  2\scal{\sigma(u^0)\theta}{U}+2\sone{G(d^0)\rho}{D}\le\Sigma_*\,\mathcal Y+\abs q_{\mathcal U}^2\le\Sigma_*\,g\,\mathcal Y+\abs q_{\mathcal U}^2 .
\]
Lemma~\ref{lem:lin-energy}(ii) with $\beta=\norm q^2_{L^2(0,T;\mathcal U)}$ and $c=\Sigma_*$ gives \eqref{eq:skel-bound} with $C=2C_*+\Sigma_*$, which depends only on the quantities listed, by \eqref{eq:kappa}.

\emph{Step 2: uniqueness, linearity and causality.} The system is linear in $(X,q)$. The difference of two solutions with the same $q$ solves the system with $q=0$, hence vanishes by \eqref{eq:skel-bound}. Linearity of $\Gam$ follows. If $q=q'$ on $[0,t]$, then $\Gam(q)$ and $\Gam(q')$ solve the same system on $[0,t]$, and the a priori estimate on $[0,t]$ shows that they coincide there.

\emph{Step 3: Galerkin approximation.} For $x\in H_n\times L_n$, Lemma~\ref{lem:coercive}(ii) applied to the constant function $x$ shows that $t\mapsto\Lin(t)x$ belongs to $L^2(0,T;V'\times\LL)$. Hence the Galerkin system
\begin{equation}\label{eq:skel-galerkin}
  \partial_tX_n+\Pi_n\Lin(t)X_n=\Pi_n\Sigma q,\qquad X_n(0)=0,
\end{equation}
is a linear system of ordinary differential equations in $H_n\times L_n$ whose coefficients and forcing are in $L^2(0,T)$. It has a unique absolutely continuous solution, given by
\[
  X_n(t)=\Gam_n(q)(t):=\int_0^t\mathcal T_n(t,s)\Pi_n\Sigma(s)q(s)\dd s,
\]
where $\mathcal T_n(t,s)$ is the fundamental matrix of $\dot x=-\Pi_n\Lin(t)x$, that is, $\partial_t\mathcal T_n(t,s)=-\Pi_n\Lin(t)\mathcal T_n(t,s)$ and $\mathcal T_n(s,s)=I$. By the Galerkin version of Lemma~\ref{lem:lin-energy} and Step~1, $\Ec_T(X_n)\le C_\Gam\norm q^2$ for every $n$. By \ref{P1}, \ref{P2} and Lemma~\ref{lem:coercive}(ii), $(\partial_tX_n)$ is bounded in $L^2(0,T;V'\times\LL)$.

\emph{Step 4: passage to the limit.} By Lemma~\ref{lem:AL}, applied in the triples $V\subset H\subset V'$ and $\DA\subset\HH\subset\LL$, a subsequence of $(X_n)$ converges weakly in $L^2(0,T;V\times\DA)$, weakly-$*$ in $L^\infty(0,T;H\times\HH)$ and strongly in $L^2(0,T;H\times\HH)$ to some $X=(U,D)$. To identify the limit we test with the eigenfunctions. Using $b(u,v,w)=-b(u,w,v)$, \eqref{eq:M}, Lemma~\ref{lem:f}(iii) and $\varphi_j,\varsigma_j\in W^{1,\infty}$, we find for $(U',D')\in V\times\DA$
\begin{equation}\label{eq:test}
  \abs{\dual{\Lin_1(t)(U',D')}{\varphi_j}}+\abs{\scal{\Lin_2(t)(U',D')}{\varsigma_j}}\le C_j\Big(1+\abs{u^0(t)}+\hone{d^0(t)}^{2N}\Big)\Big(\abs{U'}+\hone{D'}\Big),
\end{equation}
for instance
\begin{align*}
  \dual{B(U',u^0)}{\varphi_j}&=-b(U',\varphi_j,u^0)\le\norm{\nabla\varphi_j}_{L^\infty}\abs{U'}\abs{u^0},\\
  \scal{Df(d^0)D'}{\varsigma_j}&\le C\norm{\varsigma_j}_{L^\infty}\norm{1+\abs{d^0}^{2N}}_{L^2}\abs{D'}.
\end{align*} Since $\sup_t(\abs{u^0}+\hone{d^0})<\infty$, the maps $X\mapsto\int_0^T\dual{\Lin_1X}{\varphi_j}\psi\dd t$ and $X\mapsto\int_0^T\scal{\Lin_2X}{\varsigma_j}\psi\dd t$, for $\psi\in C^1_c([0,T))$, are continuous linear functionals on $L^2(0,T;H\times\HH)$. Testing \eqref{eq:skel-galerkin} with $\varphi_j\psi$ and $\varsigma_j\psi$ for $n\ge j$, integrating by parts in time and letting $n\to\infty$, we see that $X$ solves \eqref{eq:skel-u}--\eqref{eq:skel-d} in the sense of distributions on $[0,T)$ with $X(0)=0$, tested against every $\varphi_j$ and $\varsigma_j$. Since $\{\varphi_j\}$ spans a dense subspace of $V$ and $\{\varsigma_j\}$ of $\LL$, $X$ solves the system in $V'\times\LL$. By Lemma~\ref{lem:coercive}(ii), $\partial_tX\in L^2(0,T;V'\times\LL)$, and Lemma~\ref{lem:Ito} shows that $X\in\KT$.

By uniqueness, every subsequence of $(\Gam_n(q))$ has a further subsequence converging to $\Gam(q)$ in the above sense. Hence
\begin{equation}\label{eq:galerkin-skel-conv}
  \Gam_n(q)\rightharpoonup\Gam(q)\quad\text{weakly in }L^2(0,T;V\times\DA)\qquad\text{for every }q\in L^2(0,T;\mathcal U),
\end{equation}
along the whole sequence.
\end{proof}

\begin{proposition}\label{prop:skel-cont}
Let $\Lambda<\infty$, and let $\Sb_\Lambda=\{q\in L^2(0,T;\mathcal U):\norm q^2\le \Lambda\}$ be endowed with the weak topology. If $q_n\to q$ in $\Sb_\Lambda$, then $\Gam(q_n)\to\Gam(q)$ in $\KT$. In particular $\Gam(\Sb_\Lambda)$ is a compact subset of $\KT$, and $\Gam\colon L^2(0,T;\mathcal U)\to\KT$ is a compact operator.
\end{proposition}

\begin{proof}
By linearity, $Z_n=(Z_n^u,Z_n^d):=\Gam(q_n)-\Gam(q)=\Gam(p_n)$, where $p_n=(\theta_n',\rho_n')=q_n-q$ satisfies $p_n\rightharpoonup0$ weakly in $L^2(0,T;\mathcal U)$ and $\norm{p_n}^2\le4\Lambda$.

\emph{Step 1: $Z_n\to0$ in $L^2(0,T;H\times\HH)$.} By \eqref{eq:skel-bound}, $(Z_n)$ is bounded in $\KT$, and by Lemma~\ref{lem:coercive}(ii) and the equations, $(\partial_tZ_n)$ is bounded in $L^2(0,T;V'\times\LL)$. By Lemma~\ref{lem:AL}, every subsequence has a further subsequence converging strongly in $L^2(0,T;H\times\HH)$ and weakly in $L^2(0,T;V\times\DA)$ to some $Z$. Multiplication by $\Sigma$ is a bounded linear operator from $L^2(0,T;\mathcal U)$ to $L^2(0,T;H\times\HH)$, by \eqref{eq:kappa}; hence it is weakly continuous, and the forcing terms $\Sigma p_n$ converge weakly to $0$. Passing to the limit as in Step~4 of the proof of Proposition~\ref{prop:skeleton}, $Z$ solves the skeleton system with zero control, so $Z=0$. Therefore the whole sequence $(Z_n)$ converges to $0$ in $L^2(0,T;H\times\HH)$.

\emph{Step 2: convergence in $\KT$.} By the Cauchy--Schwarz inequality, for all $t\in[0,T]$,
\[
  2\Big|\int_0^t\Big(\scal{\sigma(u^0)\theta_n'}{Z_n^u}+\sone{G(d^0)\rho_n'}{Z_n^d}\Big)\dd s\Big|\le2\,\Sigma_*^{1/2}\norm{p_n}\norm{Z_n}_{L^2(0,T;H\times\HH)}=:\beta_n ,
\]
and $\beta_n\to0$ by Step~1. Lemma~\ref{lem:lin-energy}(ii) with $\beta=\beta_n$ and $c=0$ gives $\Ec_T(Z_n)\le2\beta_n\exp(2C_*\int_0^Tg)\to0$.

Finally, $\Sb_\Lambda$ is compact and metrisable in the weak topology, so its continuous image $\Gam(\Sb_\Lambda)$ is compact in $\KT$. Since every bounded set of $L^2(0,T;\mathcal U)$ is contained in some $\Sb_\Lambda$, $\Gam$ is compact.
\end{proof}

\subsection{Well-posedness of the linearised SPDE}\label{sec:lin-wp}

\begin{proof}[Proof of Theorem~\ref{thm:linear}\textup{(i)}]
We follow the proof of Proposition~\ref{prop:skeleton}, with the forcing $\Sigma q$ replaced by the additive noise $\Sigma\dd W$.

\emph{Step 1: a priori estimate.} Let $X=(U,D)$ be a solution of \eqref{eq:lin-u}--\eqref{eq:lin-d} with paths in $\KT$. By Lemma~\ref{lem:Ito}, Lemma~\ref{lem:coercive}(i) and \eqref{eq:kappa},
\begin{equation}\label{eq:lin-ito}
  \mathcal Y(t)+\int_0^t\Dc\dd s\le2C_*\int_0^tg\,\mathcal Y\dd s+\Sigma_*T+\mathcal M(t),
\end{equation}
where $\mathcal M(t)=2\int_0^t\scal{U}{\sigma(u^0)\dd W_1}+2\int_0^t\sone{G(d^0)\dd W_2}{D}$ is a continuous local martingale with $\dd\langle\mathcal M\rangle_t\le4\Sigma_*\mathcal Y(t)\dd t$. For $R>0$ let $\tau_R=\inf\{t:\Ec_t(X)\ge R\}\wedge T$ and $Z_R(t)=\Ec_{t\wedge\tau_R}(X)\le R$. Evaluating \eqref{eq:lin-ito} at $s\wedge\tau_R$ and taking the supremum over $s\le t$, as in Step~1 of the proof of Lemma~\ref{lem:LLN-local},
\[
  Z_R(t)\le4C_*\int_0^tg\,Z_R\dd s+2\Sigma_*T+2\sup_{s\le t}\abs{\mathcal M(s\wedge\tau_R)} .
\]
By the Burkholder--Davis--Gundy and Young inequalities, for every $\gamma>0$,
\[
  \E\sup_{s\le t}\abs{\mathcal M(s\wedge\tau_R)}^p\le C_p\,\E\Big(4\Sigma_*\int_0^{t\wedge\tau_R}\mathcal Y\dd s\Big)^{p/2}\le C_p(\Sigma_*T)^{p/2}\,\E Z_R(t)^{p/2}\le\gamma\,\E Z_R(t)^p+C_{p,\gamma}.
\]
Since $\E Z_R(t)^p\le R^p$, we may absorb the term with $\gamma$ small and apply Gronwall's lemma as in Step~4 of the proof of Lemma~\ref{lem:LLN-local}. This gives $\E Z_R(T)^p\le C_p$ with $C_p$ independent of $R$. Since $X$ has paths in $\KT$, $\tau_R=T$ for $R>\Ec_T(X)$, and Fatou's lemma yields
\begin{equation}\label{eq:lin-moment}
  \E\,\Ec_T(X)^p\le C_p,\qquad p\ge1,
\end{equation}
where $C_p$ depends only on $p$, $T$, $C_*$, $\Sigma_*$ and $\int_0^Tg\dd t$. Since $\norm X_{\KT}^2\le(1+T)\Ec_T(X)$, this is the moment bound of Theorem~\ref{thm:linear}(i).

\emph{Step 2: uniqueness.} The difference $\bar X$ of two solutions satisfies, for a.e.\ $\omega$, the deterministic system of Lemma~\ref{lem:lin-energy} with $F=0$, so $\bar X=0$ by \eqref{eq:lin-gronwall} with $\beta=c=0$.

\emph{Step 3: Galerkin approximation.} With $\mathcal T_n$ as in Step~3 of the proof of Proposition~\ref{prop:skeleton}, the Galerkin system $\dd X_n+\Pi_n\Lin(t)X_n\dd t=\Pi_n\Sigma\dd W$, $X_n(0)=0$, is a linear stochastic differential equation in $H_n\times L_n$ with additive noise. It has the unique solution
\begin{equation}\label{eq:galerkin-repr}
  X_n(t)=\int_0^t\mathcal T_n(t,s)\Sigma_n(s)\dd W(s),\qquad \Sigma_n=\Pi_n\Sigma,
\end{equation}
as follows from It\^o's product rule applied to $\mathcal T_n(0,t)X_n(t)$. The process $X_n$ is Gaussian with continuous paths in a finite-dimensional space, so all its moments are finite. By \ref{P1} and \ref{P2}, $\abs{\Sigma_n(t)}_{\LHS(\mathcal U,H\times\HH)}^2\le\Sigma_*$, and the computation of Step~1, done as in the Galerkin part of Lemma~\ref{lem:lin-energy}, gives
\begin{equation}\label{eq:galerkin-moment}
  \E\,\Ec_T(X_n)^p\le C_p\qquad\text{for all }n\text{ and }p\ge1 .
\end{equation}

\emph{Step 4: passage to the limit.} By \eqref{eq:galerkin-moment}, $(X_n)$ is bounded in $L^2(\Omega\times[0,T];V\times\DA)$, so a subsequence converges weakly in this space to a limit $X^1$. The limit is progressively measurable, since progressively measurable processes form a closed, hence weakly closed, subspace. The functional $Z\mapsto\E\operatorname*{ess\,sup}_{t}\big(\abs{Z_u(t)}^2+\hone{Z_d(t)}^2\big)$ is convex and lower semicontinuous for the strong topology of $L^2(\Omega\times[0,T];H\times\HH)$, hence weakly lower semicontinuous. By \eqref{eq:galerkin-moment} it follows that $X^1\in L^\infty(0,T;H\times\HH)$ almost surely.

For $n\ge j$ we have $P_n\varphi_j=\varphi_j$ and $\widetilde P_n\varsigma_j=\varsigma_j$, so the Galerkin equations tested with $\varphi_j$ and $\varsigma_j$ read
\begin{align*}
  \scal{U_n(t)}{\varphi_j}+\int_0^t\dual{\Lin_1X_n}{\varphi_j}\dd s&=\int_0^t\scal{\varphi_j}{\sigma(u^0)\dd W_1},\\
  \scal{D_n(t)}{\varsigma_j}+\int_0^t\scal{\Lin_2X_n}{\varsigma_j}\dd s&=\int_0^t\scal{\varsigma_j}{G(d^0)\dd W_2}.
\end{align*}
The right-hand sides do not depend on $n$. Multiplying by $\psi\in L^\infty(\Omega\times[0,T])$ and integrating over $\Omega\times[0,T]$, the left-hand sides become continuous linear functionals of $X_n$ in $L^2(\Omega\times[0,T];H\times\HH)$, by \eqref{eq:test}. Passing to the limit and using that $\psi$ is arbitrary, we find that $X^1$ satisfies the tested equations for $\dd t\otimes\P$-a.e.\ $(t,\omega)$, for every $j$. By Lemma~\ref{lem:coercive}(ii), applied pathwise, the drift $\Lin X^1$ lies a.s.\ in $L^2(0,T;V'\times\LL)$. Since $\{\varphi_j\}$ spans a dense subspace of $V$ and $\{\varsigma_j\}$ of $\LL$, the integral forms of \eqref{eq:lin-u}--\eqref{eq:lin-d} hold in $V'\times\LL$ for $\dd t\otimes\P$-a.e.\ $(t,\omega)$. By the It\^o formula in the Gelfand triples (Lemma~\ref{lem:Ito}; in the form of \cite{Krylov13} and \cite[Section~4.2]{LiuRockner15}, which requires the equation only for a.e.\ $t$), $X^1$ has a version with continuous paths in $H\times\HH$ that satisfies the equations for all $t\in[0,T]$, a.s. Thus $X^1$ is a solution with paths in $\KT$, and Step~1 applies to it.

Finally, by Step~2 every subsequence of $(X_n)$ has a further subsequence converging weakly to $X^1$, so
\begin{equation}\label{eq:galerkin-conv}
  X_n\rightharpoonup X^1\quad\text{weakly in }L^2(\Omega\times[0,T];V\times\DA)
\end{equation}
along the whole sequence.
\end{proof}

\subsection{Gaussian structure and the link with the skeleton}\label{sec:gaussian}
For $k\in L^2(0,T;\mathcal U)$ let
\[
  \mathbb W(k)=\int_0^T\langle k(t),\dd W(t)\rangle_{\mathcal U}
\]
denote the Wiener integral with respect to $W$, which is well defined since $\mathcal U$ is the reproducing kernel Hilbert space of $W$ (Section~\ref{sec:noise}), and let $\mathcal H_1=\{\mathbb W(k):k\in L^2(0,T;\mathcal U)\}$ be the first Wiener chaos of $W$; see e.g.\ \cite{Nualart06}.

\begin{lemma}\label{lem:chaos}
The following hold.
\begin{enumerate}[label=\textup{(\roman*)}]
\item The map $k\mapsto\mathbb W(k)$ is a linear isometry from $L^2(0,T;\mathcal U)$ onto $\mathcal H_1\subset L^2(\Omega)$. In particular $\mathcal H_1$ is a closed subspace of $L^2(\Omega)$, and each of its elements is a centred Gaussian random variable.
\item $\mathcal H_1$ is closed under almost sure convergence. More precisely, if $\mathbb W(k_m)\to\vartheta$ almost surely, then $\vartheta\in\mathcal H_1$ and $\mathbb W(k_m)\to\vartheta$ in $L^2(\Omega)$.
\end{enumerate}
\end{lemma}

\begin{proof}
(i) is the It\^o isometry. (ii) Let $\vartheta_m=\mathbb W(k_m)$ and $v_m=\norm{k_m}^2$. By dominated convergence, $\exp(-\frac12v_ms^2)=\E e^{\mathrm is\vartheta_m}\to\E e^{\mathrm is\vartheta}$ for all $s\in\R$. If $(v_m)$ were unbounded, the limit would vanish along a subsequence for every $s\ne0$, which contradicts its continuity at $s=0$, where it equals $1$. Hence $(v_m)$ is bounded, and, considering $s=1$, convergent. The fourth moments $\E\vartheta_m^4=3v_m^2$ are then bounded, so $(\vartheta_m^2)$ is uniformly integrable and $\vartheta_m\to\vartheta$ in $L^2(\Omega)$ by Vitali's theorem. Since $\mathcal H_1$ is closed, $\vartheta\in\mathcal H_1$.
\end{proof}

\begin{lemma}\label{lem:dual}
Let $\KT^0=\{z\in\KT:z(0)=0\}$ and let $\mathbb F\subset\KT^*$ be the space of functionals
\[
  \ell_\varpi(v,e)=\int_0^T\Big(\scal{\varpi_1(t)}{v(t)}+\sone{\varpi_2(t)}{e(t)}\Big)\dd t,\qquad \varpi_1\in L^\infty(0,T;H),\ \varpi_2\in L^\infty(0,T;\HH).
\]
Every $\ell_\varpi\in\mathbb F$ is continuous on $L^2(0,T;H\times\HH)$. For every $\ell\in\KT^*$ there is a sequence $(\ell_m)\subset\mathbb F$ with $\ell_m(z)\to\ell(z)$ for every $z\in\KT^0$.
\end{lemma}

\begin{proof}
The first statement is the Cauchy--Schwarz inequality. For the second, the map $(v,e)\mapsto(v,v,e,e)$ embeds $\KT$ isomorphically onto a closed subspace of
\[
  C([0,T];H)\times L^2(0,T;V)\times C([0,T];\HH)\times L^2(0,T;\DA).
\]
By the Hahn--Banach theorem, $\ell=\ell^1+\ell^2+\ell^3+\ell^4$, where each $\ell^i$ is a continuous linear functional on the $i$-th factor.

\emph{The $L^2$ parts.} The functional $\ell^2$ has the form $v\mapsto\int_0^T\dual{\varpi(t)}{v(t)}\dd t$ with $\varpi\in L^2(0,T;V')$. Approximating $\varpi$ in $L^2(0,T;V')$ by elements $\varpi_m$ of $L^\infty(0,T;H)$, which is possible since $H$ is dense in $V'$, we obtain functionals $v\mapsto\int_0^T\scal{\varpi_m}{v}\dd t$ in $\mathbb F$ that converge to $\ell^2$ in norm. Since $I+\An$ is an isomorphism from $\DA$ onto $\LL$, the functional $\ell^4$ has the form $e\mapsto\int_0^T\scal{\varpi(t)}{(I+\An)e(t)}\dd t$ with $\varpi\in L^2(0,T;\LL)$. We approximate $\varpi$ in $L^2(0,T;\LL)$ by elements $\varpi_m$ of $L^\infty(0,T;\HH)$ and use $\scal{\varpi_m}{(I+\An)e}=\sone{\varpi_m}{e}$.

\emph{The $C$ parts.} By Singer's representation theorem (see e.g.\ \cite{Dinculeanu67}), $\ell^1(v)=\int_{[0,T]}\scal{v(t)}{\mu(\dd t)}$ for a regular $H$-valued measure $\mu$ of bounded variation; since $H$ is a Hilbert space, $\mu=h\abs\mu$ with $h\in L^\infty(\abs\mu;H)$. Let $v\in C([0,T];H)$ with $v(0)=0$, extend $v$ by $0$ to negative times, and let $v_m(t)=m\int_{t-1/m}^tv(s)\dd s$. Then $v_m\to v$ uniformly on $[0,T]$, and by Fubini's theorem
\[
  \ell^1(v)=\lim_{m\to\infty}\int_{[0,T]}\scal{v_m(t)}{\mu(\dd t)}=\lim_{m\to\infty}\int_0^T\scal{v(s)}{\varpi_m(s)}\dd s,
\]
where $\varpi_m(s)=m\,\mu\big([s,s+\tfrac1m]\cap[0,T]\big)$ satisfies $\abs{\varpi_m(s)}\le m\abs\mu([0,T])$, so that $\varpi_m\in L^\infty(0,T;H)$. The functional $\ell^3$ is treated in the same way, with an $\HH$-valued measure and the inner product of $\HH$.
\end{proof}

\begin{proof}[Proof of Theorem~\ref{thm:linear}\textup{(ii)}]
Let $X_n$ be the Galerkin approximations \eqref{eq:galerkin-repr}. For every $\ell\in\mathbb F$, the map $Z\mapsto\E[\xi\,\ell(Z)]$ is a continuous linear functional on $L^2(\Omega\times[0,T];H\times\HH)$ for every $\xi\in L^2(\Omega)$, by Lemma~\ref{lem:dual}. Hence \eqref{eq:galerkin-conv} implies
\begin{equation}\label{eq:ell-conv}
  \ell(X_n)\rightharpoonup\ell(X^1)\quad\text{weakly in }L^2(\Omega),\qquad\ell\in\mathbb F .
\end{equation}

\emph{Step 1: duality at the Galerkin level.} Let $\ell=\ell_\varpi\in\mathbb F$ with $\varpi=(\varpi_1,\varpi_2)$, and write $\langle\cdot,\cdot\rangle$ for the inner product of $H\times\HH$. Since $\Pi_n$ is the orthogonal projection onto $H_n\times L_n$ in $H\times\HH$ by \ref{P1} and \ref{P2}, we have $\langle\varpi(t),x\rangle=\langle\Pi_n\varpi(t),x\rangle$ for $x\in H_n\times L_n$. By \eqref{eq:galerkin-repr} and the stochastic Fubini theorem \cite{DPZ14}, whose hypotheses hold because the integrand is deterministic and bounded,
\[
  \ell(X_n)=\int_0^T\Big\langle\Pi_n\varpi(t),\int_0^t\mathcal T_n(t,s)\Sigma_n(s)\dd W(s)\Big\rangle\dd t=\mathbb W(k_n),
\]
with
\[
  k_n(s)=\int_s^T\Sigma_n(s)^*\mathcal T_n(t,s)^*\Pi_n\varpi(t)\dd t,
\]
where $\mathcal T_n(t,s)^*$ is the adjoint in $H_n\times L_n$ with the inner product of $H\times\HH$, and $\Sigma_n(s)^*\colon H_n\times L_n\to\mathcal U$ is the adjoint of $\Sigma_n(s)$. In particular $\ell(X_n)\in\mathcal H_1$. Let $q\in L^2(0,T;\mathcal U)$ and let $\Gam_n(q)$ be as in Step~3 of the proof of Proposition~\ref{prop:skeleton}. By the It\^o isometry and Fubini's theorem,
\begin{equation}\label{eq:duality-n}
  \E\big[\ell(X_n)\,\mathbb W(q)\big]=\langle k_n,q\rangle_{L^2(0,T;\mathcal U)}=\int_0^T\Big\langle\varpi(t),\int_0^t\mathcal T_n(t,s)\Sigma_n(s)q(s)\dd s\Big\rangle\dd t=\ell\big(\Gam_n(q)\big).
\end{equation}

\emph{Step 2: passage to the limit for $\ell\in\mathbb F$.} By \eqref{eq:ell-conv}, and since the closed subspace $\mathcal H_1$ is weakly closed, $\ell(X^1)\in\mathcal H_1$. By \eqref{eq:galerkin-skel-conv} and Lemma~\ref{lem:dual}, $\ell(\Gam_n(q))\to\ell(\Gam(q))$. Letting $n\to\infty$ in \eqref{eq:duality-n}, we obtain
\begin{equation}\label{eq:duality}
  \E\big[\ell(X^1)\,\mathbb W(q)\big]=\ell\big(\Gam(q)\big)\qquad\text{for all }q\in L^2(0,T;\mathcal U).
\end{equation}

\emph{Step 3: general $\ell\in\KT^*$.} Let $(\ell_m)\subset\mathbb F$ be as in Lemma~\ref{lem:dual}. Since $X^1\in\KT^0$ almost surely, $\ell_m(X^1)\to\ell(X^1)$ almost surely. By Step~2 and Lemma~\ref{lem:chaos}(ii), $\ell(X^1)\in\mathcal H_1$ and the convergence also holds in $L^2(\Omega)$. Since $\Gam(q)\in\KT^0$, also $\ell_m(\Gam(q))\to\ell(\Gam(q))$. Applying \eqref{eq:duality} to $\ell_m$ and letting $m\to\infty$, we see that \eqref{eq:duality} holds for every $\ell\in\KT^*$.

\emph{Step 4: conclusion.} Let $\ell\in\KT^*$ and write $\ell(X^1)=\mathbb W(\tilde k_\ell)$ with $\tilde k_\ell\in L^2(0,T;\mathcal U)$. By the It\^o isometry and \eqref{eq:duality},
\[
  \langle\tilde k_\ell,q\rangle_{L^2(0,T;\mathcal U)}=\E\big[\mathbb W(\tilde k_\ell)\,\mathbb W(q)\big]=\ell(\Gam(q))\qquad\text{for every }q\in L^2(0,T;\mathcal U),
\]
so $\tilde k_\ell=k_\ell$ is the Riesz representative of $\ell\circ\Gam$. In particular $\ell(X^1)$ is a centred Gaussian random variable with variance $\norm{k_\ell}^2=\norm{\ell\circ\Gam}^2$. Since $\KT$ is separable and every continuous linear functional of the $\KT$-valued random variable $X^1$ is centred Gaussian, $X^1$ is a centred Gaussian random variable in $\KT$.
\end{proof}

Theorem~\ref{thm:linear}(ii) identifies the law of $X^1$ completely. Let $\Gam^*\colon\KT^*\to L^2(0,T;\mathcal U)$ denote the adjoint of $\Gam$, after the Riesz identification of $L^2(0,T;\mathcal U)$ with its dual, so that $\Gam^*\ell=k_\ell$ is the function of Theorem~\ref{thm:linear}(ii). We use the definition of the Cameron--Martin space of a centred Gaussian measure $\mu$ on a separable Banach space $E$ from \cite[Section~2.4]{Bogachev98}:
\[
  H(\mu)=\{h\in E:\abs h_{H(\mu)}<\infty\},\qquad \abs h_{H(\mu)}=\sup\Big\{\ell(h):\ \ell\in E^*,\ \int_E\ell^2\dd\mu\le1\Big\}.
\]

\begin{corollary}\label{cor:gaussian}
Let $\mu=\Law(X^1)$, which is a centred Gaussian measure on $\KT$. Let $\mathscr H=\Gam(L^2(0,T;\mathcal U))$, with the norm $\norm z_{\mathscr H}=\inf\{\norm q:\Gam q=z\}$.
\begin{enumerate}[label=\textup{(\roman*)}]
\item \emph{Covariance.} $\E[\ell_1(X^1)\,\ell_2(X^1)]=\langle\Gam^*\ell_1,\Gam^*\ell_2\rangle_{L^2(0,T;\mathcal U)}$ for $\ell_1,\ell_2\in\KT^*$; that is, the covariance operator of $\mu$ is $\Gam\Gam^*\colon\KT^*\to\KT$.
\item \emph{Cameron--Martin space.} $H(\mu)=\mathscr H$ and $\abs z_{H(\mu)}=\norm z_{\mathscr H}$ for $z\in\mathscr H$. The space $(\mathscr H,\norm\cdot_{\mathscr H})$ is a Hilbert space, isometric to $(\ker\Gam)^\perp$ via $\Gam$, and continuously embedded in $\KT$. Unless $\mathscr H$ is finite-dimensional, it is not closed in $\KT$.
\item \emph{Support.} The topological support of $\mu$ is the closure of $\mathscr H$ in $\KT$:
\[
  \operatorname{supp}\Law(X^1)=\overline{\Gam\big(L^2(0,T;\mathcal U)\big)}^{\,\KT}.
\]
\item \emph{The MDP rate function.} The rate function of Theorem~\ref{thm:MDP} is the Cameron--Martin form of $\mu$:
\[
  I(z)=\tfrac12\norm z_{\mathscr H}^2\quad\text{for }z\in\mathscr H,\qquad I(z)=+\infty\quad\text{for }z\in\KT\setminus\mathscr H .
\]
\end{enumerate}
\end{corollary}

\begin{proof}
(i) follows from Theorem~\ref{thm:linear}(ii) and the It\^o isometry.

(ii) By (i), $\abs h_{H(\mu)}=\sup\{\ell(h):\norm{\Gam^*\ell}\le1\}$. Let $h=\Gam q$. Then $\ell(h)=\langle\Gam^*\ell,q\rangle\le\norm{\Gam^*\ell}\norm q$, so $\abs h_{H(\mu)}\le\norm h_{\mathscr H}$. Conversely, let $q_0$ be the preimage of $h$ of minimal norm, so that $q_0\in(\ker\Gam)^\perp$ and $\norm{q_0}=\norm h_{\mathscr H}$. Since $\ker\Gam=(\Gam^*(\KT^*))^\perp$, the range of $\Gam^*$ is dense in $(\ker\Gam)^\perp$. If $q_0\ne0$, choose $\ell_n$ with $\Gam^*\ell_n\to q_0/\norm{q_0}$. Then $\ell_n(h)=\langle\Gam^*\ell_n,q_0\rangle\to\norm{q_0}$, so $\abs h_{H(\mu)}\ge\norm{q_0}$. Now let $h\in\KT$ with $c=\abs h_{H(\mu)}<\infty$. The map $\Gam^*\ell\mapsto\ell(h)$ is well defined and bounded by $c\norm{\Gam^*\ell}$ on the range of $\Gam^*$, so by the Riesz theorem there is $k\in(\ker\Gam)^\perp$ with $\ell(h)=\langle\Gam^*\ell,k\rangle=\ell(\Gam k)$ for all $\ell\in\KT^*$. By the Hahn--Banach theorem $h=\Gam k\in\mathscr H$. Hence $H(\mu)=\mathscr H$ with equal norms. The map $\Gam\colon(\ker\Gam)^\perp\to\mathscr H$ is a bijective isometry, so $\mathscr H$ is a Hilbert space, and $\norm z_{\KT}\le\norm\Gam\norm z_{\mathscr H}$. Finally, if $\mathscr H$ were closed in $\KT$, the open mapping theorem would make $\Gam\colon(\ker\Gam)^\perp\to\mathscr H$ an isomorphism of Banach spaces. Since $\Gam$ is compact (Proposition~\ref{prop:skel-cont}), the unit ball of $(\ker\Gam)^\perp$ would then be compact, so $(\ker\Gam)^\perp$, and hence $\mathscr H$, would be finite-dimensional.

(iii) Since $\KT$ is a separable Banach space, $\mu$ is a Radon centred Gaussian measure, and its topological support is the closure of its Cameron--Martin space \cite[Theorem~3.6.1]{Bogachev98}. We conclude with (ii).

(iv) By definition, $I(z)=\inf\{\frac12\norm q^2:\Gam q=z\}$, which equals $\frac12\norm z_{\mathscr H}^2$ if $z\in\mathscr H$ and $+\infty$ otherwise.
\end{proof}
\section{Central limit theorem}\label{sec:CLT}

In this section we prove Theorem~\ref{thm:CLT} by estimating the remainder $\widehat X^\eps=(X^\eps-X^0)/\sqrt\eps-X^1$. It solves the linearised equation with the deterministic operator $\Lin$, forced by terms of order $\sqrt\eps$ whose moments are controlled by the law of large numbers. Since these terms are forcing terms rather than functions of the unknown, no localisation is needed.

Throughout this section, $X^1=(U,D)$ is the solution of \eqref{eq:lin-u}--\eqref{eq:lin-d} given by Theorem~\ref{thm:linear}, on the same stochastic basis and driven by the same Wiener process as $X^\eps$, and $(w,e)=X^\eps-X^0$ is as in Section~\ref{sec:LLN}.

\subsection{The equation for the remainder}\label{sec:CLT-eq}
Define the rescaled fluctuation and the remainder
\[
  (U^\eps,D^\eps)=\frac{(w,e)}{\sqrt\eps},\qquad \widehat X^\eps=(\widehat U,\widehat D)=(U^\eps-U,\,D^\eps-D),
\]
and the rescaled Taylor remainder of the potential (see Lemma~\ref{lem:fspace}(d))
\[
  \mathcal R^\eps=\frac{1}{\sqrt\eps}\,\mathcal R(d^0,e)=\frac{f(d^\eps)-f(d^0)}{\sqrt\eps}-Df(d^0)D^\eps .
\]
We divide \eqref{eq:w}--\eqref{eq:e} by $\sqrt\eps$. By bilinearity, $\eps^{-1/2}B(w,w)=\sqrt\eps\,B(U^\eps,U^\eps)$, and similarly for $M$ and $\Bt$. Moreover $\Lin_2^\flat(U^\eps,D^\eps)+Df(d^0)D^\eps=\Lin_2(U^\eps,D^\eps)$, and, since $G$ is linear, $G(d^\eps)=G(d^0)+\sqrt\eps\,G(D^\eps)$. Hence $(U^\eps,D^\eps)$ solves
\begin{align*}
  \dd U^\eps+\Lin_1(U^\eps,D^\eps)\dd t&=-\sqrt\eps\big[B(U^\eps,U^\eps)+M(D^\eps,D^\eps)\big]\dd t+\sigma(u^\eps)\dd W_1,\\
  \dd D^\eps+\Lin_2(U^\eps,D^\eps)\dd t&=\Big[-\sqrt\eps\,\Bt(U^\eps,D^\eps)-\mathcal R^\eps+\frac{\sqrt\eps}2G^2(d^\eps)\Big]\dd t+\big[G(d^0)+\sqrt\eps\,G(D^\eps)\big]\dd W_2 .
\end{align*}
Subtracting \eqref{eq:lin-u}--\eqref{eq:lin-d}, we find that $\widehat X^\eps$ satisfies $\widehat U(0)=0$, $\widehat D(0)=0$ and
\begin{align}
  \dd\widehat U+\Lin_1(\widehat U,\widehat D)\dd t&=-\sqrt\eps\big[B(U^\eps,U^\eps)+M(D^\eps,D^\eps)\big]\dd t+\big[\sigma(u^\eps)-\sigma(u^0)\big]\dd W_1, \label{eq:hatU}\\
  \dd\widehat D+\Lin_2(\widehat U,\widehat D)\dd t&=\Big[-\sqrt\eps\,\Bt(U^\eps,D^\eps)-\mathcal R^\eps+\frac{\sqrt\eps}2G^2(d^\eps)\Big]\dd t+\sqrt\eps\,G(D^\eps)\dd W_2. \label{eq:hatD}
\end{align}
Every term that is linear in $\widehat X^\eps$ is contained in $\Lin$, whose coefficients are deterministic, and the right-hand sides do not contain $\widehat X^\eps$. Moreover, the noise $G(d^0)\dd W_2$ of the linearised equation cancels exactly against the leading part of $G(d^\eps)\dd W_2$, and the It\^o correction appears only at order $\sqrt\eps$.

Since $X^\eps$, $X^0$ and $X^1$ have paths in $\KT$ and solve their respective equations, the drifts of \eqref{eq:hatU}--\eqref{eq:hatD} belong a.s.\ to $L^2(0,T;V')$ and $L^2(0,T;\LL)$, and the diffusion coefficients to $L^2(0,T;\LHS)$ and $L^2(0,T;\LHSd)$; see Section~\ref{sec:LLN-eq} and Lemma~\ref{lem:coercive}(ii). Hence Lemma~\ref{lem:Ito} applies to $\widehat X^\eps$.

\subsection{Bounds on the forcing terms}\label{sec:CLT-forcing}
Let
\[
  \Theta^\eps=\sup_{t\le T}\big(\abs{U^\eps}^2+\hone{D^\eps}^2\big)+\int_0^T\big(\norm{U^\eps}^2+\hda{D^\eps}^2\big)\dd t=\frac1\eps\norm{X^\eps-X^0}_{\KT}^2
\]
be the energy of the rescaled fluctuation.

\begin{lemma}\label{lem:CLT-forcing}
The following hold.
\begin{enumerate}[label=\textup{(\roman*)}]
\item For every $q\ge1$, $\sup_{\eps\in(0,1]}\E(\Theta^\eps)^q<\infty$.
\item For all $t\in[0,T]$,
\begin{equation}\label{eq:R-bound}
  \abs{\mathcal R^\eps(t)}\le C\sqrt\eps\Big(1+\hone{d^0(t)}^{2N-1}+\hone{e(t)}^{2N-1}\Big)\hone{D^\eps(t)}^2 .
\end{equation}
\item Let
\[
  \Phi^\eps=\abs{B(U^\eps,U^\eps)}_{V'}^2+\abs{M(D^\eps,D^\eps)}_{V'}^2+\abs{\Bt(U^\eps,D^\eps)}^2+\frac1\eps\abs{\mathcal R^\eps}^2+\abs{G^2(d^\eps)}^2 .
\]
Then, almost surely,
\begin{equation}\label{eq:Xi}
  \int_0^T\Phi^\eps\dd t\le\Xi^\eps:=C\Big(1+(\Theta^\eps)^{2N+1}+\norm{X^\eps}_{\KT}^2\Big),
\end{equation}
and $\sup_{\eps\in(0,1]}\E(\Xi^\eps)^q<\infty$ for every $q\ge1$.
\end{enumerate}
Here $C$ depends only on $\Ocal$, $\Ft$, $\mathsf h$, $T$ and $\sup_t\hone{d^0(t)}$.
\end{lemma}

\begin{proof}
(i) Since $\Theta^\eps=\eps^{-1}\norm{X^\eps-X^0}_{\KT}^2$, this is Theorem~\ref{thm:LLN}.

(ii) Lemma~\ref{lem:fspace}(d) with $d=d^0$ and $D=e=\sqrt\eps D^\eps$ gives
\[
  \abs{\mathcal R(d^0,e)}\le C\big(1+\hone{d^0}^{2N-1}+\hone e^{2N-1}\big)\,\eps\hone{D^\eps}^2 .
\]
Divide by $\sqrt\eps$.

(iii) We bound the terms of $\Phi^\eps$ in turn, writing $\Theta=\Theta^\eps$.

\emph{Quadratic terms.} By Lemma~\ref{lem:nonlinear}(ii)--(iv),
\[
  \abs{B(U^\eps,U^\eps)}_{V'}^2+\abs{M(D^\eps,D^\eps)}_{V'}^2+\abs{\Bt(U^\eps,D^\eps)}^2\le C\big(\abs{U^\eps}^2+\abs{\nabla D^\eps}^2\big)\big(\norm{U^\eps}^2+\hda{D^\eps}^2\big),
\]
whose integral over $[0,T]$ is at most $C\Theta^2$.

\emph{Taylor remainder.} By (ii) and $\sup_t\hone{e}^2=\eps\sup_t\hone{D^\eps}^2\le\Theta$,
\[
  \frac1\eps\int_0^T\abs{\mathcal R^\eps}^2\dd t\le CT\big(1+\sup_t\hone{d^0}^{4N-2}+\Theta^{2N-1}\big)\Theta^2 .
\]

\emph{It\^o correction.} By Lemma~\ref{lem:G}(i), $\int_0^T\abs{G^2(d^\eps)}^2\dd t\le C_hT\norm{X^\eps}_{\KT}^2$.

Since $\Theta^2\le1+\Theta^{2N+1}$, \eqref{eq:Xi} follows. The moment bound for $\Xi^\eps$ follows from (i) and Proposition~\ref{prop:moments}.
\end{proof}

\subsection{Proof of Theorem~\ref{thm:CLT}}\label{sec:CLT-proof}
Let
\[
  \widehat{\mathcal Y}=\abs{\widehat U}^2+\hone{\widehat D}^2,\qquad \widehat\Dc=\norm{\widehat U}^2+\abs{\nabla\widehat D}^2+\abs{\An\widehat D}^2,\qquad \widehat Z(t)=\Ec_t(\widehat X^\eps)=\sup_{s\le t}\widehat{\mathcal Y}(s)+\int_0^t\widehat\Dc\dd s .
\]

\begin{proof}[Proof of Theorem~\ref{thm:CLT}]
\emph{Step 1: It\^o's formula.} Applying Lemmas~\ref{lem:Ito} and \ref{lem:coercive}(i) to \eqref{eq:hatU}--\eqref{eq:hatD}, we get
\begin{equation}\label{eq:CLT-ito}
  \widehat{\mathcal Y}(t)+\int_0^t\widehat\Dc\dd s\le2C_*\int_0^tg\,\widehat{\mathcal Y}\dd s+\sum_{i=1}^5\int_0^tJ_i\dd s+\widehat{\mathcal M}(t),
\end{equation}
where
\begin{align*}
  J_1&=-2\sqrt\eps\dual{B(U^\eps,U^\eps)+M(D^\eps,D^\eps)}{\widehat U},&
  J_2&=-2\sqrt\eps\scal{\Bt(U^\eps,D^\eps)}{(I+\An)\widehat D},\\
  J_3&=-2\scal{\mathcal R^\eps}{(I+\An)\widehat D},&
  J_4&=\sqrt\eps\scal{G^2(d^\eps)}{(I+\An)\widehat D},\\
  J_5&=\abs{\sigma(u^\eps)-\sigma(u^0)}_{\LHS}^2+\eps\abs{G(D^\eps)}_{\LHSd}^2,
\end{align*}
and $\widehat{\mathcal M}$ is the continuous local martingale
\[
  \widehat{\mathcal M}(t)=2\int_0^t\scal{\widehat U}{\big(\sigma(u^\eps)-\sigma(u^0)\big)\dd W_1}+2\sqrt\eps\int_0^t\sone{G(D^\eps)\dd W_2}{\widehat D}.
\]

\emph{Step 2: the forcing terms.} By Young's inequality, $\abs{(I+\An)\widehat D}\le\abs{\widehat D}+\abs{\An\widehat D}$ and the definition of $\Phi^\eps$ in Lemma~\ref{lem:CLT-forcing}(iii),
\begin{align*}
  J_1&\le\tfrac18\norm{\widehat U}^2+C\eps\,\Phi^\eps,\qquad
  J_2\le\tfrac18\abs{\An\widehat D}^2+\abs{\widehat D}^2+C\eps\,\Phi^\eps,\\
  J_3&\le\tfrac18\abs{\An\widehat D}^2+\abs{\widehat D}^2+C\abs{\mathcal R^\eps}^2\le\tfrac18\abs{\An\widehat D}^2+\abs{\widehat D}^2+C\eps\,\Phi^\eps,\\
  J_4&\le\tfrac18\abs{\An\widehat D}^2+\abs{\widehat D}^2+C\eps\,\Phi^\eps .
\end{align*}
For the noise terms, \ref{A:sigma} and Lemma~\ref{lem:G}(i) give
\[
  J_5\le L\abs{u^\eps-u^0}^2+\eps C_h\hone{D^\eps}^2=\eps\big(L\abs{U^\eps}^2+C_h\hone{D^\eps}^2\big)\le C\eps\,\Theta^\eps .
\]
The terms with $\frac18$ add up to at most $\frac12\widehat\Dc$, and $3\abs{\widehat D}^2\le3g\widehat{\mathcal Y}$ since $g\ge1$. By Lemma~\ref{lem:CLT-forcing}(iii), for all $t\in[0,T]$,
\begin{equation}\label{eq:CLT-J}
  \sum_{i=1}^5\int_0^tJ_i\dd s\le\frac12\int_0^t\widehat\Dc\dd s+3\int_0^tg\,\widehat{\mathcal Y}\dd s+C\eps\big(\Xi^\eps+\Theta^\eps\big).
\end{equation}

\emph{Step 3: the martingale term.} By \ref{A:sigma} and Lemma~\ref{lem:G}(i),
\[
  \dd\langle\widehat{\mathcal M}\rangle_t\le4\Big(L\abs{u^\eps-u^0}^2\abs{\widehat U}^2+\eps C_h\hone{D^\eps}^2\hone{\widehat D}^2\Big)\dd t\le C\eps\,\Theta^\eps\,\widehat{\mathcal Y}(t)\dd t .
\]
By the Burkholder--Davis--Gundy and Young inequalities, for every $\gamma>0$ and $p\ge1$,
\[
  \E\sup_{s\le t}\abs{\widehat{\mathcal M}(s)}^p\le C_p\eps^{p/2}\,\E\Big[\widehat Z(t)^{p/2}(\Theta^\eps)^{p/2}\Big]\le\gamma\,\E\widehat Z(t)^p+C_{p,\gamma}\,\eps^p\,\E(\Theta^\eps)^p .
\]

\emph{Step 4: Gronwall's lemma.} Inserting \eqref{eq:CLT-J} into \eqref{eq:CLT-ito}, moving $\frac12\int_0^t\widehat\Dc$ to the left-hand side, and taking the supremum over $s\le t$ as in Step~1 of the proof of Lemma~\ref{lem:LLN-local}, we obtain
\[
  \widehat Z(t)\le C\int_0^tg\,\widehat Z\dd s+C\eps\big(\Xi^\eps+\Theta^\eps\big)+3\sup_{s\le t}\abs{\widehat{\mathcal M}(s)} .
\]
Since $\widehat Z(T)\le2(1+T)\big(\Theta^\eps+\norm{X^1}_{\KT}^2\big)$, Theorem~\ref{thm:linear}(i) and Lemma~\ref{lem:CLT-forcing}(i) show that $m_p(t)=\E\widehat Z(t)^p$ is finite. Raising to the power $p$, applying H\"older's inequality with respect to $g(s)\dd s$, taking expectations, inserting Step~3 with $\gamma$ small and absorbing, we find
\[
  m_p(t)\le C_p\Big(\int_0^Tg\dd s\Big)^{p-1}\int_0^tg\,m_p\dd s+C_p\,\eps^p\,\E\big[(\Xi^\eps)^p+(\Theta^\eps)^p\big],\qquad t\in[0,T].
\]
The expectation on the right-hand side is bounded uniformly in $\eps$ by Lemma~\ref{lem:CLT-forcing}. Since $g$ is deterministic and integrable, Gronwall's lemma gives $\E\widehat Z(T)^p\le C_p\,\eps^p$. Finally, $\norm{\widehat X^\eps}_{\KT}^2\le(1+T)\widehat Z(T)$, so
\[
  \E\Big\|\frac{X^\eps-X^0}{\sqrt\eps}-X^1\Big\|_{\KT}^{2p}=\E\norm{\widehat X^\eps}_{\KT}^{2p}\le(1+T)^p\,\E\widehat Z(T)^p\le C_p\,\eps^p . \qedhere
\]
\end{proof}

\begin{remark}\label{rem:CLT}
\begin{enumerate}[label=(\alph*),leftmargin=*]
\item \emph{On the rate.} Formally one expects an expansion $X^\eps=X^0+\sqrt\eps X^1+\eps X^2+o(\eps)$. Here $X^2$ solves the linearised equation forced by the quadratic terms $B(U,U)$, $M(D,D)$, $\Bt(U,D)$ and $\frac12D^2f(d^0)[D,D]$, by the It\^o correction $\frac12G^2(d^0)$, by the noise term $G(D)\dd W_2$ and, when $\sigma$ is differentiable, by the noise term $D\sigma(u^0)[U]\dd W_1$. These are precisely the limits of the forcing terms of \eqref{eq:hatU}--\eqref{eq:hatD} divided by $\sqrt\eps$. Thus $\widehat X^\eps$ is formally of exact order $\sqrt\eps$ unless $X^2$ vanishes; we do not prove a matching lower bound. By contrast, the rate of the law of large numbers is optimal, see Remark~\ref{rem:LLN}(b).
\item \emph{Regularity of $\sigma$.} Differentiability of $\sigma$ is not needed: the limit $X^1$ involves only $\sigma(u^0)$, and the difference $\sigma(u^\eps)-\sigma(u^0)$ is controlled by the Lipschitz condition in \ref{A:sigma} and by Theorem~\ref{thm:LLN}.
\end{enumerate}
\end{remark}
\section{Moderate deviations}\label{sec:MDP}

\subsection{The weak convergence criterion}
Recall from Section~\ref{sec:noise} that $W$ is a Wiener process with values in $U_1=H\times U_2$ and reproducing kernel Hilbert space $\mathcal U=\Hz\times\ell^2$, so that $W$ is a $C([0,T];U_1)$-valued random variable with Cameron--Martin space $\{\int_0^\cdot q\dd s:q\in L^2(0,T;\mathcal U)\}$. For $\Lambda<\infty$ let $\Sb_\Lambda$ be as in Proposition~\ref{prop:skel-cont}, a compact Polish space in the weak topology, and let
\[
  \Ac_\Lambda=\big\{q\colon[0,T]\times\Omega\to\mathcal U\ \text{predictable}:\ q(\cdot,\omega)\in\Sb_\Lambda\ \text{for }\P\text{-a.e. }\omega\big\}.
\]
The weak convergence approach has been used for moderate deviations of several hydrodynamic SPDEs, see e.g.\ \cite{WZZ15,ZhangZhou17,ZZZ16,Slavik22,LiShangWang25}. We use the sufficient condition of Matoussi, Sabbagh and Zhang \cite{MSZ21} for the Laplace principle (see also \cite[Theorem~2.3]{PSZZ22}), which is based on the variational representation of Budhiraja and Dupuis \cite{BD00,BDM08}.

\begin{proposition}\label{prop:criterion}
Let $\mathbb E$ be a Polish space with metric $\mathrm d_{\mathbb E}$, let $a(\eps)\to\infty$, and for $\eps\in(0,1]$ let $\mathcal G^\eps,\mathcal G^0\colon C([0,T];U_1)\to\mathbb E$ be Borel maps. Assume:
\begin{enumerate}[label=\textup{(\alph*)}]
\item for every $\Lambda<\infty$, every family $\{q^\eps\}\subset\Ac_\Lambda$ and every $\eta>0$,
\[
  \lim_{\eps\to0}\P\Big(\mathrm d_{\mathbb E}\Big(\mathcal G^\eps\Big(W+a(\eps)\int_0^\cdot q^\eps\dd s\Big),\ \mathcal G^0\Big(\int_0^\cdot q^\eps\dd s\Big)\Big)>\eta\Big)=0;
\]
\item for every $\Lambda<\infty$ and every sequence $q_n\to q$ in $\Sb_\Lambda$, $\mathcal G^0(\int_0^\cdot q_n\dd s)\to\mathcal G^0(\int_0^\cdot q\dd s)$ in $\mathbb E$.
\end{enumerate}
Then $\{\mathcal G^\eps(W)\}$ satisfies the large deviation principle on $\mathbb E$ with speed $a(\eps)^2$ and good rate function $I(z)=\inf\{\frac12\norm q^2:z=\mathcal G^0(\int_0^\cdot q\dd s)\}$.
\end{proposition}

\begin{proof}
In \cite[Theorem~3.2]{MSZ21} the criterion is formulated for the family $\mathcal G^\eps(W)$ with controls shifted by $\eps^{-1/2}\int_0^\cdot q^\eps$ and speed $\eps^{-1}$; in the normalisation of \cite{PSZZ22} (Theorem~2.3 of the arXiv version) the shift is $\eps^{-1}\int_0^\cdot q^\eps$ and the speed $\eps^{-2}$. In both cases the shift is the square root of the speed. The proof uses only the representation $-\log\E e^{-F(W)}=\inf_{v}\E\big[\frac12\norm v^2+F(W+\int_0^\cdot v)\big]$ of \cite{BD00}, valid for every bounded measurable $F$, applied for fixed $\eps$ to $F=\lambda\,\Psi\circ\mathcal G^\eps$ with $\Psi\in C_b(\mathbb E)$ and $\lambda$ the speed, followed by the substitution $v=\lambda^{1/2}q$ and a passage to the limit along arbitrary sequences $\eps_n\to0$, in which only $\lambda(\eps_n)\to\infty$ is used. The proof therefore applies without change with speed $\lambda(\eps)=a(\eps)^2$ and shift $\lambda(\eps)^{1/2}=a(\eps)$, and no monotonicity of $a$ is needed. The rate function is good: if $I(z)\le M$, a minimising sequence $q_k$ with $\mathcal G^0(\int_0^\cdot q_k)=z$ has a weak limit $q$ with $\frac12\norm q^2\le I(z)$ and, by (b), $\mathcal G^0(\int_0^\cdot q)=z$; hence $\{I\le M\}$ is the image of the compact set $\Sb_{2M}$ under a continuous map. The Laplace principle so obtained is equivalent to the LDP because $\mathbb E$ is Polish and $I$ is good \cite{BD19}.
\end{proof}

We apply Proposition~\ref{prop:criterion} with $\mathbb E=\KT$,
\[
  \mathcal G^\eps(\varphi)=\frac{J^\eps(\varphi)-X^0}{\delta_\eps},\qquad
  \mathcal G^0(\varphi)=\begin{cases}\Gam(\dot\varphi)&\text{if }\varphi=\int_0^\cdot q\dd s\text{ for some }q\in L^2(0,T;\mathcal U),\\ 0&\text{otherwise},\end{cases}
\]
where $J^\eps$ is the map of Theorem~\ref{thm:wellposed}(ii) and $\delta_\eps=\sqrt\eps\,a(\eps)$. The map $\mathcal G^\eps$ is Borel because $J^\eps$ is, and $\mathcal G^0$ is Borel by Lemma~\ref{lem:G0-Borel} below. Moreover $Y^\eps=\mathcal G^\eps(W)$, and $\mathcal G^0(\int_0^\cdot q)=\Gam(q)$, so the rate function of Proposition~\ref{prop:criterion} is the function $I$ of Theorem~\ref{thm:MDP}. Condition~(b) is Proposition~\ref{prop:skel-cont}. It remains to verify condition~(a). Without loss of generality we assume $a(\eps)\ge1$ and $\delta_\eps\le1$.

\begin{lemma}\label{lem:G0-Borel}
The map $\mathcal G^0\colon C([0,T];U_1)\to\KT$ is Borel measurable.
\end{lemma}

\begin{proof}
Let $\mathcal I\colon L^2(0,T;\mathcal U)\to C([0,T];U_1)$, $(\mathcal Iq)(t)=\int_0^tq(s)\dd s$. Since the inclusions $\Hz\hookrightarrow H$ and $\ell^2\hookrightarrow U_2$ are continuous, $\mathcal I$ is continuous, with $\sup_t\abs{(\mathcal Iq)(t)}_{U_1}\le C\sqrt T\norm q_{L^2(0,T;\mathcal U)}$. It is also injective: if $\int_0^tq\dd s=0$ in $U_1$ for all $t$, then $q=0$ a.e.\ by the Lebesgue differentiation theorem, since $\mathcal U\hookrightarrow U_1$ is injective. Both spaces are Polish. By the Lusin--Souslin theorem \cite[Theorem~15.1]{Kechris95}, $\mathcal I$ maps Borel sets onto Borel sets; in particular its range $\mathcal I(L^2(0,T;\mathcal U))$ is Borel, and $\mathcal I^{-1}$ is Borel measurable on it. On this range $\mathcal G^0=\Gam\circ\mathcal I^{-1}$, and $\mathcal G^0=0$ outside it. For a Borel set $F\subset\KT$,
\[
  (\mathcal G^0)^{-1}(F)=\mathcal I\big(\Gam^{-1}(F)\big)\cup\begin{cases}C([0,T];U_1)\setminus\mathcal I(L^2(0,T;\mathcal U))&\text{if }0\in F,\\ \emptyset&\text{otherwise}.\end{cases}
\]
Since $\Gam$ is continuous (Proposition~\ref{prop:skeleton}), $\Gam^{-1}(F)$ is Borel, so this set is Borel.
\end{proof}

\subsection{The controlled fluctuation equation}
Let $\Lambda<\infty$ and $q^\eps=(\theta^\eps,\rho^\eps)\in\Ac_\Lambda$.

\begin{lemma}\label{lem:Girsanov}
The process $X^{\eps,q^\eps}=J^\eps(W+a(\eps)\int_0^\cdot q^\eps\dd s)$ has paths in $\KT$ and solves, $\P$-a.s.,
\begin{align*}
  \dd u+\big[\Ast u+B(u)+M(d)\big]\dd t&=\delta_\eps\,\sigma(u)\theta^\eps\dd t+\sqrt\eps\,\sigma(u)\dd W_1,\\
  \dd d+\big[\An d+\Bt(u,d)+f(d)\big]\dd t&=\delta_\eps\,G(d)\rho^\eps\dd t+\frac\eps2G^2(d)\dd t+\sqrt\eps\,G(d)\dd W_2,
\end{align*}
with $u(0)=u_0$, $d(0)=d_0$.
\end{lemma}

\begin{proof}
Since $a(\eps)^2\int_0^T\abs{q^\eps}_{\mathcal U}^2\dd t\le a(\eps)^2\Lambda$ a.s., Novikov's condition holds, and by Girsanov's theorem \cite{DPZ14} $W+a(\eps)\int_0^\cdot q^\eps$ is a Wiener process under an equivalent probability measure $\mathbb Q^\eps$. Under $\mathbb Q^\eps$, $J^\eps(W+a(\eps)\int_0^\cdot q^\eps)$ solves \eqref{eq:u}--\eqref{eq:d} driven by this Wiener process, which is the stated system because $\sqrt\eps\,a(\eps)=\delta_\eps$; note that the Stratonovich correction is produced by the martingale part only. Since $\mathbb Q^\eps\sim\P$, the equations hold $\P$-a.s. See \cite[Theorem~5.4]{BMP20} for the details of this argument.
\end{proof}

Let $Y=(Y_u,Y_d)=(X^{\eps,q^\eps}-X^0)/\delta_\eps=\mathcal G^\eps(W+a(\eps)\int_0^\cdot q^\eps)$; we suppress the dependence on $\eps$. Expanding the nonlinearities around $X^0$ as in Section~\ref{sec:LLN}, $Y(0)=0$ and
\begin{align}
  &\dd Y_u+\big[\Lin_1(Y)+\delta_\eps B(Y_u,Y_u)+\delta_\eps M(Y_d,Y_d)\big]\dd t \notag\\
  &\qquad=\sigma(u^0+\delta_\eps Y_u)\theta^\eps\dd t+\frac1{a(\eps)}\sigma(u^0+\delta_\eps Y_u)\dd W_1, \label{eq:Yu}\\
  &\dd Y_d+\Big[\Lin_2^\flat(Y)+\delta_\eps\Bt(Y_u,Y_d)+\frac{f(d^0+\delta_\eps Y_d)-f(d^0)}{\delta_\eps}\Big]\dd t \notag\\
  &\qquad=G(d^0+\delta_\eps Y_d)\rho^\eps\dd t+\frac{\delta_\eps}{2a(\eps)^2}G^2(d^0+\delta_\eps Y_d)\dd t+\frac1{a(\eps)}G(d^0+\delta_\eps Y_d)\dd W_2. \label{eq:Yd}
\end{align}
Here we used $\sqrt\eps/\delta_\eps=a(\eps)^{-1}$ and $\eps/(2\delta_\eps)=\delta_\eps/(2a(\eps)^2)$. Thus the deterministic control enters at order one, while the noise is multiplied by $a(\eps)^{-1}\to0$ and the nonlinear corrections by $\delta_\eps\to0$.

\begin{lemma}\label{lem:MDP-apriori}
Let $\tau^\eps=\inf\{t\in[0,T]:\delta_\eps\hone{Y_d(t)}\ge1\}\wedge T$. There is $C_\Lambda$, independent of $\eps$ and of the family $\{q^\eps\}\subset\Ac_\Lambda$, such that
\[
  \E\,\Ec_{\tau^\eps}(Y)\le C_\Lambda\qquad\text{and}\qquad\P(\tau^\eps<T)\le C_\Lambda\delta_\eps^2 ,
\]
where $\Ec_\tau(Y)=\sup_{t\le\tau}(\abs{Y_u}^2+\hone{Y_d}^2)+\int_0^\tau(\norm{Y_u}^2+\abs{\nabla Y_d}^2+\abs{\An Y_d}^2)\dd t$.
\end{lemma}

\begin{proof}
Let $\mathcal Y=\abs{Y_u}^2+\hone{Y_d}^2$, $\Dc=\norm{Y_u}^2+\abs{\nabla Y_d}^2+\abs{\An Y_d}^2$, and write $\delta=\delta_\eps$, $a=a(\eps)$. By Lemma~\ref{lem:Ito},
\begin{align*}
  \dd\mathcal Y&=-2\Big[\dual{\Lin_1(Y)}{Y_u}+\scal{\Lin_2^\flat(Y)}{(I+\An)Y_d}\Big]\dd t\\
  &\quad-2\delta\Big[\dual{B(Y_u,Y_u)}{Y_u}+\dual{M(Y_d,Y_d)}{Y_u}+\scal{\Bt(Y_u,Y_d)}{(I+\An)Y_d}\Big]\dd t\\
  &\quad-\frac2\delta\scal{f(d^0+\delta Y_d)-f(d^0)}{(I+\An)Y_d}\dd t\\
  &\quad+2\scal{\sigma(u^0+\delta Y_u)\theta^\eps}{Y_u}\dd t+2\sone{G(d^0+\delta Y_d)\rho^\eps}{Y_d}\dd t\\
  &\quad+\frac\delta{a^2}\sone{G^2(d^0+\delta Y_d)}{Y_d}\dd t+\frac1{a^2}\Big[\abs{\sigma(u^0+\delta Y_u)}_{\LHS}^2+\abs{G(d^0+\delta Y_d)}_{\LHSd}^2\Big]\dd t+\dd\mathcal M,
\end{align*}
with $\mathcal M=\frac2a\int_0^\cdot\scal{Y_u}{\sigma(u^0+\delta Y_u)\dd W_1}+\frac2a\int_0^\cdot\sone{G(d^0+\delta Y_d)\dd W_2}{Y_d}$. We estimate the terms on $[0,\tau^\eps]$.
\begin{itemize}[leftmargin=*]
\item The first bracket is $\ge\frac12\Dc-C_*g\mathcal Y$ by Lemma~\ref{lem:coercive}(i).
\item The second bracket vanishes by Lemma~\ref{lem:nonlinear}(i) and Lemma~\ref{lem:cancel}.
\item With $e=\delta Y_d$, Lemma~\ref{lem:fspace}(b),(c) gives
\begin{align*}
  -\frac2\delta\scal{f(d^0+e)-f(d^0)}{(I+\An)Y_d}&=-\frac2{\delta^2}\scal{f(d^0+e)-f(d^0)}{(I+\An)e}\\
  &\le2\kappa\hone{Y_d}^2+Cg\big(\hone{Y_d}^2+\delta^{2N-1}\hone{Y_d}^{2N+1}\big),
\end{align*}
and $\delta^{2N-1}\hone{Y_d}^{2N+1}=(\delta\hone{Y_d})^{2N-1}\hone{Y_d}^2\le\hone{Y_d}^2$ for $t\le\tau^\eps$.
\item By \ref{A:sigma}, $\abs{\sigma(v)}_{\mathscr L(\Hz,H)}\le\sqrt K(1+\abs v)$; since $\delta\le1$ and $\sup_t\abs{u^0}<\infty$, $2\scal{\sigma(u^0+\delta Y_u)\theta^\eps}{Y_u}\le C(1+\abs{\theta^\eps}_{\Hz}^2)(1+\mathcal Y)$.
\item By Lemma~\ref{lem:G}(ii), $\scal{G(Y_d)\rho}{Y_d}=0$ for $\rho\in\ell^2$, so that $\sone{G(d^0+\delta Y_d)\rho^\eps}{Y_d}=\sone{G(d^0)\rho^\eps}{Y_d}+\delta\scal{\nabla(G(Y_d)\rho^\eps)}{\nabla Y_d}$. By Lemma~\ref{lem:G}(i) and the Cauchy--Schwarz inequality,
\[
  2\sone{G(d^0+\delta Y_d)\rho^\eps}{Y_d}\le2\big(\abs{G(d^0)}_{\LHSd}+\delta\abs{G(Y_d)}_{\LHSd}\big)\hone{Y_d}\abs{\rho^\eps}_{\ell^2}\le C(1+\abs{\rho^\eps}_{\ell^2}^2)(1+\mathcal Y).
\]
\item The It\^o-correction and quadratic-variation terms are bounded by $C(1+\mathcal Y)$, and $\dd\langle\mathcal M\rangle\le Ca^{-2}(1+\mathcal Y)^2\dd t$.
\end{itemize}
Hence, on $[0,\tau^\eps]$,
\[
  \dd\mathcal Y+\Dc\dd t\le\beta(t)(1+\mathcal Y)\dd t+\dd\mathcal M,\qquad\beta=C\big(g+\abs{q^\eps}_{\mathcal U}^2\big),
\]
and $\int_0^T\beta\dd t\le C(\int_0^Tg\dd t+\Lambda)=:b_\Lambda$ almost surely.
Let $\chi(t)=\exp(-\int_0^t\beta\dd s)\in[e^{-b_\Lambda},1]$, an adapted continuous process of bounded variation. By It\^o's product rule, $\dd[\chi(1+\mathcal Y)]\le-\chi\Dc\dd t+\chi\dd\mathcal M$ on $[0,\tau^\eps]$. Let $\tau'_R=\tau^\eps\wedge\inf\{t:\mathcal Y(t)+\int_0^t\Dc\ge R\}$. Integrating up to $t\wedge\tau'_R$ and using the Burkholder--Davis--Gundy inequality,
\begin{align*}
  \psi(t):=\E\sup_{s\le t}\big[\chi(1+\mathcal Y)\big](s\wedge\tau'_R)+\E\int_0^{t\wedge\tau'_R}\chi\Dc\dd s
  &\le2+C\E\Big(\int_0^{t\wedge\tau'_R}\frac{\chi^2(1+\mathcal Y)^2}{a^2}\dd s\Big)^{1/2}\\
  &\le2+\frac14\psi(t)+\frac C{a^2}\int_0^t\psi(s)\dd s,
\end{align*}
where we used $\chi^2(1+\mathcal Y)^2\le\sup[\chi(1+\mathcal Y)]\cdot\chi(1+\mathcal Y)$. Since $\psi$ is finite, Gronwall's lemma gives $\psi(T)\le C$. As $\chi\ge e^{-b_\Lambda}$, letting $R\to\infty$ (Fatou's lemma) yields $\E\,\Ec_{\tau^\eps}(Y)\le Ce^{b_\Lambda}=:C_\Lambda$. Finally, on $\{\tau^\eps<T\}$ we have $\delta^2\hone{Y_d(\tau^\eps)}^2=1$, whence $\P(\tau^\eps<T)\le\delta^2\E\,\Ec_{\tau^\eps}(Y)\le C_\Lambda\delta^2$.
\end{proof}

\subsection{Verification of condition (a)}
Let $Z=Z^\eps=\Gam(q^\eps)$, defined pathwise. Since $q^\eps$ is predictable and $\Gam$ is causal (Proposition~\ref{prop:skeleton}), $Z$ is adapted with continuous paths, and $\Ec_T(Z)\le C_\Gam\Lambda$ a.s.\ by \eqref{eq:skel-bound}.

\begin{lemma}\label{lem:MDP-conv}
For every $\eta>0$, $\lim_{\eps\to0}\P\big(\norm{Y-Z}_{\KT}>\eta\big)=0$.
\end{lemma}

\begin{proof}
Write $\delta=\delta_\eps$, $a=a(\eps)$, and $\widetilde Y=(\widetilde Y_u,\widetilde Y_d)=Y-Z$. Using Lemma~\ref{lem:fspace}(d) we write
\[
  \frac{f(d^0+\delta Y_d)-f(d^0)}{\delta}=Df(d^0)Y_d+\mathcal R_\delta,\qquad\mathcal R_\delta=\frac1\delta\,\mathcal R(d^0,\delta Y_d).
\]
Subtracting \eqref{eq:skel-u}--\eqref{eq:skel-d} from \eqref{eq:Yu}--\eqref{eq:Yd}, and using $G(d^0+\delta Y_d)-G(d^0)=\delta G(Y_d)$, we obtain $\widetilde Y(0)=0$ and
\begin{align*}
  \dd\widetilde Y_u+\Lin_1(\widetilde Y)\dd t&=\big[\sigma(u^0+\delta Y_u)-\sigma(u^0)\big]\theta^\eps\dd t-\delta\big[B(Y_u,Y_u)+M(Y_d,Y_d)\big]\dd t\\
  &\quad+\frac1a\sigma(u^0+\delta Y_u)\dd W_1,\\
  \dd\widetilde Y_d+\Lin_2(\widetilde Y)\dd t&=\Big[\delta G(Y_d)\rho^\eps-\delta\Bt(Y_u,Y_d)-\mathcal R_\delta+\frac\delta{2a^2}G^2(d^0+\delta Y_d)\Big]\dd t\\
  &\quad+\frac1aG(d^0+\delta Y_d)\dd W_2 .
\end{align*}
For $R>0$ let $\tau_R=\inf\{t:\Ec_t(Y)\ge R\}\wedge T$ and $\tau_*=\tau^\eps\wedge\tau_R$, with $\tau^\eps$ as in Lemma~\ref{lem:MDP-apriori}. Let $\widetilde{\mathcal Y}=\abs{\widetilde Y_u}^2+\hone{\widetilde Y_d}^2$ and $\widetilde\Dc=\norm{\widetilde Y_u}^2+\abs{\nabla\widetilde Y_d}^2+\abs{\An\widetilde Y_d}^2$. On $[0,\tau_*]$ we have $\Ec_t(Y)\le R$ and $\delta\hone{Y_d}\le1$, and Lemma~\ref{lem:fspace}(d) gives $\abs{\mathcal R_\delta}\le C\delta\big(1+\hone{d^0}^{2N-1}+1\big)\hone{Y_d}^2\le C\delta R$. By Lemma~\ref{lem:Ito}, Lemma~\ref{lem:coercive}(i), Lemma~\ref{lem:nonlinear}, \ref{A:sigma}, Lemma~\ref{lem:G} and Young's inequality, on $[0,\tau_*]$,
\begin{align*}
  2\scal{[\sigma(u^0+\delta Y_u)-\sigma(u^0)]\theta^\eps}{\widetilde Y_u}&\le2\sqrt L\,\delta\abs{Y_u}\abs{\theta^\eps}_{\Hz}\abs{\widetilde Y_u}\le\abs{\theta^\eps}_{\Hz}^2\widetilde{\mathcal Y}+L\delta^2R,\\
  2\delta\abs{\dual{B(Y_u,Y_u)+M(Y_d,Y_d)}{\widetilde Y_u}}&\le\tfrac14\norm{\widetilde Y_u}^2+C\delta^2\big(\abs{Y_u}^2\norm{Y_u}^2+\abs{\nabla Y_d}^2\hda{Y_d}^2\big),\\
  2\delta\sone{G(Y_d)\rho^\eps}{\widetilde Y_d}&\le\abs{\rho^\eps}_{\ell^2}^2\widetilde{\mathcal Y}+C\delta^2R,\\
  2\abs{\scal{\delta\Bt(Y_u,Y_d)+\mathcal R_\delta}{(I+\An)\widetilde Y_d}}&\le\tfrac14\abs{\An\widetilde Y_d}^2+2\widetilde{\mathcal Y}+C\delta^2\big(\abs{Y_u}\norm{Y_u}\abs{\nabla Y_d}\hda{Y_d}+R^2\big),\\
  \frac\delta{a^2}\sone{G^2(d^0+\delta Y_d)}{\widetilde Y_d}&\le\widetilde{\mathcal Y}+C\delta^2a^{-4}(1+R),
\end{align*}
and the quadratic variation terms are bounded by $Ca^{-2}(1+R)$. Integrating in time and using $\Ec_{\tau_*}(Y)\le R$, we obtain on $[0,\tau_*]$
\[
  \dd\widetilde{\mathcal Y}+\tfrac12\widetilde\Dc\dd t\le\widetilde\beta\,\widetilde{\mathcal Y}\dd t+\zeta\dd t+\dd\widetilde{\mathcal M},\qquad\widetilde\beta=2C_*g+4+\abs{q^\eps}_{\mathcal U}^2,
\]
for a non-negative process $\zeta$ with $\int_0^{\tau_*}\zeta\dd t\le C_R(\delta^2+a^{-2})$ a.s.,
where $\widetilde{\mathcal M}=\frac2a\int_0^\cdot\scal{\widetilde Y_u}{\sigma(u^0+\delta Y_u)\dd W_1}+\frac2a\int_0^\cdot\sone{G(d^0+\delta Y_d)\dd W_2}{\widetilde Y_d}$ satisfies $\dd\langle\widetilde{\mathcal M}\rangle\le Ca^{-2}(1+R)\widetilde{\mathcal Y}\dd t$ on $[0,\tau_*]$. Let $\widetilde\chi(t)=\exp(-\int_0^t\widetilde\beta)\in[e^{-\tilde b_\Lambda},1]$ with $\tilde b_\Lambda=\int_0^T(2C_*g+4)+\Lambda$. As in the proof of Lemma~\ref{lem:MDP-apriori} (here all quantities are bounded on $[0,\tau_*]$, since $\Ec_{\tau_*}(Y)\le R$ and $\Ec_T(Z)\le C_\Gam\Lambda$),
\begin{align*}
  \E\sup_{t\le\tau_*}\widetilde\chi\widetilde{\mathcal Y}+\frac12\E\int_0^{\tau_*}\widetilde\chi\widetilde\Dc\dd t&\le C_R(\delta^2+a^{-2})+C\E\Big(\int_0^{\tau_*}\frac{1+R}{a^2}\widetilde\chi\widetilde{\mathcal Y}\dd t\Big)^{1/2}\\
  &\le C_R(\delta^2+a^{-2})+\frac14\E\sup_{t\le\tau_*}\widetilde\chi\widetilde{\mathcal Y}+\frac{C(1+R)T}{a^2},
\end{align*}
whence $\E\,\Ec_{\tau_*}(\widetilde Y)\le C_{\Lambda,R}(\delta_\eps^2+a(\eps)^{-2})$. Now, by Lemma~\ref{lem:MDP-apriori} and Chebyshev's inequality,
\begin{align*}
  \P\big(\norm{\widetilde Y}_{\KT}>\eta\big)&\le\P(\tau^\eps<T)+\P\big(\tau_R<T,\ \tau^\eps=T\big)+\P\big(\Ec_{\tau_*}(\widetilde Y)>\eta^2/(1+T),\ \tau_*=T\big)\\
  &\le C_\Lambda\delta_\eps^2+\P\big(\Ec_{\tau^\eps}(Y)\ge R\big)+\frac{(1+T)C_{\Lambda,R}}{\eta^2}\big(\delta_\eps^2+a(\eps)^{-2}\big)\\
  &\le C_\Lambda\delta_\eps^2+\frac{C_\Lambda}R+\frac{(1+T)C_{\Lambda,R}}{\eta^2}\big(\delta_\eps^2+a(\eps)^{-2}\big).
\end{align*}
Letting first $\eps\to0$ and then $R\to\infty$ proves the lemma.
\end{proof}

\begin{proof}[Proof of Theorem~\ref{thm:MDP}]
By Lemma~\ref{lem:MDP-conv}, condition~(a) of Proposition~\ref{prop:criterion} holds, since $\mathcal G^\eps(W+a(\eps)\int_0^\cdot q^\eps)=Y$ and $\mathcal G^0(\int_0^\cdot q^\eps)=\Gam(q^\eps)=Z$. Condition~(b) is Proposition~\ref{prop:skel-cont}. Proposition~\ref{prop:criterion} yields the theorem.
\end{proof}

\subsection{Observables}
\begin{proof}[Proof of Corollary~\ref{cor:observables}]
Let $\ell=\ell_\Phi$ and let $k_\ell$ be the Riesz representative of $\ell\circ\Gam$, so that $\norm{k_\ell}^2=s_\Phi^2$ by Theorem~\ref{thm:linear}(ii).

\emph{Step 1.} By the contraction principle \cite[Theorem~4.2.1]{DZ98}, $\ell(Y^\eps)$ satisfies the LDP on $\R$ with speed $a(\eps)^2$ and rate function
\[
  \mathcal I(s)=\inf\{I(z):\ell(z)=s\}=\inf\Big\{\tfrac12\norm q^2:\langle k_\ell,q\rangle=s\Big\}=\frac{s^2}{2s_\Phi^2},
\]
where for $s_\Phi>0$ the infimum is attained at $q=sk_\ell/s_\Phi^2$ by the Cauchy--Schwarz inequality, and for $s_\Phi=0$ we have $\mathcal I(s)=+\infty$ for $s\ne0$.

\emph{Step 2.} Let $\widetilde\Phi^\eps=(\Phi(X^\eps)-\Phi(X^0))/\delta_\eps$. On $\{\delta_\eps\norm{Y^\eps}_{\KT}\le1\}$ we have $\abs{\widetilde\Phi^\eps-\ell(Y^\eps)}\le C_\Phi\delta_\eps\norm{Y^\eps}_{\KT}^2$. Hence, for $\eta>0$ and $R>0$, if $\eps$ is so small that $\delta_\eps^{-1}\ge R$ and $\eta/(C_\Phi\delta_\eps)\ge R^2$, then $\{\abs{\widetilde\Phi^\eps-\ell(Y^\eps)}>\eta\}\subset\{\norm{Y^\eps}_{\KT}\ge R\}$. By \eqref{eq:skel-bound}, $I(z)\ge\norm z_{\KT}^2/(2(1+T)C_\Gam)$, and the upper bound of Theorem~\ref{thm:MDP} applied to the closed set $\{\norm z_{\KT}\ge R\}$ gives
\[
  \limsup_{\eps\to0}\frac1{a(\eps)^2}\log\P\big(\abs{\widetilde\Phi^\eps-\ell(Y^\eps)}>\eta\big)\le-\frac{R^2}{2(1+T)C_\Gam}.
\]
Since $R$ is arbitrary, $\widetilde\Phi^\eps$ and $\ell(Y^\eps)$ are exponentially equivalent at speed $a(\eps)^2$, so by \cite[Theorem~4.2.13]{DZ98} $\widetilde\Phi^\eps$ satisfies the LDP with the rate function $\mathcal I$.

\emph{Step 3.} If $s_\Phi>0$, $\mathcal I$ is continuous, so its infimum over the closed set $[r,\infty)$ and over the open set $(r,\infty)$ both equal $r^2/(2s_\Phi^2)$; the LDP bounds give the stated limit, since $\{\widetilde\Phi^\eps\ge r\}=\{\Phi(X^\eps)-\Phi(X^0)\ge r\sqrt\eps\,a(\eps)\}$. If $s_\Phi=0$, the upper bound gives $-\infty$.
\end{proof}

\section{Uniformity in the initial data}\label{sec:uniform}

In this section we prove the uniform versions of our results announced in Section~\ref{sec:main}. Uniformity over \emph{bounded} sets of initial data is the form needed for exit problems in infinite dimensions, where closed balls are not compact \cite{SBD19}; for a recent uniform LDP for a hydrodynamic model, see \cite{SunYuLiu26}. We use the equivalences between the various notions of uniform large deviation principles established by Salins \cite{Salins19}.

\subsection{Notation}
Let $\Xs=H\times\HH$ with $\norm x_{\Xs}^2=\abs{u_0}^2+\hone{d_0}^2$ for $x=(u_0,d_0)$, and let $\mathfrak B_\varrho=\{x\in\Xs:\norm x_{\Xs}\le\varrho\}$. For $x\in\Xs$ we denote by $X^{\eps,x}=J^\eps(x,W)$ and $X^{0,x}=(u^{0,x},d^{0,x})$ the solutions of \eqref{eq:u}--\eqref{eq:d} with initial datum $x$ (Theorem~\ref{thm:app-wellposed}), and by $g_x$, $\Lin^x$, $\Gam_x$, $X^{1,x}$ and $I_x$ the objects of Section~\ref{sec:results} built on $X^{0,x}$, namely the function \eqref{eq:g}, the linearised operator, the skeleton map, the solution of \eqref{eq:lin-u}--\eqref{eq:lin-d} and the rate function of Theorem~\ref{thm:MDP}. Finally
\[
  Y^{\eps,x}=\frac{X^{\eps,x}-X^{0,x}}{\delta_\eps},\qquad \delta_\eps=\sqrt\eps\,a(\eps),\qquad \Upsilon_x(s)=\{z\in\KT:I_x(z)\le s\}.
\]

\subsection{Uniform bounds and Lipschitz dependence}
\begin{lemma}\label{lem:unif-bounds}
For every $\varrho>0$ there is $C_\varrho$ such that
\[
  \sup_{x\in\mathfrak B_\varrho}\Big(\norm{X^{0,x}}_{\KT}^2+\int_0^Tg_x(t)\dd t\Big)\le C_\varrho .
\]
Consequently, the constants $C_\Gam$ of Proposition~\ref{prop:skeleton} and $C_\Lambda$, $C_{\Lambda,R}$ of Lemmas~\ref{lem:MDP-apriori} and \ref{lem:MDP-conv} can be chosen uniformly for $x\in\mathfrak B_\varrho$; they depend on $q^\eps$ only through $\Lambda$.
\end{lemma}

\begin{proof}
By Proposition~\ref{prop:moments} with $\eps=0$, $\norm{X^{0,x}}_{\KT}^2\le C(1+\varrho^2)$ for $x\in\mathfrak B_\varrho$. Writing $X^{0,x}=(u,d)$ and using \eqref{eq:g},
\begin{align*}
  \int_0^Tg_x\dd t&\le T+\Big(1+\sup_t\abs u^2\Big)\int_0^T\norm u^2\dd t+\sup_t\abs{\nabla d}^2\int_0^T\hda d^2\dd t\\
  &\quad+\big(1+\sup_t\hone d^{4N}\big)\Big(T+\sqrt T\Big(\int_0^T\hda d^2\dd t\Big)^{1/2}\Big),
\end{align*}
which is bounded by a constant depending only on $\varrho$. For the last statement we inspect the proofs. The constants of Lemma~\ref{lem:coercive} do not depend on $X^0$. The constants in Proposition~\ref{prop:skeleton} and Lemmas~\ref{lem:MDP-apriori}, \ref{lem:MDP-conv} depend on the initial datum only through $\int_0^Tg\dd t$, $\sup_t\abs{u^0(t)}$ and $\sup_t\hone{d^0(t)}$, and on the controls only through the bound $\Lambda$ on $\int_0^T\abs{q^\eps}_{\mathcal U}^2\dd t$.
\end{proof}

As a first consequence, the quantitative law of large numbers and central limit theorem hold uniformly over bounded sets of initial data.

\begin{corollary}\label{cor:unif-LLN-CLT}
For all $\varrho>0$ and $p\ge1$ there is $C_{\varrho,p}$ such that for all $\eps\in(0,1]$
\begin{align*}
  \sup_{x\in\mathfrak B_\varrho}\E\norm{X^{\eps,x}-X^{0,x}}_{\KT}^{2p}&\le C_{\varrho,p}\,\eps^p,\qquad
  \sup_{x\in\mathfrak B_\varrho}\E\norm{X^{1,x}}_{\KT}^{2p}\le C_{\varrho,p},\\
  \sup_{x\in\mathfrak B_\varrho}\E\Big\|\frac{X^{\eps,x}-X^{0,x}}{\sqrt\eps}-X^{1,x}\Big\|_{\KT}^{2p}&\le C_{\varrho,p}\,\eps^p .
\end{align*}
\end{corollary}

\begin{proof}
We inspect the constants in Theorems~\ref{thm:LLN}, \ref{thm:linear}(i) and \ref{thm:CLT}.
\begin{itemize}[leftmargin=*]
\item By Remark~\ref{rem:LLN}(a), the constant of Theorem~\ref{thm:LLN} depends on $x$ only through $\int_0^Tg_x\dd t$ and $\norm x_{\Xs}^2$.
\item By \eqref{eq:lin-moment}, the constant of Theorem~\ref{thm:linear}(i) depends on $x$ only through $\int_0^Tg_x\dd t$ and the quantity $\Sigma_*$ of \eqref{eq:kappa}. The latter is bounded by $K(1+\sup_t\abs{u^{0,x}(t)}^2)+C_h\sup_t\hone{d^{0,x}(t)}^2\le C(1+\norm{X^{0,x}}_{\KT}^2)$.
\item In Section~\ref{sec:CLT}, the constant of Lemma~\ref{lem:CLT-forcing} depends on $x$ only through $\sup_t\hone{d^{0,x}(t)}$. The moments of $\Xi^\eps$ are controlled by those of $\Theta^\eps$, hence by the constant of Theorem~\ref{thm:LLN}, and by $\E\norm{X^{\eps,x}}_{\KT}^{2q}\le C_q(1+\norm x_{\Xs}^2)^q$ (Proposition~\ref{prop:moments}). The Gronwall argument in Step~4 of the proof of Theorem~\ref{thm:CLT} adds only the factor $\exp(C_p(\int_0^Tg_x)^p)$.
\end{itemize}
By Lemma~\ref{lem:unif-bounds}, $\int_0^Tg_x\dd t$ and $\norm{X^{0,x}}_{\KT}$ are bounded uniformly for $x\in\mathfrak B_\varrho$, which proves the three bounds.
\end{proof}

\begin{lemma}\label{lem:lipschitz}
For all $\varrho,\Lambda<\infty$ there are constants $C_\varrho$ and $C_{\varrho,\Lambda}$ such that for all $x,y\in\mathfrak B_\varrho$:
\begin{enumerate}[label=\textup{(\roman*)}]
\item $\norm{X^{0,x}-X^{0,y}}_{\KT}\le C_\varrho\norm{x-y}_{\Xs}$;
\item $\sup_{q\in\Sb_\Lambda}\norm{\Gam_x(q)-\Gam_y(q)}_{\KT}\le C_{\varrho,\Lambda}\norm{x-y}_{\Xs}$.
\end{enumerate}
\end{lemma}

\begin{proof}
(i) Let $Z=(Z_u,Z_d)=X^{0,x}-X^{0,y}$. Expanding the nonlinearities around $X^{0,y}$ as in Section~\ref{sec:LLN}, $Z$ solves
\[
  \partial_tZ+\Lin^{y,\flat}Z+\big(B(Z_u,Z_u)+M(Z_d,Z_d),\ \Bt(Z_u,Z_d)+f(d^{0,y}+Z_d)-f(d^{0,y})\big)=0,\qquad Z(0)=x-y,
\]
where $\Lin^{y,\flat}=(\Lin_1^y,\Lin_2^{y,\flat})$. Let $\mathcal Y=\abs{Z_u}^2+\hone{Z_d}^2$ and $\Dc=\norm{Z_u}^2+\abs{\nabla Z_d}^2+\abs{\An Z_d}^2$. As in the proof of Lemma~\ref{lem:LLN-energy}, using Lemmas~\ref{lem:coercive}(i), \ref{lem:cancel} and \ref{lem:fspace}(b),(c),
\[
  \frac{\dd}{\dd t}\mathcal Y+\Dc\le Cg_y\big(\mathcal Y+\hone{Z_d}^{2N+1}\big).
\]
By Lemma~\ref{lem:unif-bounds}, $\hone{Z_d(t)}\le\norm{X^{0,x}}_{\KT}+\norm{X^{0,y}}_{\KT}\le2C_\varrho^{1/2}$, so that $\hone{Z_d}^{2N+1}\le C_\varrho\mathcal Y$. Gronwall's lemma and Lemma~\ref{lem:unif-bounds} give $\Ec_T(Z)\le2\norm{x-y}_{\Xs}^2\exp(C_\varrho\int_0^Tg_y)\le C_\varrho\norm{x-y}_{\Xs}^2$.

(ii) Let $q=(\theta,\rho)\in\Sb_\Lambda$, $P=\Gam_x(q)$, $E=X^{0,x}-X^{0,y}$ and $\mathcal V=\Gam_x(q)-\Gam_y(q)$. Then $\mathcal V(0)=0$ and
\[
  \partial_t\mathcal V+\Lin^y\mathcal V=-\big(\Lin^x-\Lin^y\big)P+\big(\big[\sigma(u^{0,x})-\sigma(u^{0,y})\big]\theta,\ G(E_d)\rho\big).
\]
The operator $\Lin^x-\Lin^y$ is bilinear in $(E,P)$:
\begin{align*}
  \big[(\Lin^x-\Lin^y)P\big]_u&=B(P_u,E_u)+B(E_u,P_u)+M(P_d,E_d)+M(E_d,P_d),\\
  \big[(\Lin^x-\Lin^y)P\big]_d&=\Bt(P_u,E_d)+\Bt(E_u,P_d)+\big[Df(d^{0,x})-Df(d^{0,y})\big]P_d .
\end{align*}
By Lemma~\ref{lem:nonlinear}(ii)--(iv), each of the first six terms is bounded in $V'$ or in $\LL$ by $C\big(\abs a\norm a\abs b\norm b\big)^{1/2}$, where $(a,b)$ is a pair of components of $P$ and $E$, measured in $H$ and $V$ for velocity components and in $\HH$ and $\DA$ for director components. By Lemma~\ref{lem:f}(iii), H\"older's inequality and \eqref{eq:GN},
\[
  \abs{\big[Df(d^{0,x})-Df(d^{0,y})\big]P_d}\le C\big(1+\hone{d^{0,x}}^{2N-1}+\hone{d^{0,y}}^{2N-1}\big)\hone{E_d}\hone{P_d}.
\]
Integrating in time and using the Cauchy--Schwarz inequality, we obtain
\[
  \int_0^T\abs{(\Lin^x-\Lin^y)P}_{V'\times\LL}^2\dd t\le C_\varrho\norm E_{\KT}^2\norm P_{\KT}^2 .
\]
Moreover, by \ref{A:sigma} and Lemma~\ref{lem:G}(i),
\[
  \abs{[\sigma(u^{0,x})-\sigma(u^{0,y})]\theta}\le\sqrt L\abs{E_u}\abs\theta_{\Hz},\qquad \hone{G(E_d)\rho}\le C_h^{1/2}\hone{E_d}\abs\rho_{\ell^2}.
\] With $\mathcal Y=\abs{\mathcal V_u}^2+\hone{\mathcal V_d}^2$ and $\Dc$ as before, Lemma~\ref{lem:coercive}(i) and Young's inequality give
\[
  \frac{\dd}{\dd t}\mathcal Y+\frac12\Dc\le\big(2C_*g_y+3+\abs q_{\mathcal U}^2\big)\mathcal Y+C\abs{(\Lin^x-\Lin^y)P}_{V'\times\LL}^2+C\norm E_{\KT}^2 .
\]
By Gronwall's lemma, Proposition~\ref{prop:skeleton} (which gives $\norm P_{\KT}^2\le(1+T)C_\Gam\Lambda$), Lemma~\ref{lem:unif-bounds} and part (i),
\[
  \Ec_T(\mathcal V)\le C_\varrho\norm E_{\KT}^2\big(\norm P_{\KT}^2+T\big)\exp\Big(\int_0^T(2C_*g_y+3)\dd t+\Lambda\Big)\le C_{\varrho,\Lambda}\norm{x-y}_{\Xs}^2 .\qedhere
\]
\end{proof}

\subsection{The uniform moderate deviation principle}
\begin{theorem}\label{thm:uniform}
Assume \eqref{eq:a} and let $\varrho>0$.
\begin{enumerate}[label=\textup{(\roman*)}]
\item \emph{Equicontinuous uniform Laplace principle.} For every family $\mathscr L\subset C_b(\KT)$ of equibounded and uniformly equicontinuous functions, i.e.\ such that
\[
  \sup_{\Psi\in\mathscr L}\sup_{z\in\KT}\abs{\Psi(z)}<\infty\qquad\text{and}\qquad\lim_{\eta\to0}\ \sup_{\Psi\in\mathscr L}\ \sup_{\norm{z-z'}_{\KT}\le\eta}\abs{\Psi(z)-\Psi(z')}=0,
\]
we have
\[
  \lim_{\eps\to0}\ \sup_{x\in\mathfrak B_\varrho}\ \sup_{\Psi\in\mathscr L}\ \Big|\frac1{a(\eps)^2}\log\E\exp\big(-a(\eps)^2\Psi(Y^{\eps,x})\big)+\inf_{z\in\KT}\big\{\Psi(z)+I_x(z)\big\}\Big|=0 .
\]
\item \emph{Freidlin--Wentzell uniform bounds.} For all $s_0>0$ and $\eta>0$,
\begin{align*}
  &\liminf_{\eps\to0}\ \inf_{x\in\mathfrak B_\varrho}\ \inf_{z\in\Upsilon_x(s_0)}\Big[\frac1{a(\eps)^2}\log\P\big(\norm{Y^{\eps,x}-z}_{\KT}<\eta\big)+I_x(z)\Big]\ge0,\\
  &\limsup_{\eps\to0}\ \sup_{x\in\mathfrak B_\varrho}\ \sup_{s\in[0,s_0]}\Big[\frac1{a(\eps)^2}\log\P\big(\mathrm{dist}_{\KT}(Y^{\eps,x},\Upsilon_x(s))\ge\eta\big)+s\Big]\le0 .
\end{align*}
\item \emph{Dembo--Zeitouni uniform bounds on compact sets.} For every compact set $\mathcal K\subset\Xs$, every closed $F\subset\KT$ and every open $O\subset\KT$,
\begin{align*}
  \limsup_{\eps\to0}\ \sup_{x\in\mathcal K}\ \frac1{a(\eps)^2}\log\P(Y^{\eps,x}\in F)&\le-\inf_{x\in\mathcal K}\ \inf_{z\in F}I_x(z),\\
  \liminf_{\eps\to0}\ \inf_{x\in\mathcal K}\ \frac1{a(\eps)^2}\log\P(Y^{\eps,x}\in O)&\ge-\sup_{x\in\mathcal K}\ \inf_{z\in O}I_x(z).
\end{align*}
\item \emph{Continuity of the Gaussian limit.} For every $\ell\in\KT^*$ the map $x\mapsto\operatorname{Var}\ell(X^{1,x})=\norm{\ell\circ\Gam_x}^2$ is Lipschitz continuous on $\mathfrak B_\varrho$, and $x\mapsto\Upsilon_x(s)$ is Lipschitz continuous on $\mathfrak B_\varrho$ for the Hausdorff distance on closed subsets of $\KT$, for every $s\ge0$.
\end{enumerate}
\end{theorem}

\begin{proof}
\emph{Step 1: the uniform convergence in probability.} Fix $\Lambda<\infty$ and $\eta>0$. For $x\in\mathfrak B_\varrho$ and $q\in\Ac_\Lambda$, let 
$$
Y^{\eps,x,q}=\mathcal G_x^\eps(W+a(\eps)\int_0^\cdot q\dd s),
$$
 where $\mathcal G^\eps_x(\varphi)=(J^\eps(x,\varphi)-X^{0,x})/\delta_\eps$, and let $\mathcal G^0_x(\int_0^\cdot q\dd s)=\Gam_x(q)$ for $q\in L^2(0,T;\mathcal U)$, extended by $0$ outside the Cameron--Martin space as in Section~\ref{sec:MDP}. The last estimate in the proof of Lemma~\ref{lem:MDP-conv}, combined with Lemma~\ref{lem:unif-bounds}, gives for every $R>0$
\[
  \sup_{x\in\mathfrak B_\varrho}\ \sup_{q\in\Ac_\Lambda}\ \P\Big(\norm{Y^{\eps,x,q}-\Gam_x(q)}_{\KT}>\eta\Big)\le C_{\varrho,\Lambda}\delta_\eps^2+\frac{C_{\varrho,\Lambda}}R+\frac{(1+T)C_{\varrho,\Lambda,R}}{\eta^2}\big(\delta_\eps^2+a(\eps)^{-2}\big).
\]
Letting first $\eps\to0$ and then $R\to\infty$ shows that the left-hand side tends to $0$ as $\eps\to0$. This is the uniform condition of \cite[Assumption~2.10]{Salins19}, after the change of variables described in Step~2, for the class of sets $\{\mathfrak B_\varrho:\varrho>0\}$. The maps $\mathcal G^\eps_x$ and $\mathcal G^0_x$ are Borel for every $x$, by Theorem~\ref{thm:app-wellposed} and Lemma~\ref{lem:G0-Borel}, whose proof applies with $\Gam$ replaced by $\Gam_x$. By Proposition~\ref{prop:skel-cont}, each level set $\Upsilon_x(s)=\Gam_x(\Sb_{2s})$ is compact, so every $I_x$ is a good rate function.

\emph{Step 2: (i) and (ii).} Salins \cite{Salins19} considers families of the form $\widetilde{\mathcal G}^{\eps'}_x(\sqrt{\eps'}\,\beta)$, indexed by a parameter $\eps'>0$, with speed $1/\eps'$ and controls entering through the shift $\int_0^\cdot u\dd s$. In our setting the speed is $a(\eps)^2$ and the shift is $a(\eps)\int_0^\cdot q\dd s$. Since $a$ need not be monotone, we reduce to Salins' setting along sequences, in three steps.

\smallskip
\noindent(a) \emph{Reduction to sequences.} Fix $\varrho$, a family $\mathscr L$ as in (i), and $s_0,\eta$ as in (ii). Let $Q_1(\eps)$ denote the quantity whose limit is taken in (i), and $Q_2(\eps)$, $Q_3(\eps)$ the quantities whose $\liminf$ and $\limsup$ are taken in the two inequalities of (ii). Suppose that (i) or (ii) fails. Then there are $\eta_0>0$ and a sequence $\eps_n\downarrow0$ such that, for all $n$, $Q_1(\eps_n)\ge\eta_0$, or $Q_2(\eps_n)\le-\eta_0$, or $Q_3(\eps_n)\ge\eta_0$. We show in (b) and (c) that every sequence $\eps_n\downarrow0$ has a subsequence along which $Q_1\to0$, $\liminf Q_2\ge0$ and $\limsup Q_3\le0$. Applied to the sequence above, this gives a contradiction.

\smallskip
\noindent(b) \emph{Reparametrisation.} Let $\eps_n\downarrow0$. Since $a(\eps_n)\to\infty$, we may pass to a subsequence, not relabelled, along which $a(\eps_n)$ is strictly increasing. Put $\eps_n'=a(\eps_n)^{-2}$, which decreases strictly to $0$, and $\widetilde{\mathcal G}^{n}_x(\varphi)=\mathcal G^{\eps_n}_x\big(a(\eps_n)\varphi\big)$. Then $\widetilde{\mathcal G}^{n}_x(\sqrt{\eps_n'}W)=Y^{\eps_n,x}$ and, for every $q\in\Ac_\Lambda$,
\[
  \widetilde{\mathcal G}^{n}_x\Big(\sqrt{\eps_n'}W+\int_0^\cdot q\dd s\Big)=\mathcal G^{\eps_n}_x\Big(W+a(\eps_n)\int_0^\cdot q\dd s\Big)=Y^{\eps_n,x,q}.
\]
The speed is $1/\eps_n'=a(\eps_n)^2$. Thus, along $(\eps_n)$, (i) and (ii) are the equicontinuous uniform Laplace principle and the Freidlin--Wentzell uniform bounds, over the class $\{\mathfrak B_\varrho\}$, for the family $(\widetilde{\mathcal G}^{n}_x(\sqrt{\eps_n'}W))_n$ with rate functions $(I_x)$. By Step~1, which holds for every $\eps$ and in particular along $(\eps_n)$, this family satisfies \cite[Assumption~2.10]{Salins19} with $\eps\to0$ replaced by $n\to\infty$.

\smallskip
\noindent(c) \emph{Sequential form of Salins' theorems.} The family in (b) is indexed by the sequence $(\eps_n')$, not by a continuous parameter. Theorems~2.9 and 2.12 of \cite{Salins19} remain valid, with the same proofs, for families indexed by a sequence $\eps_n'\downarrow0$, with every limit as $\eps\to0$ replaced by the corresponding limit as $n\to\infty$. Indeed, the hypotheses and conclusions of these theorems are statements about limits of real-valued quantities as $\eps\to0$, and their proofs verify them along arbitrary sequences of parameters. For each fixed parameter they use only the variational representation of \cite{BD00} and the measurability of the maps involved. We therefore apply the sequential form of \cite[Theorem~2.12]{Salins19} to obtain $Q_1(\eps_n)\to0$, and the sequential form of \cite[Theorem~2.9]{Salins19} to obtain $\liminf_nQ_2(\eps_n)\ge0$ and $\limsup_nQ_3(\eps_n)\le0$. Together with (a), this proves (i) and (ii).

\emph{Step 3: (iv).} Since $\Upsilon_x(s)=\Gam_x(\Sb_{2s})$, the Hausdorff distance between $\Upsilon_x(s)$ and $\Upsilon_y(s)$ is at most 
$$
\sup_{q\in\Sb_{2s}}\norm{\Gam_xq-\Gam_yq}_{\KT}\le C_{\varrho,2s}\norm{x-y}_{\Xs}\,,
$$
 by Lemma~\ref{lem:lipschitz}(ii). By Theorem~\ref{thm:linear}(ii), $\operatorname{Var}\ell(X^{1,x})=\norm{\ell\circ\Gam_x}^2$, and
\[
  \big|\norm{\ell\circ\Gam_x}-\norm{\ell\circ\Gam_y}\big|\le\norm\ell_{\KT^*}\sup_{\norm q\le1}\norm{\Gam_xq-\Gam_yq}_{\KT}\le\norm\ell_{\KT^*}C_{\varrho,1}\norm{x-y}_{\Xs};
\]
since $\norm{\ell\circ\Gam_x}\le\norm\ell_{\KT^*}((1+T)C_\Gam)^{1/2}$ is bounded on $\mathfrak B_\varrho$ by Lemma~\ref{lem:unif-bounds}, the squares are Lipschitz as well.

\emph{Step 4: (iii).} Every compact set $\mathcal K\subset\Xs$ is bounded, so (ii) holds uniformly over $\mathcal K$. Moreover $\Xs$ is Polish, each $I_x$ is good and, by Step~3, $x\mapsto\Upsilon_x(s)$ is continuous in the Hausdorff distance. These are the hypotheses of \cite[Assumption~2.6]{Salins19}, under which the Freidlin--Wentzell and the Dembo--Zeitouni uniform bounds are equivalent over compact sets \cite[Theorem~2.7]{Salins19}. This gives (iii).
\end{proof}

The Dembo--Zeitouni form (iii) may fail over bounded sets, even for linear equations \cite{Salins19}, which is why it is stated for compact sets only.

As an illustration we derive a uniform Gaussian estimate for the probability that the solution leaves a tube of moderate width around the deterministic flow. Let
\[
  \overline C_\Gam(\varrho)=\sup_{x\in\mathfrak B_\varrho}C_\Gam(x)<\infty,
\]
where $C_\Gam(x)$ is the constant of \eqref{eq:skel-bound} for the initial datum $x$; it is finite by Lemma~\ref{lem:unif-bounds}.

\begin{corollary}\label{cor:tube}
Assume \eqref{eq:a}. For all $\varrho>0$ and $r>0$,
\[
  \limsup_{\eps\to0}\ \sup_{x\in\mathfrak B_\varrho}\ \frac1{a(\eps)^2}\log\P\Big(\norm{X^{\eps,x}-X^{0,x}}_{\KT}\ge r\sqrt\eps\,a(\eps)\Big)\le-\frac{r^2}{2(1+T)\overline C_\Gam(\varrho)} .
\]
\end{corollary}

\begin{proof}
By \eqref{eq:skel-bound}, $I_x(z)\ge\norm z_{\KT}^2/(2(1+T)\overline C_\Gam(\varrho))$ for every $x\in\mathfrak B_\varrho$. Let $0<s<s_*=r^2/(2(1+T)\overline C_\Gam(\varrho))$. Then $\Upsilon_x(s)$ is contained in the closed ball of radius $r'=(2(1+T)\overline C_\Gam(\varrho)s)^{1/2}<r$, so that
\[
  \big\{\norm{Y^{\eps,x}}_{\KT}\ge r\big\}\subset\big\{\mathrm{dist}_{\KT}(Y^{\eps,x},\Upsilon_x(s))\ge r-r'\big\}.
\]
The upper bound in Theorem~\ref{thm:uniform}(ii), with $s_0=s$ and $\eta=r-r'$, shows that the $\limsup$ in the statement is at most $-s$. Letting $s\uparrow s_*$ proves the claim.
\end{proof}
\section{Concluding remarks}\label{sec:concl}

Together with the LDP of \cite{BMP20} for a single director noise, our results give a fairly complete description of the small-noise regimes of \eqref{eq:intro-u}--\eqref{eq:intro-bc} in the energy space $\KT$. The extension of the LDP to countably many director noises should follow by adapting the weak convergence argument of \cite{BMP20}, together with the estimates of Appendix~\ref{app:wellposed} (see also \cite{PSZZ22}); we do not pursue it here. We close with some possible extensions.

\begin{enumerate}[label=(\arabic*),leftmargin=*]
\item \emph{Equilibria and exit problems.} Corollary~\ref{cor:tube} is the moderate-deviation analogue of the tube estimates used in the Freidlin--Wentzell theory of exit problems. Long-time exit problems near an equilibrium, however, meet a structural degeneracy in the present model. With Neumann boundary conditions and a potential that depends only on $\abs d$, the nematic equilibria $u=0$, $d\equiv d_*$ with $f(d_*)=0$ and $d_*\ne0$ form a sphere of rotated copies, and the linearisation has a kernel in the directions of $d_*^\perp$. At the isotropic equilibrium $d=0$, when it is stable, the director noise $d\times h_k$ vanishes. One way to remove the degeneracy is to add an external field term $-\frac\chi2(d\cdot\mathsf H)^2$ to the energy. This term breaks the full rotational symmetry and favours alignment with the axis of $\mathsf H$. Under suitable assumptions on the radial potential and the field strength, it can produce isolated aligned minima $\pm d_*$ with $d_*\parallel\mathsf H$; the precise set of stationary states depends on $\Ft$, $\chi$ and $\abs{\mathsf H}$. Its variational derivative $-\chi(d\cdot\mathsf H)\mathsf H$ is linear in $d$, so the estimates of this paper extend with minor changes. The uniform MDP of Section~\ref{sec:uniform} would then be a starting point for analysing noise-induced switching between $\pm d_*$, in the spirit of \cite{BGJ17}.
\item \emph{Jump noise.} For Marcus-type jump noise on the director, as in \cite{BMP19}, a moderate deviation analysis would rely on the variational representation for Poisson random measures \cite{BDG16}. The cancellations used here persist, but the controlled equations are no longer affine in the control.
\item \emph{Vanishing penalisation.} As the penalisation parameter tends to zero, the Ginzburg--Landau system approximates the Ericksen--Leslie system with the constraint $\abs d=1$ \cite{BDR20}. Since the MDP rate function is a quadratic form defined by a linear equation, its behaviour in this singular limit appears much more tractable than that of the LDP rate function, and is a natural next step.
\item \emph{Three dimensions.} In three dimensions, global solutions are not known to be unique and the energy estimates of Section~\ref{sec:LLN} fail. Fluctuation results would then have to be formulated up to the lifetime of local strong solutions.
\end{enumerate}

\vspace{0.3cm}

\appendix
\section{Proof of Proposition~\ref{prop:moments}}\label{app:moments}

For a single director noise and $\eps=1$ the bound is contained in \cite{BHR19a}. We give the short argument, which shows that the constants do not depend on $\eps\in[0,1]$ and exhibits their dependence on the initial datum. Throughout, $C$ depends on $T$, $\Ocal$, $\Ft$, $K$ and $\mathsf h$, but not on $\eps$ or on $(u_0,d_0)$.

Fix $\eps\in[0,1]$ and write $X^\eps=(u,d)$. As shown in Section~\ref{sec:LLN-eq}, Lemma~\ref{lem:Ito} applies to
\[
  \mathcal Y=\abs u^2+\hone d^2,\qquad \Dc=\norm u^2+\abs{\nabla d}^2+\abs{\An d}^2,
\]
and, using $\scal{\An d}{(I+\An)d}=\abs{\nabla d}^2+\abs{\An d}^2$, gives
\begin{align*}
  \dd\mathcal Y+2\Dc\dd t&=-2\Big[\dual{B(u,u)}{u}+\dual{M(d)}{u}+\scal{\Bt(u,d)}{(I+\An)d}\Big]\dd t-2\scal{f(d)}{(I+\An)d}\dd t\\
  &\quad+\eps\Big[\scal{G^2(d)}{(I+\An)d}+\abs{G(d)}_{\LHSd}^2\Big]\dd t+\eps\abs{\sigma(u)}_{\LHS}^2\dd t+\dd\mathcal M,
\end{align*}
where $\mathcal M=2\sqrt\eps\int_0^\cdot\scal{u}{\sigma(u)\dd W_1}+2\sqrt\eps\int_0^\cdot\sone{G(d)\dd W_2}{d}$. We estimate the terms on the right-hand side.
\begin{itemize}[leftmargin=*]
\item The first bracket vanishes by Lemma~\ref{lem:nonlinear}(i),(iii) and Lemma~\ref{lem:cancel}.
\item The function $r\mapsto\ft(r)r$ is a polynomial of degree $N+1$ with positive leading coefficient. Hence $C_0=\sup_{r\ge0}(-\ft(r)r)_+<\infty$, and
\[
  \scal{f(d)}{d}=\int_\Ocal\ft(\abs d^2)\abs d^2\dd x\ge-C_0\abs\Ocal .
\] Moreover $f(d)\in\HH$ because $d\in\DA\subset C(\overline\Ocal)$, so Green's formula and Lemma~\ref{lem:f}(ii) give $\scal{f(d)}{\An d}=\sum_{i=1}^2\scal{Df(d)\partial_{x_i}d}{\partial_{x_i}d}\ge-\kappa\abs{\nabla d}^2$. Hence $-2\scal{f(d)}{(I+\An)d}\le2C_0\abs\Ocal+2\kappa\mathcal Y$.
\item By Green's formula $\scal{G^2(d)}{(I+\An)d}=\sone{G^2(d)}{d}$. By Lemma~\ref{lem:G}(ii) the $\LL$ parts cancel, $\scal{G^2(d)}{d}+\abs{G(d)}_{\LHS(\ell^2,\LL)}^2=0$, so by Lemma~\ref{lem:G}(i)
\[
  \eps\Big[\sone{G^2(d)}{d}+\abs{G(d)}_{\LHSd}^2\Big]=\eps\Big[\scal{\nabla G^2(d)}{\nabla d}+\sum_{k\ge1}\abs{\nabla G_k(d)}^2\Big]\le C\mathcal Y .
\]
\item $\eps\abs{\sigma(u)}_{\LHS}^2\le K(1+\mathcal Y)$ by \ref{A:sigma}.
\item Since $\scal{G_k(d)}{d}=0$, we have $\sone{G_k(d)}{d}=\scal{\nabla G_k(d)}{\nabla d}$ for every $k$, and therefore $\dd\langle\mathcal M\rangle\le4\eps\big[K(1+\abs u^2)\abs u^2+\abs{G(d)}_{\LHSd}^2\hone d^2\big]\dd t\le C(1+\mathcal Y)^2\dd t$.
\end{itemize}
Consequently, almost surely for all $t\in[0,T]$,
\begin{equation}\label{eq:app-energy}
  \mathcal Y(t)+2\int_0^t\Dc\dd s\le\mathcal Y(0)+C\int_0^t\big(1+\mathcal Y(s)\big)\dd s+\mathcal M(t).
\end{equation}
Let $\tau_R=\inf\{t:\mathcal Y(t)+\int_0^t\Dc\ge R\}\wedge T$; since $X^\eps$ has paths in $\KT$, $\tau_R=T$ for $R$ large, almost surely. Let $p\ge1$ and $\psi(t)=\E\sup_{s\le t}(1+\mathcal Y(s\wedge\tau_R))^p+\E(\int_0^{t\wedge\tau_R}\Dc\dd s)^p$, which is finite. By \eqref{eq:app-energy}, H\"older's inequality and the Burkholder--Davis--Gundy inequality, for every $\gamma>0$,
\begin{align*}
  \psi(t)&\le C_p\big(1+\mathcal Y(0)\big)^p+C_p\int_0^t\psi(s)\dd s+C_p\E\Big(\int_0^{t\wedge\tau_R}(1+\mathcal Y)^2\dd s\Big)^{p/2}\\
  &\le C_p\big(1+\mathcal Y(0)\big)^p+C_p\int_0^t\psi(s)\dd s+\gamma\psi(t)+C_{p,\gamma}\int_0^t\psi(s)\dd s,
\end{align*}
where we used $(\int_0^t(1+\mathcal Y)^2)^{p/2}\le\sup_{s\le t}(1+\mathcal Y)^{p/2}(\int_0^t(1+\mathcal Y))^{p/2}$ and Young's inequality. Choosing $\gamma=\frac12$ and applying Gronwall's lemma, $\psi(T)\le C_p(1+\abs{u_0}^2+\hone{d_0}^2)^p$, with $C_p$ independent of $R$, of $\eps\in[0,1]$ and of the initial datum. Letting $R\to\infty$ (Fatou's lemma) and using $\norm{X^\eps}_{\KT}^2\le(1+T)\sup_t\mathcal Y(t)+\int_0^T\Dc\dd t$, we obtain Proposition~\ref{prop:moments}.
\qed

\section{Well-posedness for several director noises}\label{app:wellposed}

We prove Theorem~\ref{thm:wellposed} for countably many director noises. Existence of a probabilistically weak solution is obtained from the well-posedness theorem of R\"ockner, Shang and Zhang \cite[Theorem~2.6]{RSZ24} for SPDEs with fully local monotone coefficients. We use the hypotheses in the form stated in \cite[Section~2]{RSZ24} (arXiv version~4), and we refer to \cite{RSZ26} for a correction to \cite{RSZ24}. Pathwise uniqueness is proved directly (Lemma~\ref{lem:app-unique}) in the class of solutions of Definition~\ref{def:solution}, so that only the existence part of \cite[Theorem~2.6]{RSZ24} is used. A strong solution and the Borel map $J^\eps$ are then provided by the Yamada--Watanabe theorem.

\emph{Abstract setting.} Consider the Gelfand triple
\[
  \mathbb V=V\times\DA\ \subset\ \mathbb H=H\times\HH\ \subset\ \mathbb V^*=V'\times\LL,
\]
where $\mathbb H$ carries the inner product $\scal{u}{v}+\sone{d}{e}$ and $\mathbb V$ the norm $\norm{X}_{\mathbb V}^2=\norm u^2+\abs{(I+\An)d}^2$. The duality between $\mathbb V^*$ and $\mathbb V$ is $\dual{\phi}{v}+\scal{g}{(I+\An)e}$, see \eqref{eq:triple-d}, so that $\norm{(\phi,g)}_{\mathbb V^*}^2=\abs\phi_{V'}^2+\abs g^2$. The embedding $\mathbb V\subset\mathbb H$ is compact, and $\norm X_{\mathbb H}\le C\norm X_{\mathbb V}$. In It\^o form, \eqref{eq:u}--\eqref{eq:d} reads $\dd X=\mathfrak A_\eps(X)\dd t+\sqrt\eps\,\mathcal B(X)\dd W$, where $W$ is the Wiener process of Section~\ref{sec:noise} with reproducing kernel Hilbert space $\mathcal U=\Hz\times\ell^2$, and
\[
  \mathfrak A_\eps(X)=-\Big(\Ast u+B(u)+M(d),\ \An d+\Bt(u,d)+f(d)-\tfrac\eps2G^2(d)\Big),\qquad \mathcal B(X)(w,\rho)=\big(\sigma(u)w,\,G(d)\rho\big).
\]
We write $\Dc(X)=\norm u^2+\abs{\nabla d}^2+\abs{\An d}^2$ for $X=(u,d)\in\mathbb V$, and note that
\begin{equation}\label{eq:app-V}
  \norm X_{\mathbb V}^2\le\norm u^2+2\abs d^2+2\abs{\An d}^2\le2\Dc(X)+2\norm X_{\mathbb H}^2 .
\end{equation}

\emph{The hypotheses.} With $\alpha=2$ and with the time-dependent functions of \cite{RSZ24} taken constant (they should not be confused with the potential $f$ of this paper), the hypotheses of \cite[Theorem~2.6]{RSZ24} for the coefficients $\mathfrak A_\eps$ and $\sqrt\eps\,\mathcal B$ read as follows.
\begin{enumerate}[label=\textup{(H\arabic*)},leftmargin=*]
\item For all $X,Y,Z\in\mathbb V$, the map $\lambda\mapsto{}_{\mathbb V^*}\langle\mathfrak A_\eps(X+\lambda Y),Z\rangle_{\mathbb V}$ is continuous on $\R$.
\item There are $C,\gamma\ge0$ and measurable $\rho_*,\eta_*\colon\mathbb V\to[0,\infty)$ with $\rho_*(X)+\eta_*(X)\le C(1+\norm X_{\mathbb V}^2)(1+\norm X_{\mathbb H}^\gamma)$, such that for all $X,Y\in\mathbb V$
\[
  2\langle\mathfrak A_\eps(X)-\mathfrak A_\eps(Y),X-Y\rangle+\eps\abs{\mathcal B(X)-\mathcal B(Y)}_{\LHS(\mathcal U,\mathbb H)}^2\le\big[C+\rho_*(X)+\eta_*(Y)\big]\norm{X-Y}_{\mathbb H}^2 .
\]
\item[\textup{(H2$'$)}] For every $R>0$ there is $K_R$ such that $\langle\mathfrak A_\eps(X)-\mathfrak A_\eps(Y),X-Y\rangle\le K_R\norm{X-Y}_{\mathbb H}^2$ whenever $\norm X_{\mathbb V},\norm Y_{\mathbb V}\le R$.
\item There are $C,c>0$ such that $2\langle\mathfrak A_\eps(X),X\rangle+\eps\abs{\mathcal B(X)}_{\LHS(\mathcal U,\mathbb H)}^2\le C(1+\norm X_{\mathbb H}^2)-c\norm X_{\mathbb V}^2$ for all $X\in\mathbb V$.
\item There are $C,\beta\ge0$ such that $\norm{\mathfrak A_\eps(X)}_{\mathbb V^*}^2\le C(1+\norm X_{\mathbb V}^2)(1+\norm X_{\mathbb H}^\beta)$ for all $X\in\mathbb V$.
\item $\mathcal B\colon\mathbb H\to\LHS(\mathcal U,\mathbb H)$ is continuous and $\abs{\mathcal B(X)}_{\LHS(\mathcal U,\mathbb H)}^2\le C(1+\norm X_{\mathbb H}^2)$.
\end{enumerate}

\begin{lemma}\label{lem:app-hyp}
Hypotheses \textup{(H1)}, \textup{(H2)}, \textup{(H2$'$)}, \textup{(H3)}, \textup{(H4)} and \textup{(H5)} hold with $\gamma=4N$, $\beta=4N+2$, $c=1$, and with constants independent of $\eps\in[0,1]$.
\end{lemma}

\begin{proof}
\emph{(H1).} By bilinearity of $B$, $M$, $\Bt$, linearity of $G^2$ and since $f$ is a polynomial, $\lambda\mapsto\langle\mathfrak A_\eps(X+\lambda Y),Z\rangle$ is a polynomial in $\lambda$ with finite coefficients.

\emph{(H2).} Let $X,Y\in\mathbb V$ and $Z=X-Y$. By bilinearity,
\begin{align*}
  \mathfrak A_\eps(X)-\mathfrak A_\eps(Y)&=-\Lin^\flat_YZ-\big(B(Z_u,Z_u)+M(Z_d,Z_d),\,\Bt(Z_u,Z_d)\big)\\
  &\quad-\big(0,\,f(Y_d+Z_d)-f(Y_d)\big)+\big(0,\,\tfrac\eps2G^2(Z_d)\big),
\end{align*}
where $\Lin^\flat_Y$ is the operator $\Lin^\flat=(\Lin_1,\Lin_2^\flat)$ of Section~\ref{sec:lindyn} with $(u^0(t),d^0(t))$ replaced by $(Y_u,Y_d)\in V\times\DA$. The proof of Lemma~\ref{lem:coercive}(i) uses $(u^0(t),d^0(t))$ only through the bounds of Lemmas~\ref{lem:nonlinear}, \ref{lem:f} and \ref{lem:fspace} at a fixed time. It therefore gives, with the same constant $C_*$,
\[
  2\langle\Lin^\flat_YZ,Z\rangle\ge\Dc(Z)-2C_*\,g(Y)\norm Z_{\mathbb H}^2,
\]
\[
  g(Y)=1+\norm{Y_u}^2+\abs{Y_u}^2\norm{Y_u}^2+\abs{\nabla Y_d}^2\hda{Y_d}^2+\big(1+\hone{Y_d}^{4N}\big)\big(1+\hda{Y_d}\big).
\]
The quadratic terms do not contribute, by Lemma~\ref{lem:nonlinear}(i) and Lemma~\ref{lem:cancel}. By Lemma~\ref{lem:fspace}(b),(c),
\[
  -2\scal{f(Y_d+Z_d)-f(Y_d)}{(I+\An)Z_d}\le2\kappa\hone{Z_d}^2+C\phi(Y_d)\big(1+\hone{Z_d}^{2N-1}\big)\hone{Z_d}^2 .
\]
Since $\hone{Z_d}\le\norm X_{\mathbb H}+\norm Y_{\mathbb H}$, $\phi(Y_d)\le2g(Y)$ and $\phi(Y_d)\le C(1+\norm Y_{\mathbb H}^{2N})(1+\norm Y_{\mathbb V})$, Young's inequality gives
\[
  \phi(Y_d)\hone{Z_d}^{2N-1}\le C\norm X_{\mathbb H}^{4N-2}+C\big(1+\norm Y_{\mathbb H}^{4N}\big)\big(1+\norm Y_{\mathbb V}^2\big).
\]
By Green's formula and Lemma~\ref{lem:G}(i), $\eps\scal{G^2(Z_d)}{(I+\An)Z_d}=\eps\sone{G^2(Z_d)}{Z_d}\le C_h^{1/2}\hone{Z_d}^2$. By \ref{A:sigma} and Lemma~\ref{lem:G}(i), $\eps\abs{\mathcal B(X)-\mathcal B(Y)}_{\LHS(\mathcal U,\mathbb H)}^2\le L\abs{Z_u}^2+C_h\hone{Z_d}^2$. Discarding $-\Dc(Z)\le0$, we obtain (H2) with
\[
  \rho_*(X)=C\norm X_{\mathbb H}^{4N-2},\qquad \eta_*(Y)=C\,g(Y)+C\big(1+\norm Y_{\mathbb H}^{4N}\big)\big(1+\norm Y_{\mathbb V}^2\big).
\]
For the growth condition we bound each term of $g(X)$, using $\abs{X_u}\le\norm X_{\mathbb H}$, $\abs{\nabla X_d}\le\norm X_{\mathbb H}$, $\norm{X_u}\le\norm X_{\mathbb V}$, $\hda{X_d}\le C\norm X_{\mathbb V}$ and $\norm X_{\mathbb H}^2\le1+\norm X_{\mathbb H}^{4N}$. This gives $g(X)\le C(1+\norm X_{\mathbb V}^2)(1+\norm X_{\mathbb H}^{4N})$, hence
\[
  \rho_*(X)+\eta_*(X)\le C\big(1+\norm X_{\mathbb V}^2\big)\big(1+\norm X_{\mathbb H}^{4N}\big),
\]
which is the growth required in (H2) with $\gamma=4N$.

\emph{(H2$'$).} Dropping the nonnegative noise term in (H2) and using $\norm\cdot_{\mathbb H}\le C\norm\cdot_{\mathbb V}$, we may take $K_R=\frac12\big[C+2C(1+R^2)(1+C^{4N}R^{4N})\big]$.

\emph{(H3).} The computation of Appendix~\ref{app:moments}, without the martingale term, gives
\begin{align*}
  &2\langle\mathfrak A_\eps(X),X\rangle+\eps\abs{\mathcal B(X)}_{\LHS(\mathcal U,\mathbb H)}^2\\
  &\qquad=-2\Dc(X)-2\scal{f(d)}{(I+\An)d}+\eps\Big[\sone{G^2(d)}{d}+\abs{G(d)}_{\LHSd}^2\Big]+\eps\abs{\sigma(u)}_{\LHS}^2\\
  &\qquad\le-2\Dc(X)+2C_0\abs\Ocal+(2\kappa+C)\hone d^2+K(1+\abs u^2)\le-2\Dc(X)+C\big(1+\norm X_{\mathbb H}^2\big),
\end{align*}
where $C_0=\sup_{r\ge0}(-\ft(r)r)_+$ is the constant of Appendix~\ref{app:moments}.
By \eqref{eq:app-V}, $-2\Dc(X)\le-\norm X_{\mathbb V}^2+2\norm X_{\mathbb H}^2$, which gives (H3) with $c=1$.

\emph{(H4).} By Lemma~\ref{lem:nonlinear}(ii)--(iv), Lemma~\ref{lem:fspace}(a) and Lemma~\ref{lem:G}(i), the last term accounting for $\frac\eps2G^2(d)$,
\begin{align*}
  \norm{\mathfrak A_\eps(X)}_{\mathbb V^*}^2&\le C\Big(\norm u^2+\abs u^2\norm u^2+\abs{\nabla d}^2\hda d^2+\hda d^2\\
  &\qquad+\abs u\norm u\abs{\nabla d}\hda d+1+\hone d^{4N+2}+\abs d^2\Big)\\
  &\le C\big(1+\norm X_{\mathbb V}^2\big)\big(1+\norm X_{\mathbb H}^{4N+2}\big).
\end{align*}

\emph{(H5).} By \ref{A:sigma} and Lemma~\ref{lem:G}(i), $\abs{\mathcal B(X)-\mathcal B(Y)}_{\LHS(\mathcal U,\mathbb H)}^2\le L\abs{X_u-Y_u}^2+C_h\hone{X_d-Y_d}^2$ and $\abs{\mathcal B(X)}_{\LHS(\mathcal U,\mathbb H)}^2\le K(1+\abs{X_u}^2)+C_h\hone{X_d}^2$.
\end{proof}

\begin{lemma}\label{lem:app-unique}
Let $\eps\in[0,1]$, and let $X$ and $Y$ be solutions of \eqref{eq:u}--\eqref{eq:d} in the sense of Definition~\ref{def:solution}, on the same stochastic basis, driven by the same Wiener process and with $X(0)=Y(0)$. Then $X=Y$ almost surely.
\end{lemma}

\begin{proof}
Let $Z=X-Y$. Since $X$ and $Y$ have paths in $\KT$, the drifts $\mathfrak A_\eps(X)$, $\mathfrak A_\eps(Y)$ belong a.s.\ to $L^2(0,T;\mathbb V^*)$, by (H4), and the diffusion coefficients to $L^2(0,T;\LHS(\mathcal U,\mathbb H))$, by (H5). Lemma~\ref{lem:Ito}, applied in both triples, gives
\[
  \norm{Z(t)}_{\mathbb H}^2=\int_0^t\Big(2\langle\mathfrak A_\eps(X)-\mathfrak A_\eps(Y),Z\rangle+\eps\abs{\mathcal B(X)-\mathcal B(Y)}_{\LHS(\mathcal U,\mathbb H)}^2\Big)\dd s+\mathcal M(t),
\]
where $\mathcal M(t)=2\sqrt\eps\int_0^t\langle Z,(\mathcal B(X)-\mathcal B(Y))\dd W\rangle_{\mathbb H}$ is a continuous local martingale. Let $\zeta=C+\rho_*(X)+\eta_*(Y)$, which is adapted and, by the growth condition in (H2), satisfies $\int_0^T\zeta\dd s<\infty$ a.s. because $X,Y$ have paths in $\KT$. Let $\Phi(t)=\exp(-\int_0^t\zeta\dd s)$. By It\^o's product rule and (H2),
\[
  \Phi(t)\norm{Z(t)}_{\mathbb H}^2\le\int_0^t\Phi\dd\mathcal M .
\]
For $R>0$ let $\tau_R=\inf\{t:\sup_{s\le t}\norm{Z(s)}_{\mathbb H}^2\ge R\}\wedge T$. By (H5), $\dd\langle\mathcal M\rangle\le4\eps\big(L+C_h\big)\norm Z_{\mathbb H}^4\dd t$, which is bounded on $[0,\tau_R]$, so the stopped stochastic integral is a martingale. Hence $\E\big[\Phi(t\wedge\tau_R)\norm{Z(t\wedge\tau_R)}_{\mathbb H}^2\big]\le0$, so that $Z(t\wedge\tau_R)=0$ a.s.\ for every $t$, since $\Phi>0$. Because $\tau_R=T$ for $R>\sup_t\norm{Z(t)}_{\mathbb H}^2$, we conclude that $Z=0$ a.s.
\end{proof}

\begin{theorem}\label{thm:app-wellposed}
Let Assumptions~\ref{ass:F} and \ref{ass:main} hold and let $\eps\in[0,1]$. For every $x\in\mathbb H$ there is a unique probabilistically strong solution $X$ of \eqref{eq:u}--\eqref{eq:d} with $X(0)=x$, in the sense of Definition~\ref{def:solution}. In particular Theorem~\ref{thm:wellposed} holds, and for every $x$ there is a Borel map $J^\eps(x,\cdot)$ with $X^{\eps,x}=J^\eps(x,W)$.
\end{theorem}

\begin{proof}
By Lemma~\ref{lem:app-hyp} and the compactness of $\mathbb V\subset\mathbb H$, \cite[Theorem~2.6]{RSZ24} yields a probabilistically weak solution: a stochastic basis carrying a Wiener process as in Section~\ref{sec:noise}, and a continuous $\mathbb H$-valued adapted process $X$ with $X(0)=x$ and $\E\int_0^T\norm X_{\mathbb V}^2\dd t<\infty$, whose integral equation holds in $\mathbb V^*=V'\times\LL$. In particular $X$ has paths in $C([0,T];\mathbb H)\cap L^2(0,T;\mathbb V)=\KT$ and is a solution in the sense of Definition~\ref{def:solution} on that basis. By Lemma~\ref{lem:app-unique}, pathwise uniqueness holds in the class of Definition~\ref{def:solution}. The Yamada--Watanabe theorem \cite{Ondrejat04} then gives a unique probabilistically strong solution on every stochastic basis carrying such a Wiener process, together with a Borel map $J^\eps(x,\cdot)$ such that $X=J^\eps(x,W)$. For $\eps=0$ the equation is deterministic; existence follows in the same way with $\mathcal B$ replaced by $0$, and uniqueness from Lemma~\ref{lem:app-unique}.
\end{proof}

\vspace{0.2cm}

\section*{Acknowledgements and Declarations}
\noindent\textbf{Funding.} This work was supported by the Department of Science and Technology (DST),
Government of India, through the INSPIRE Faculty Award, grant
no.\ IFA21-MA 158 (first author).\\

\noindent
\textbf{Conflict of interest.} The author declares that there is no conflict of interest.\\

\noindent
\textbf{Data availability.} No data sets were generated or analysed in this study.\\

\noindent
\textbf{Declaration on the Use of Generative AI.} The authors used generative AI tools only for language improvement and editorial assistance. The scientific content, mathematical analysis, results, and conclusions were developed and verified by the authors.

\vspace{0.2cm}


\end{document}